\documentclass[pdf,aos,preprint]{imsart}
\usepackage{xr}
\usepackage{amsmath, mathrsfs}
\RequirePackage[numbers]{natbib}
\usepackage{graphicx}
\usepackage{color} 
\usepackage{verbatim}
\usepackage{hyperref}
\newcommand{\qed}{{\unskip\nobreak\hfil\penalty50\hskip2em\vadjust{}
            \nobreak\hfil$\Box$\parfillskip=0pt\finalhyphendemerits=0\par}}

\newtheorem{thm}{Theorem}[section] %(If you want theorem numbered
\newtheorem{lemma}{Lemma}[section] %%    with section number.

\newtheorem{definition}{Definition}[section]
\newcommand{\bed}{\begin{definition}}
\newcommand{\eed}{\end{definition}}
\newcommand{\beq}{\begin{equation}}
\newcommand{\eeq}{\end{equation}}

\newcommand{\eps}{\epsilon}

\newcommand{\bitem}{\begin{itemize}}
\newcommand{\eitem}{\end{itemize}}

\newcommand{\goto}{\rightarrow}

\newcommand{\margmax}{\mathrm{argmax}}

\newcommand{\beqn}{\begin{equation}}
\newcommand{\eeqn}{\end{equation}}
\newcommand{\balign}{\begin{align}}
\newcommand{\ealign}{\end{align}}

\newcommand{\sgn}{\mathrm{sgn}}

\newcommand{\hamm}{\mathrm{Hamm}}

\newcommand{\bel}{\begin{eqnarray}\label}
\newcommand{\eel}{\end{eqnarray}}
\newcommand{\bes}{\begin{eqnarray*}}
\newcommand{\ees}{\end{eqnarray*}}
\newcommand{\hell}{\hat{\ell}}

\def\sgn{\hbox{\rm sgn}}

\newcommand{\muhat}{\hat{\mu}} 
\newcommand{\hs}{\hat{S}} 
\newcommand{\ang}{\mathrm{cos}}

\usepackage{amssymb}
\begin{document}
\begin{frontmatter}

\title{Phase Transitions for High Dimensional Clustering and Related Problems}

\begin{aug}
\author{\fnms{Jiashun} \snm{Jin}\thanksref{t1}\ead[label=e1]{jiashun@stat.cmu.edu}},
\author{\fnms{Zheng Tracy} \snm{Ke}\thanksref{t2}\ead[label=e2]{zke@galton.uchicago.edu}},
\and
\author{\fnms{Wanjie} \snm{Wang}\thanksref{t3}\ead[label=e3]{wanjiew@wharton.upenn.edu}}

%\thankstext{t1}{}
%\thankstext{t2}{Supported in part by NSF Grant DMS-1208315.}
%\thankstext{t3}{}

\affiliation{Carnegie Mellon University\thanksmark{t1},  
University of Chicago\thanksmark{t2} and University of 
Pennsylvania\thanksmark{t3}} 
    
         \address{J. Jin\\
		Department of Statistics\\
		Carnegie Mellon University\\
		Pittsburgh, Pennsylvania, 15213\\
		USA\\
		\printead{e1}}
	
	\address{Z. Ke\\
		Department of Statistics\\
		University of Chicago\\
		Chicago, Illinois, 60637\\
		USA\\
		\printead{e2}}

	\address{W. Wang\\
		Department of Biostatistics and Epidemiology\\
		University of Pennsylvania, Perelman School of Medicine \\
		Philadelphia, Pennsylvania, 19104\\
		USA\\
		\printead{e3}}
\end{aug}

%%%%%%%%%%%%
%%%%%%%%%%%%
%%%%%%%%%%%%
\begin{abstract} 
Consider a two-class clustering problem where we observe  $X_i  = \ell_i \mu + Z_i$,   $Z_i \stackrel{iid}{\sim} N(0, I_p)$, $1 \leq i \leq n$.  
The feature vector $\mu\in R^p$ is   unknown but  is presumably sparse.    The class labels $\ell_i\in\{-1, 1\}$   are  also unknown and the main interest is to estimate them. 

We are interested in the statistical limits. 
In the two-dimensional phase space calibrating the rarity and strengths of useful features,   we find the precise demarcation for the {\it Region of Impossibility} and   
{\it Region of Possibility}. 
In the former, useful features are too rare/weak for  successful clustering. 
In the latter, useful features are strong enough to allow  successful clustering. 
The results are extended to the  case of colored noise using Le Cam's idea 
on comparison of experiments. 

We also extend the study on statistical limits for clustering to that for signal recovery and 
that for global testing. We compare the statistical limits for three problems and  expose some interesting insight.

We propose classical PCA and Important Features PCA (IF-PCA) for clustering. 
For a threshold $t > 0$, IF-PCA clusters by applying classical PCA to all columns of 
$X$ with an $L^2$-norm  larger than $t$. 
We also propose  two aggregation methods.   
For any parameter in the Region of Possibility, some  of these  methods 
yield successful clustering.

We discover a phase transition for IF-PCA. 
For any threshold $t > 0$, let $\xi^{(t)}$ be the first  left singular vector of the post-selection data matrix.
The phase space partitions into two different regions. 
In one region, there is a $t$ such that  $\mathrm{cos}(\xi^{(t)}, \ell)  \goto 1$ and IF-PCA yields successful clustering. In the other, 
$\mathrm{cos}(\xi^{(t)}, \ell)  \leq c_0 < 1$ for all $t > 0$.

Our results require delicate analysis, especially on post-selection Random Matrix Theory and on  lower bound arguments. 
%%%%%%%%%%%
%%%%%%%%%%%
%%%%%%%%%%%
\end{abstract}

 

\end{frontmatter}

\section{Introduction}
\setcounter{equation}{0} 
Motivated by the interest on gene microarray study,  
we consider a clustering problem where we have $n$ subjects from two 
different classes (e.g., normal and diseased), measured on the same set of 
$p$ features (i.e., gene expression level).  
To facilitate the analysis, we assume that  two classes are equally likely so the class labels satisfy 
\begin{equation} \label{model1} 
\ell_i \stackrel{iid}{\sim} 2 \mathrm{Bernoulli}(1/2) -1, \qquad 1\leq i\leq n.  
\end{equation} 
We also assume that  the $p$-dimensional data vectors $X_i$'s are standardized, so that for a  contrast mean vector $\mu \in R^p$, 
%%%%%%%%%
%%%%%%%%%
%%%%%%%%%
%%%%%%%%%
\begin{equation} \label{model2} 
X_i = \ell_i  \mu + Z_i,   \qquad Z_i \stackrel{iid}{\sim} N(0, I_p), \qquad 1 \leq i \leq n. 
\end{equation}  
Throughout this paper,  we call feature $j$, $1 \leq j \leq p$,  a ``useless feature" or ``noise" if $\mu(j) = 0$  
and a ``useful feature" or ``signal" otherwise.

The paper focuses on the problem of clustering (i.e., estimating the class labels $\ell_i$).  Such 
a problem is of interest, especially in the study of complex disease \cite{Lee}. 
In the two-dimensional phase space calibrating the signal rarity and signal strengths,   
we   are interested in  the following limits. \footnote{All limits in this paper are with respect to the ARW model introduced in Section \ref{subsec:ARW}.} 
\begin{itemize} 
\item {\it Statistical limits}. This is the precise boundary that separates the Region of Impossibility and Region of 
Possibility. In the former, the signals are so rare and (individually) weak that it is impossible for any method to correctly identify most of the class labels. In the latter, the signals are strong enough to allow successful clustering, and it is desirable to develop methods that  cluster successfully. 
\item {\it Computationally tractable statistical limits}.  This is similar to the boundary above,  except  that for both Possibility and Impossibility,  we only consider statistical methods that are computationally tractable.  
\end{itemize}  
We use Region of Possibility and Region of Impossibility as generic terms, which may vary from occurrence to occurrence.  %In this paper, we precisely characterize the statistical limits, and establish upper bounds   
%for the computationally tractable statistical limits. 

The paper also contains three closely related objectives as follows, which we discuss in Sections~\ref{subsec:IF-PCA} and \ref{sec:IF-PCA}, Section~\ref{sec:sig}, and Section~\ref{sec:hyp}, respectively. 
\begin{itemize} 
\item Performance of the recent idea of Important Features PCA (IF-PCA).   
\item Limits for recovering the support of $\mu$ (signal recovery).    
\item Limits for testing whether $X_i$'s are {\it iid} samples from  $N(0, I_p)$,  or generated from Model (\ref{model2}) (hypothesis testing). 
\end{itemize}  
Our work on sparse clustering is related to  Azizyan {\it et al} \cite{WassermanSparse} and Chan and Hall \cite{CH} (see also \cite{PanShen2007, raftery2006,Sun2012,skmeans}):  the three papers share the same spirit that we should do a  feature selection  before we cluster.  
Our work on support recovery is related to recent interest on sparse PCA  (e.g.,  Amini and Wainwright \cite{AW},  Johnstone and Lu \cite{JohnstoneLu},  Vu and  Lei \cite{Vu}, Wang {\it et al} \cite{WLL},  Arias-Castron and Verzelen \cite{Ery}),  and our work on hypothesis testing is related to recent interest  on matrix estimation and matrix testing  (e.g.,  Arias-Castro and Verzelen \cite{Ery}, Cai {\it et al} \cite{CMW}). 
However, our work is different in many important aspects, especially for our focus on the limits and on the Rare/Weak models.  See Section \ref{sec:discuss}  for more discussion.

%%%%%%%%%%%
%%%%%%%%%%%
%%%%%%%%%%%
\subsection{Four clustering methods}  \label{subsec:4method}  
Denoting the data matrix by $X$, we write 
\[
X' = [X_1, X_2, \ldots, X_n], \qquad  X = [x_1, x_2, \ldots, x_p]. 
\]
We introduce two  methods: a feature aggregation method and IF-PCA.  Each method includes a special  case, which 
can be viewed as a different method.

The first method $\hell^{(sa)}_N$ targets on the case where the signals are rare but  individually  strong (``sa": {\it Sparse Aggression}; $N$:  tuning parameter; usually,  $N \ll p$), so feature selection is desirable.   Denote the support of $\mu$ by   
\begin{equation} \label{DefineSupp} 
S(\mu) = \{1 \leq j \leq p: \mu(j) \neq 0\}. 
\end{equation} 
The procedure first estimates $S(\mu)$ by optimizing ($\|\cdot\|_1$:  vector $L^1$-norm)  
\begin{equation} \label{DefinehatS} 
\hat{S}^{(sa)}_N = \margmax_{\{ S \subset \{1, 2, \ldots, p\} : |S| = N \}}  \bigl\{   \| \sum\nolimits_{j \in S} x_j \|_1 \bigr\}, 
\end{equation} 
and then cluster  by aggregating all selected features 
$\hell^{(sa)}_{N}  = \sgn\bigl( \sum\nolimits_{j \in \hat{S}^{(sa)}_N} x_j  \bigr)$.\footnote{For any vector $x \in R^n$,  $\sgn(x) \in R^n$ is the vector where  the $i$-th entry  is $\sgn(x_i)$, $1 \leq i \leq n$ ($\sgn(x_i)  = -1, 0, 1$ according to $x_i < 0$,  $= 0$, or $> 0$).}  

An important special case is $N = p$, where $\hell_N^{(sa)}$ reduces to the method of {\it Simple Aggregation}  which we denote by $\hell_*^{(sa)}$.\footnote{The superscript ``sa" now loses its original meaning, but we keep it for consistency.}  This procedure targets on the case where the signals are weak but less sparse, so feature selection is hopeless. Note that $\hell_N^{(sa)}$ is generally NP-hard but $\hell_*^{(sa)}$ is not. 

The second method is IF-PCA, denoted by $\hell_q^{(if)}$, where $q > 0$ is a tuning parameter. The method  targets on the case where the signals are rare but individually strong.  To use $\hell_q^{(if)}$,  we first select features  using the   $\chi^2$-tests:  
\begin{equation} \label{DefineQ}
\hat{S}^{(if)}_q = \{1 \leq j  \leq p:  Q(j)  \geq \sqrt{2 q \log(p)} \},  \;\;     Q(j) = (\|x_j\|^2 - n) /\sqrt{2n}.  
\end{equation} 
We then obtain the first left singular vector $\xi^{(q)}$ of the post-selection data matrix $X^{(q)}$ (containing only columns of $X$ where the indices are in $\hat{S}^{(if)}_q$):  
\begin{equation} \label{Definexiq}  
\xi^{(q)} = \xi(X^{(q)}),  
\end{equation} 
and cluster by $\hell_{q}^{(if)}   = \sgn(\xi^{(q)})$.  
IF-PCA includes the classical PCA (denoted by $\hell_*^{(if)}$) as a special case, where the feature selection step is skipped, and $\xi^{(q)}$ reduces to the first singular vector of $X$.\footnote{The superscript ``if" now loses its original meaning, but we keep it for consistency.}   

%Note that  a tung free version of IF-PCA was introduced in Jin and Wang \cite{IFPCA} where it was applied successfully  
%to $10$ gene  microarray data sets. 

In Table~\ref{table:4method}, we compare all four methods. Note that for more complicated cases (e.g., the nonzero $\mu(j)$'s may be both positive and negative), we may consider a variant of $\hat{\ell}_N^{(sa)}$ which clusters by $\hat{\ell}_N^{(sa)} = \sgn( X \hat{\mu})$, with $\hat{\mu}$ being  $\margmax_{\{ \mu(j) \in \{-1, 0, 1\}, \|\mu\|_0 = N \}} \| X \mu\|_q$, where $q > 0$. If we let $q=1$ and restrict $\mu(j)\in\{0,1\}$, it reduces to the current $\hat{\ell}_N^{(sa)}$. Note that when $N = p$ and $q = 2$, approximately, $\hat{\mu}$ is proportional to the first right singular vector of $X$ and $\hat{\ell}_N^{(sa)}$ is approximately the classical PCA.   
Note also that  $\hat{\ell}_q^{(if)}$ can be viewed as the adaption of  IF-PCA  in Jin and Wang \cite{IFPCA}  to Model (\ref{model2}). The version in \cite{IFPCA} is a tuning free algorithm for analyzing microarray data and is much more sophisticated. The current version of IF-PCA is similar to that in Johnstone and Lu \cite{JohnstoneLu} but is also different in purpose and in implementation: the former is for estimating $\ell$ and uses the first {\it left} singular vector of the post-selection data matrix, and the latter is for estimating $\mu$ and uses the first {\it right} singular vector. 
The theory two methods entail are also very different.  See Sections \ref{subsec:spike} and  \ref{sec:discuss} for more discussion.   
%Note that IF-PCA  bears some resemblance to the method in \cite{JohnstoneLu} for both methods use $\chi^2$-%screening and post-selection PCA. However,  the focus of \cite{JohnstoneLu}  is on estimating  $\mu$ and the procedure %there is very different from IF-PCA: for example, they use the right singular vector of $X^{(q)}$ but we use the left %singular vector, and we have an additional clustering step on the eigenvector.   

%%%%%%%%%%%
%%%%%%%%%%%
\begin{table}[ht!]
\centering
\caption{Comparison of basic characteristics of four methods.  *: signals are comparably stronger but still weak. $\dagger$:  a tuning-free version exists. }
\scalebox{0.87}{
\begin{tabular}{rcccc}
\hline
Methods & Simple Aggregation    & Sparse  Aggregation  & Classical PCA   &  
IF-PCA  \\  
&  $\hell_{*}^{(sa)}$  &  $\hell_{N}^{(sa)}$  ($N \ll p$) & $\hell_{*}^{(if)}$ &  
 $\hell_q^{(if)} (q > 0)$ \\ 
\hline 
Signals  & less sparse/weak  & sparse/strong* & moderately sparse/weak & very sparse/strong \\ 
Feature selection & No & Yes  & No  & Yes \\  
Comp. complexity  &  Polynomial &  NP-hard & Polynomial & Polynomial \\ 
Need tuning & No & Yes & No & Yes$\dagger$  \\  
\hline
\end{tabular}
} 
\label{table:4method}
\end{table}

%%%%%%%%%
%%%%%%%%%
%%%%%%%%%
\subsection{Rare and Weak signal model} \label{subsec:ARW} 
To study all these limits, we invoke the Asymptotic Rare and Weak (ARW) model \citep{candes2009exact, DJ04, DJ08, IPT}. 
In ARW,  for two parameters $(\eps, \tau)$, we model the contrast mean vector $\mu$ by 
\begin{equation} \label{ARW1} 
\mu(j)  \stackrel{iid}{\sim} ( 1- \eps) \nu_0 + \eps \nu_{\tau}, \qquad 1 \leq j \leq p, 
\end{equation}
where $\nu_a$ denotes the point mass at $a$.  In Model (\ref{ARW1}),   all signals  have  the same sign and magnitude. Such an assumption can be largely relaxed;  see Sections~\ref{subsec:other-limits} and \ref{sec:discuss}. 
We use $p$ as the driving asymptotic parameter and tie $(n, \eps, \tau)$ to $p$ by fixed parameters. In detail,  fixing  $(\theta, \beta) \in (0,1)^2$ and $\alpha > 0$, we model 
\begin{equation} \label{ARW2} 
n = n_p = p^{\theta}, \qquad \eps =  \eps_p = p^{- \beta}, \qquad  \tau=\tau_p  = p^{-\alpha}. 
\end{equation} 
In our model, $n \ll p$ for we focus on the modern  ``large $n$, really large $p$" regime \citep{science}. The study can be conveniently extended  to the case of $n \gg p$.  

%%%%%%%%%%%
%%%%%%%%%%%
\subsection{Limits for clustering} 
\label{subsec:stat-limit}  
\begin{figure}[t!]
\begin{center}
\includegraphics[height = 1.9 in]{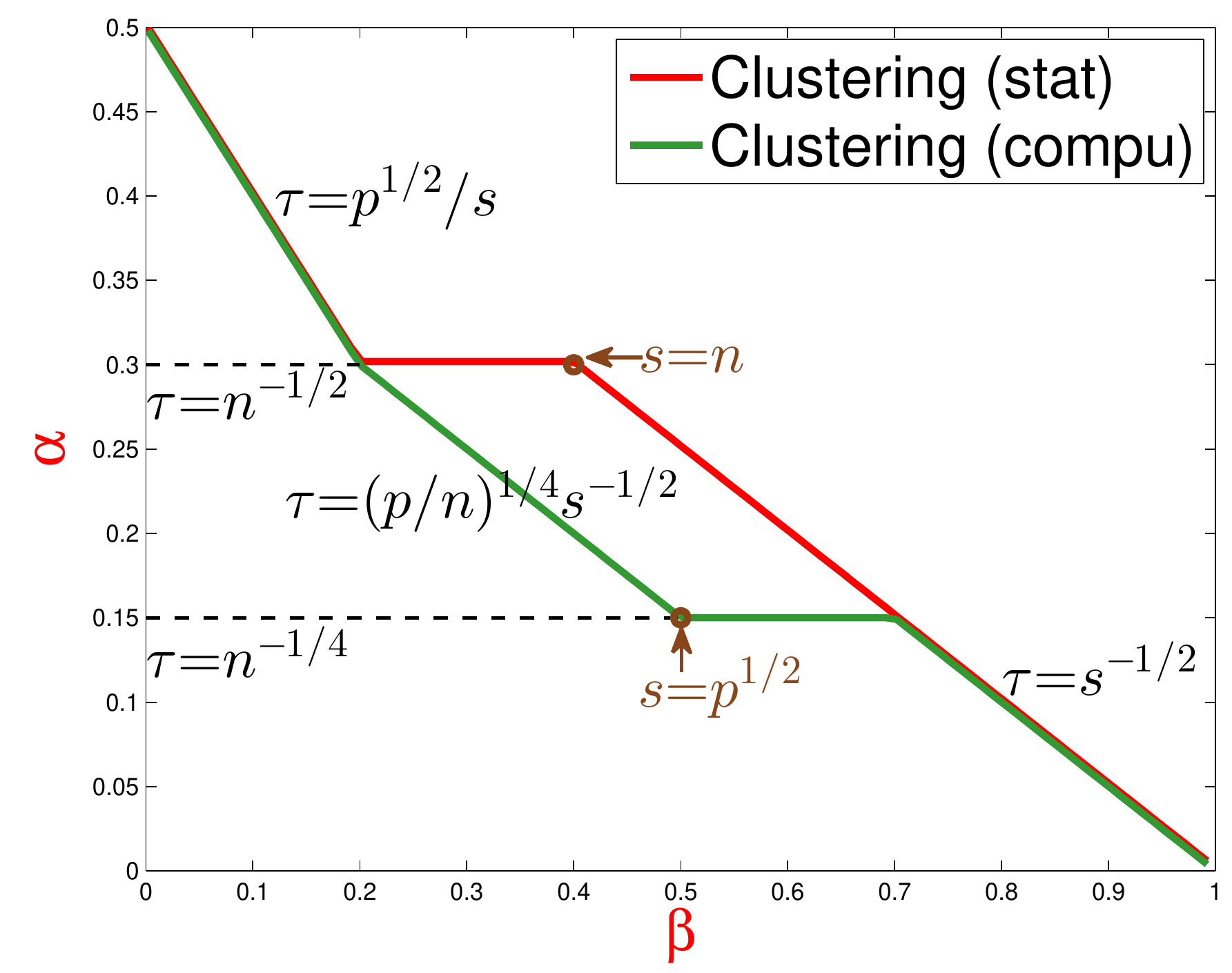} 
\includegraphics[height = 1.9 in]{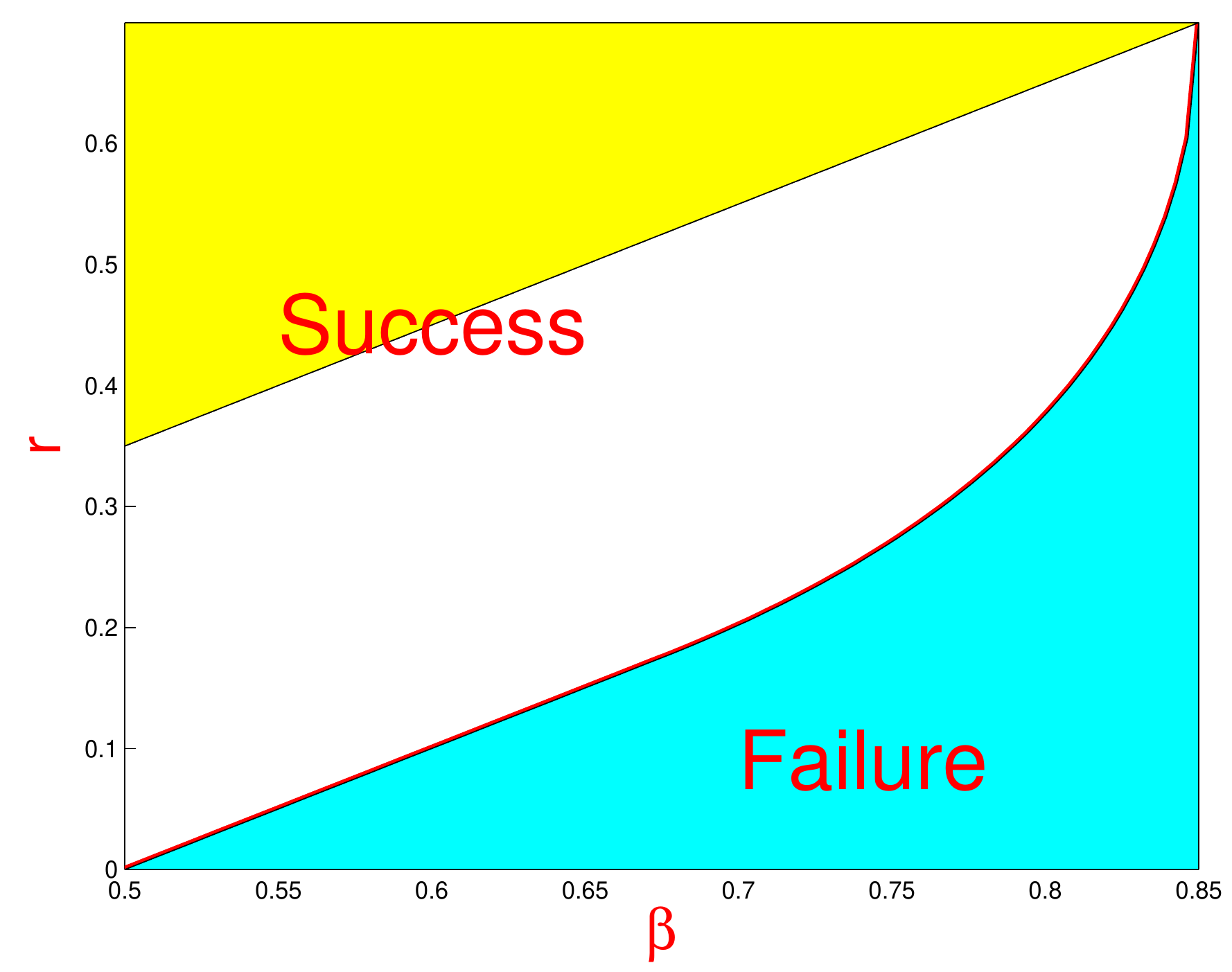} 
\end{center}
\caption{Left: the statistical limits (red) and the CTUB (green) for clustering ($s = p \eps_p$ is the expected number of signals). Right: phase transition of IF-PCA. White region: successful clustering is possible  but successful feature selection is impossible (using column-wise $\chi^2$ scores).  Yellow region: both successful clustering and feature selection are possible. } 
\label{fig:clulimit}
\end{figure}

%We first discuss statistical limits, and then the computationally tractable statistical limits.  
Let $\Pi$ be the set of all possible permutations on $\{-1,1\}$. 
For any clustering procedure $\hat{\ell}$ (where $\hat{\ell}_i$ takes values from $\{-1, 1\}$),  we measure the performance by the  Hamming distance:  
\begin{equation} \label{hamm1} 
\hamm_p(\hell, \alpha,\beta,\theta) =   n^{-1}  \inf_{\pi \in \Pi}   \biggl\{ \sum_{i = 1}^n  P(\hat{\ell}_i  \neq \pi \ell_i) \biggr\},    
\end{equation} 
where the probability is evaluated  with respective to  $(\mu,\ell, Z)$. 
Fixing $\theta \in (0,1)$, introduce a curve $\alpha = \eta^{clu}_{\theta}(\beta)$ in the $\beta$-$\alpha$ plane  by 
\[
\eta^{clu}_{\theta}(\beta)   = 
\left\{
\begin{array}{ll}
(1 - 2 \beta)/2,	 	& \qquad   \beta < (1-\theta)/2,     \\
\theta/2, 				& \qquad  (1 - \theta)/2 < \beta < (1 - \theta),     \\
(1 - \beta)/2,	&  \qquad  \beta > (1 - \theta).  
\end{array}
\right.
\]
%%%%%%%%%%
%%%%%%%%%%
%%%%%%%%%%
\begin{thm} \label{thm:clu-a} 
({\it Statistical lower bound}).\footnote{The ``lower bound" refers to the 
information lower bound as in the literature, not the lower bound for the curves in Figure \ref{fig:clulimit} (say). Same for the ``upper bound".} 
Fix $(\theta, \beta) \in (0,1)^2$ and $\alpha > 0$ such that $\alpha > \eta^{clu}_{\theta}(\beta)$. Consider the clustering problem for Models (\ref{model1})-(\ref{model2}) and (\ref{ARW1})-(\ref{ARW2}).  For any procedure $\hell$,   $\lim\inf_{p \goto \infty}\hamm_p(\hell, \alpha,\beta,\theta) \geq 1/2$. 
\end{thm}    
%The statistical lower bound is achievable by the methods of simple aggregation $\hell_*^{(sa)}$ and sparse aggregation $%\hell_N^{(sa)}$ in the 
%less sparse and the more sparse case, respectively.  
%%%%%%%%%%
%%%%%%%%%%
%%%%%%%%%%
\begin{thm} \label{thm:club} 
({\it Statistical upper bound for clustering}). 
Fix $(\theta, \beta) \in (0,1)^2$ and $\alpha > 0$ such that $\alpha < \eta^{clu}_{\theta}(\beta)$, and consider the clustering problem for Models (\ref{model1})-(\ref{model2}) and (\ref{ARW1})-(\ref{ARW2}).   As $p \goto \infty$,  
\begin{itemize} 
\item $\hamm_p(\hell_*^{(sa)}, \alpha,\beta,\theta) \goto 0$,  if $0 < \beta < (1 - \theta)/2$.  
\item $\hamm_p(\hell_N^{(sa)}, \alpha,\beta,\theta) \goto 0$,   if $(1 - \theta)/2 < \beta < 1$ and  $N=\lceil p\eps_p\rceil$.\footnote{$\lceil x \rceil$ denotes  the smallest integer that is no smaller than $x$.}  
\end{itemize} 
\end{thm}     
As a result, the curve $\alpha  = \eta^{clu}_{\theta}(\beta)$ divides the $\beta$-$\alpha$ plane 
into two regions: Region of Impossibility and Region of Possibility. In the former, the signals are so  weak that  successful 
clustering is impossible. In the latter, the signals are strong enough to allow successful clustering.

Consider computationally tractable limits.  We call a curve $r = \eta_{\theta}(\beta)$ in the $\beta$-$\alpha$ plane a {\it Computationally Tractable Upper Bound (CTUB)} if for any fixed $(\theta, \alpha, \beta)$ such that $\alpha < \eta_{\theta}(\beta)$,  there is a computationally tractable clustering method $\hell$ such that $\hamm_p(\hell, \alpha, \beta, \theta) \goto 0$. A CTUB $r = \eta_{\theta}(\beta)$ is tight if for any computationally tractable method $\hell$ and any fixed $(\theta, \alpha, \beta)$ such that $\alpha > \eta_{\theta}(\beta)$, 
$\lim\inf_{p \goto \infty}\hamm_p(\hell, \alpha,\beta,\theta) \geq 1/2$. In this case, we call $r = \eta_{\theta}(\beta)$ the  
Computationally Tractable Boundary (CTB). 
%%, and it is always below the statistical limit $\alpha = \eta_\theta^{clu}(\beta)$.  
Define 
\[
\tilde{\eta}^{clu}_{\theta}(\beta)   = 
\left\{
\begin{array}{ll}
(1 - 2 \beta)/2,	 	 	& \qquad   \beta < (1 - \theta)/2,\\ %\;\;   (\mbox{Case (a)})     \\
(1 + \theta - 2 \beta)/4,	 	 	& \qquad   (1 - \theta)/2 < \beta < 1/2,\\ %\;\;   (\mbox{Case (b)}),     \\ 
\theta/4,    &\qquad 1/2 < \beta < 1 - \theta/2,\\ % \;\;   (\mbox{Case (c)}), \\
(1 - \beta)/2,   &\qquad 1 - \theta/2  < \beta < 1.   %\;\;    (\mbox{Case (d)}). \\    
\end{array}
\right.
\]
\vspace{-2.5em}
\begin{thm} \label{thm:clu-c} 
({\it  A  CTUB  for clustering}). 
Fix $(\theta, \beta) \in (0,1)^2$ and $\alpha > 0$ such that $\alpha < \tilde{\eta}^{clu}_{\theta}(\beta)$, and consider the clustering problem for Models (\ref{model1})-(\ref{model2}) and (\ref{ARW1})-(\ref{ARW2}).   As $p \goto \infty$,  
\begin{itemize} 
\item $\hamm_p(\hell_*^{(sa)}, \alpha,\beta,\theta) \goto 0$,  if $0 < \beta < (1 - \theta)/2$.  
\item $\hamm_p(\hell_*^{(if)}, \alpha,\beta,\theta) \goto 0$,   if $(1 - \theta)/2 < \beta < 1/2$. 
\item $\hamm_p(\hell_q^{(if)}, \alpha,\beta,\theta) \goto 0$,   if $1/2 < \beta < 1$ and we take $q \geq 3$.   
\end{itemize} 
\end{thm}    
%%%%%%%%%%
%%%%%%%%%%
%%%%%%%%%%
%%%%%%%%%%
We now discuss CTB. We discuss the cases (a) $0 < \beta < (1 - \theta)/2$, (b) $(1 - \theta)/2 < \beta < 1/2$, (c) $1/2 < \theta < 1 - \theta/2$, and (d) $1 - \theta/2 < \beta < 1$ separately.   
Note that the CTB is sandwiched by two curves $\alpha=\eta^{clu}_{\theta}(\beta)$ and $\alpha = \tilde{\eta}_{\theta}^{clu}(\beta)$. In (a) and (d), $\tilde{\eta}^{clu}_{\theta}(\beta) = \eta^{clu}_{\theta}(\beta)$,  
so our  CTUB (i.e., CTUB given in Theorem \ref{thm:clu-c})  is tight.  
For (b), we are not sure but we conjecture that  our CTUB  is tight.\footnote{We know that CTB crosses two 
points $(\beta, \alpha) =(1/2, \theta/4)$ and $(\beta, \alpha) = ((1-\theta)/2, \theta/2)$.  A natural guess is that the CTB in  this part is a line segment connecting the  two points.}     
  For (c), we have good reasons to believe that  our  CTUB  is tight. 
In fact,  our model is intimately connected to the spike model \cite{JohnstoneLu}; see Section \ref{subsec:spike}.  
The tightness of our CTUB under the spike model has been well-studied (e.g., \cite{BR13learning, MW}).  Translating their results\footnote{Consider the hypothesis testing in the spike model. \cite{BR13learning} proves that, with the ``planted clique" conjecture, for $n<p$ and $s=o(\sqrt{p})$, if $\|\mu\|_0=s$ and $\|\mu\|^2\leq s\sqrt{\log(p)/n}$, 
there is no polynomial-time test that is powerful. In ARW, since $\|\mu\|^2\approx s\tau^2$, the above translates to (ignoring the logarithmic factor) $\alpha>\theta/4=\tilde{\eta}^{clu}_\theta(\beta)$.} to our setting suggests that there is a small constant $\delta > 0$ such that when $1/2 < \beta < 1/2 + \delta$,  our CTUB is tight. Note that for (c), the CTUB  $\alpha = \tilde{\eta}_{\theta}^{clu}(\beta)$ is flat.  
By the monotonicity of CTB  (see below), our CTUB is tight for (c). See Figure~\ref{fig:clulimit}. 

{\bf Remark} ({\it Monotonicity of  CTB}). We show the CTB is monotone  
in $\beta$ (with $\theta$ fixed).  Fix $\delta > 0$ and consider a new experiment, where for each column of 
the data matrix, we keep the column with probability $p^{-\delta}$ and replace it with an independent  column drawn 
from $N(0, I_n)$ with probability $1 - p^{-\delta}$.  Compare this with the original experiment. 
The parameters $(\alpha,\theta)$ are the same, but  $\beta$ has become $(\beta + \delta)$. 
The second experiment is harder, for it is the result of the original experiment by 
sub-sampling the columns. This shows that the CTB is monotone in $\beta$.  The monotonicity now follows 
by Le Cam's results on comparison of experiments \cite{LeCambook}. 

%Le Cam's  comparison of experiments is a very useful tool.  
%In Section~\ref{subsec:colored}, we use it to analyze a more difficult case where the noise is colored.
%Theorems \ref{thm:clu-a}-\ref{thm:clu-c} do not cover the most challenging case where the parameters $(\alpha, \beta)$ %fall exactly on either of  the two critical boundaries $\alpha = \eta_{\theta}^{clu}(\beta)$ and $\alpha = \tilde{\eta}_{\theta}%^{clu}(\beta)$.   We study these cases in Section \ref{sec:IF-PCA}. 

%{\bf Remark}. In the very sparse region $1<\beta<1/2$, other methods may also achieve the CTUB. However, we focus on IF-PCA due to %several reasons. First, it has been applied to a couple of gene microarray data sets and shows impressive performance \cite{IFPCA}. %Second, it is convenient to use in practice: computationally simple, with a tuning-free version and easy to modify (e.g., to incorporate %Efron's empirical null correction and to deal with heteroscedastic columns; see \cite{IFPCA}).   
%Third, for the most challenging case where parameters fall on the separating boundaries, IF-PCA is the only method known so far that is %theoretically tractable. 

%%%%%%%%%%%%%
%%%%%%%%%%%%%
%%%%%%%%%%%%%
\subsection{Phase transition for IF-PCA} \label{subsec:IF-PCA}    
IF-PCA is a flexible clustering method that is easy to use and computationally efficient.   In \cite{IFPCA}, we  developed a 
tuning free version of IF-PCA using Higher Criticism \cite{DJ04,DJ15, JinKe} and applied it to $10$  microarray 
data sets with satisfactory results. The success of IF-PCA in real data analysis 
motivates us to investigate the method in depth. To facilitate delicate analysis, we consider the version of IF-PCA in Section \ref{subsec:4method}, and reveal an interesting phase transition. 

To this end, we investigate a very  challenging case (not covered in Theorems \ref{thm:clu-a}-\ref{thm:clu-c}) where $(\alpha, \beta)$ 
fall exactly on the CTUB in Theorem \ref{thm:clu-c}:   
\begin{equation} \label{newcalib1} 
\alpha = \tilde{\eta}_{\theta}^{clu}(\beta).   \footnote{The case $\alpha < \tilde{\eta}_{\theta}^{clu}(\beta)$ is comparably easier to study, and the case 
$\alpha >  \tilde{\eta}_{\theta}^{clu}(\beta)$ belongs to the Region of Impossibility for computationally tractable methods; see our     conjectures.}  
\end{equation} 
Also, note that a key step in IF-PCA is the column-wise $\chi^2$-screening. In our model, a column $x_j$ is either distributed as $N(0, I_n)$ or $N(\tau_p \ell, I_n)$, where $\tau_p = p^{-\alpha}$. For the $\chi^2$-screening to be non-trivial, we 
further require that  
\begin{equation} \label{newcalib2}
1/2 < \beta < 1 - \theta/2. 
\end{equation} 
For $\beta$ in this range,  the curve $\alpha = \tilde{\eta}_{\theta}^{clu}(\beta)$ is flat, i.e., $\tilde{\eta}_{\theta}^{clu}(\beta) \equiv \theta/4$, and so $\tau_p = p^{-\theta/4} = n^{-1/4}$.  For $\beta$ outside this range, (\ref{newcalib1}) dictates that either  
$\tau_p \ll n^{-1/4}$ (so that the signals are too weak that the $\chi^2$-screening bounds to fail) or $\tau_p \gg n^{-1/4}$ (so that 
the signals are too strong that the $\chi^2$-screening is relatively trivial).  See Figure \ref{fig:clulimit}.

%{\color{red}We shall focus on the method IF-PCA for several reasons. First, it is the only method known so far that is %theoretically tractable when parameters fall on the separating boundary between Region of Possibility and Region of %Impossibility. Second, IF-PCA has been applied to several microarray data sets and shows impressive performance 
%(e.g., in \cite{IFPCA} and also Section~\ref{subsec:realdata}). Third, it is convenient to use in practice: computationally %simple, with a tuning-free version and easy to modify (e.g., to incorporate Efron's empirical null correction and to deal %with heteroscedastic columns; see \cite{IFPCA}).}

We now restrict our attention to (\ref{newcalib1})-(\ref{newcalib2}), where  we recall that  $\tau_p = p^{-\theta/4}$. 
To make the case more interesting, we adjust the calibration of $\tau_p$ slightly by an $O(\log^{1/4}(p))$ factor:  
\begin{equation} \label{newcalib3} 
\tau_p^* = p^{-\theta/4} (4 r \log(p))^{1/4}, \;\;\;  \mbox{where $0 < r  < 1$ is a fixed parameter}.  
\end{equation} 
With this calibration, the $\chi^2$-screening could be successful but non-trivial. 

Introduce the {\it standard phase function\footnote{It was introduced in the literature to study the phase transitions of multiple testing and classification with rare/weak signals.}} \cite{DJ04,DJ08}
\begin{equation} \label{stdphase} 
\rho^*(\beta)  = \left\{ 
\begin{array}{ll}
\beta  - 1/2, 			& \qquad  1/2 < \beta < 3/4,   \\
(1 - \sqrt{1 - \beta})^2, 	& \qquad  3/4  < \beta < 1. 
\end{array} 
\right. 
\end{equation} 
Define the {\it phase function} for IF-PCA by 
\begin{equation} \label{adjustphase} 
\rho_{\theta}^*(\beta)  = (1 - \theta) \cdot \rho^*\Big(1/2 + \frac{\beta - 1/2}{1 - \theta}\Big), \qquad 1/2 < \beta < 1 - \theta/2.  
\end{equation} 
For any two vectors $x$ and $y$ in $R^n$, let $\mathrm{cos}(x, y) = \bigl| \langle x/ \|x\|, y /\|y\|\rangle \bigr|$.  
%%%%%%%%%%
%%%%%%%%%%
%%%%%%%%%%
\begin{thm} \label{thm:IF-PCA} 
({\it Phase transition for IF-PCA}). 
Fix $(\theta,  \beta, \alpha, r) \in (0,1)^4$ and $q > 0$ such that (\ref{newcalib1})-(\ref{newcalib2}) hold. 
Consider  IF-PCA $\hell_q^{(if)}$ for Models (\ref{model1})-(\ref{model2}) and (\ref{ARW1})-(\ref{ARW2}),  where 
$\tau_p$ is replaced by the new calibration $\tau_p^*$ in (\ref{newcalib3}), and let $\xi^{(q)}$ be the leading left singular vector  as in (\ref{Definexiq}).  As $p \goto \infty$, 
\begin{itemize} 
\item If $r > \rho_{\theta}^*(\beta)$, then with probability at least $1-o(p^{-2})$, $\mathrm{cos}(\xi^{(q^*)}, \ell) \goto 1$ with $q^*=(\beta - \theta/2 + r)^2/(4r)$ for  $r > (\beta - \theta/2)/3$ and  $q^*=4r$ otherwise.
\item If $r < \rho_{\theta}^*(\beta)$, then with probability at least $1-o(n^{-1})$, there is a constant $c_0 \in (0,1)$ such that   
$\mathrm{cos}(\xi^{(q)}, \ell) \leq c_0$  for any fixed $0 < q < 1$.
\end{itemize} 
\end{thm}    
Theorem \ref{thm:IF-PCA} is proved in Section \ref{sec:IF-PCA}, using  delicate spectral analysis on the post-selection data matrix (and so the term of {\it post-selection Random Matrix Theory (RMT)}). Compared to many works on RMT where the data matrix has independent entries \cite{Vershynin}, the entries of the post-selection data matrix are complicatedly correlated, so the required analysis is more delicate. We conjecture that when $r<\rho_\theta^*(\beta)$,   $\mathrm{cos}(\xi^{(q)}, \ell) \goto 0$ for any fixed $0 < q < 1$. For now,  we can only show this for $q$ in a certain range; see the proof for details. 

Figure~\ref{fig:clulimit} (right) displays the phase diagram for IF-PCA. For fixed $(\alpha, \beta)$ in the interior of the white region, successful feature selection is impossible (by column-wise $\chi^2$-screening) but successful clustering is  possible. This shows that feature selection and clustering are related but different problems. 
   
% IF-PCA has a tuning free version \cite{IFPCA} which performs nicely with  gene microarray data. IF-PCA consists of a screening step and a post-selection PCA step.   In this paper, 
{\bf Remark}.  For the IF-PCA considered here, 
we use column-wise $\chi^2$-tests for screening which is computationally inexpensive. 
Alternatively, we may use some regularization methods for screening (e.g.,  \cite{Aspremont,Lei,ZHT}).   However, these methods are computationally more expensive, need tuning parameters that are  hard to set, and 
are designed for feature  selection, not clustering. For these reasons, it is unclear whether such alternatives  may really  help.

%%%%%%%%%
%%%%%%%%%
%%%%%%%%%
\subsection{Clustering when the noise is colored} \label{subsec:colored} 
Consider a new version of ARW where $(\ell, \mu)$ are the same as  in Models (\ref{model1}), (\ref{ARW1})-(\ref{ARW2}), but  Model (\ref{model2}) is replaced by a colored noise model 
%%%%%%%%%%%%%%%%%
%%%%%%%%%%%%%%%%%
%%%%%%%%%%%%%%%%%
\begin{equation} \label{newARW}
X = \ell \mu' +  A Z B, \qquad   Z_i(j) \stackrel{iid}{\sim} N(0,1), \qquad 1 \leq i \leq n, \; 1 \leq j \leq p,  
\end{equation} 
where $A$ and $B$ are two non-random matrices.   
\bed  \label{def:Lp}
We use $L_p > 0$ to denote a generic multi-$\log(p)$ term which may vary from occurrence to occurrence such that for any fixed $\delta  > 0$,  $L_p p^{-\delta} \goto 0$ and $L_p p^{\delta} \goto \infty$, as $p \goto \infty$.  
\eed 
%%%%%%%%%%
%%%%%%%%%%
%%%%%%%%%%
\begin{thm}  \label{thm:clu-d} 
({\it Statistical lower bound for clustering with colored noise}). Consider the ARW model (\ref{model1})-(\ref{model2}) and  (\ref{ARW1})-(\ref{ARW2}).   Theorem \ref{thm:clu-a} continues to hold if we replace the model (\ref{model2}) by (\ref{newARW})  
where $\max\{\|A\|, \|A^{-1}\| \} \leq L_p$ and $\max\{\|B\|,\|B^{-1}\| \} \leq L_p$. 
\end{thm} 
Theorem \ref{thm:clu-d} is proved in Section \ref{sec:proof}, using Le Cam's   comparison of experiments \cite{LeCambook}. 
The idea is to construct a new experiment that is easy to analyze and that the current one can be viewed as the result of 
adding noise to it. Since ``adding noise always makes the inference harder", 
analyzing the new experiment provides a lower bound we need for the current experiment.  The idea has been used in  Hall and Jin \cite{HallJin}, but for very different settings.

Consider the case $A = I_n$. In this case,  the matrix $A Z B$ has independent rows (but the columns may be correlated and  heteroscedastic), and all four methods we proposed earlier continue to work, except that in IF-PCA we need  $q\geq 3\max\{ \mathrm{diag}(B'B) \}$. 
The following theorem is proved in Section \ref{sec:sig}. 
\begin{thm}  \label{thm:clue} 
({\it Upper bounds for clustering with colored noise}). Consider the ARW model (\ref{model1})-(\ref{model2}) and  (\ref{ARW1})-(\ref{ARW2}). Theorems \ref{thm:club}-\ref{thm:clu-c} continue to hold if we replace the model (\ref{model2}) by (\ref{newARW}) with $A = I_n$ and $B$ such that $\max\{\|B\|,  \|B^{-1}\| \} \leq L_p$ and that all diagonals of $B'B$ is upper bounded by a constant $c>0$, where we set $q\geq 3c$ in IF-PCA. 
\end{thm} 
Practically, it is desirable to have a method that does not depend on the unknown parameter $c$. One way to attack this  is  to replace the column-wise $\chi^2$-test by a plug-in $\chi^2$-test where we estimate the variance column-wise by 
Median Absolute Deviation (say). However, such methods usually involve statistics of higher order  moments; see \cite{Ery} for discussions  along this line.

%%%%%%%%%
%%%%%%%%%
%%%%%%%%%
\subsection{Limits for signal recovery and hypothesis testing} \label{subsec:other-limits} 
For a more complete  picture, we study the limits for   signal recovery and  hypothesis testing.   

The goal of signal recovery is to  recover  the support of $\mu$. For any feature selector $\hat{S}$, we measure the error by the (normalized) Hamming distance $\hamm_p(\hat{S}, \alpha,\beta,\theta) = (p \eps_p)^{-1} \sum_{j = 1}^p [P(\mu(j)=0, j\in \hat{S})+ P(\mu(j)\neq 0, j\notin \hat{S})]$, where $p\eps_p$ is the expected number of signals.  
Define 
\[
\eta^{sig}_{\theta}(\beta) = 
\left\{
\begin{array}{ll}
\theta/2,		 				& \qquad  \beta < (1 - \theta),  \\
(1 + \theta - \beta)/4,	                         & \qquad  \beta > (1 - \theta),   
\end{array}
\right.
\]
and 
	\[
	\tilde{\eta}^{sig}_{\theta}(\beta) = 
	\left\{
	\begin{array}{ll}
	\theta/2,		 				& \qquad  \beta < (1 - \theta)/2,  \\
	(1 + \theta - 2 \beta)/4,   &\qquad   (1 - \theta)/2  <  \beta < 1/2,  \\ 
	\theta/4,	& \qquad  \beta >  1/2.   
	\end{array}
	\right.
	\]
The curve $r = \eta_{\theta}^{sig}(\beta)$ can be viewed as the counterpart of $r = \eta_{\theta}^{clu}(\beta)$, which 
divides the two-dimensional phase space into the Region of Impossibility and Region of Possibility. For any fixed  
$(\beta, \alpha)$ in the former and any $\hat{S}$,   $\hamm_p(\hat{S}, \alpha,\beta,\theta)\gtrsim  1$. For any fixed  $(\beta, \alpha)$ in the latter, there is an  $\hat{S}$ such that  $\hamm_p(\hat{S}, \alpha,\beta,\theta)\goto 0$. 
The curve $r = \tilde{\eta}_{\theta}^{sig}(\beta)$ can be viewed as the counterpart of $r = \tilde{\eta}_{\theta}^{clu}(\beta)$ and provides a CTUB for the signal recovery problem.  See Section \ref{sec:sig} for more discussion. 

The goal of (global) hypothesis testing is to test a null hypothesis $H_0^{(p)}$ that the data matrix $X$ has {\it iid} entries from $N(0, 1)$ against an alternative hypothesis $H_1^{(p)}$ that $X$ is generated according to Model (\ref{model2}).  
Define 
%%%%%%%%%%%%
%%%%%%%%%%%%
%%%%%%%%%%%%
\begin{equation} \label{hypphasefunction} 
\eta_{\theta}^{hyp}(\beta)  = \max\{\eta_{\theta}^{hyp, 1}(\beta), \eta_{\theta}^{hyp, 2}(\beta)\},   \;\;
\tilde{\eta}_{\theta}^{hyp}(\beta)  = \max\{\eta_{\theta}^{hyp, 1}(\beta), \frac{\theta}{4}\},
\end{equation} 
and $\eta_{\theta}^{hyp,1}(\beta)  = (2 + \theta - 4\beta)/4$,  $\eta_{\theta}^{hyp,2} = \min\{\theta/2, (1 + \theta-\beta)/4\}$. 
Similarly,  the curve $r = \eta_{\theta}^{hyp}(\beta)$ divides the two-dimensional phase space into the Region of Impossibility and Region of Possibility. Fix $(\beta, \alpha)$ in the former,  the sum of Type I and Type II errors $\gtrsim 1$ for any testing procedures. Fixing $(\beta, \alpha)$ in the latter,   
there is a test such that  the sum of Type I and Type II errors tends to $0$.  
Also, the curve $r = \tilde{\eta}_{\theta}^{hyp}(\beta)$ provides a CTUB for the hypothesis testing problem. See Section~\ref{sec:hyp} for more discussion.   

The statistical limits for hypothesis testing here are different from those in Arias-Castro and Verzelen \cite{Ery}. For the less sparse case ($\beta < \theta/2$), the signal strength needed in our model is weaker, because all signals have the same sign. More interestingly, we find a phase transition phenomenon that is not seen in \cite{Ery}:  when $\theta < 2/3$, there are three segments for the statistical limits; when $\theta > 2/3$, there are only two segments.
\footnote{The curve $r = \eta_{\theta}^{hyp}(\beta)$  is the maximum of the boundary achievable by Simple Aggregation (a line segment)  and that by Sparse Aggregation (two line segments). Depending on where two boundaries cross each other,  $r = \eta_{\theta}^{hyp}(\beta)$  may consist of $2$ or $3$ line segments.} 
 
The tightness of CTUB for signal recovery and hypothesis testing can be addressed similarly to that for clustering. 
For signal recovery, 
the CTUB is tight in the less sparse case ($0<\beta<(1-\theta)/2$) for it matches the statistical limits; we have good reasons to believe it is tight in the sparse case ($1/2<\beta<1$), due to results in \cite{BR13learning,MW}; we are not sure for the moderate sparse case $((1-\theta)/2<\beta<1/2)$. For hypothesis testing, we have similar arguments except that the cases of ``less sparse" and ``moderate sparse" refer to that of  $0<\beta<(2-\theta)/4$ and that of  $(2-\theta)/4<\beta<1/2$, respectively.

\begin{figure}[t!]
\begin{center}
\includegraphics[height = 1.85 in]{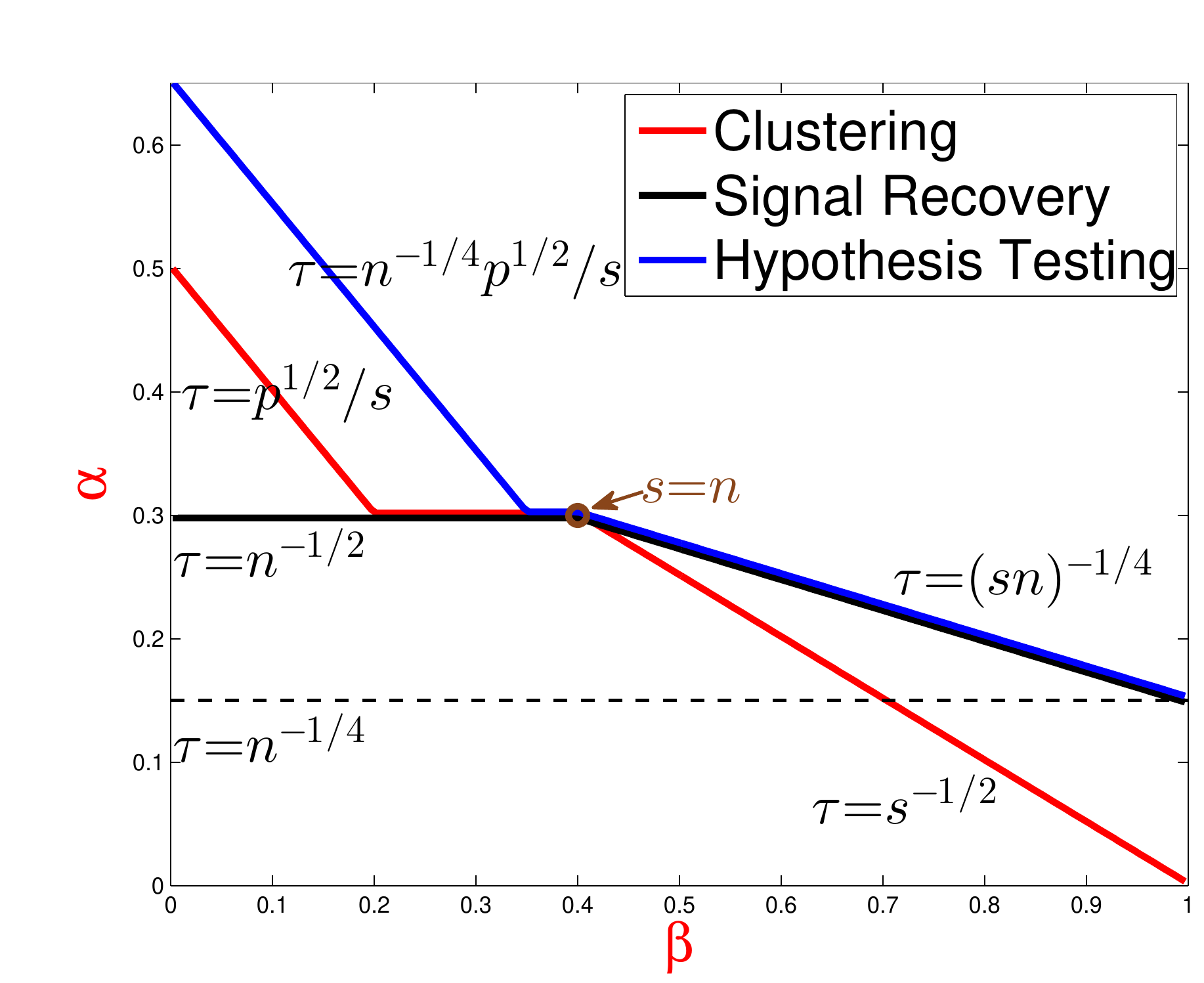}
\includegraphics[height = 1.85 in]{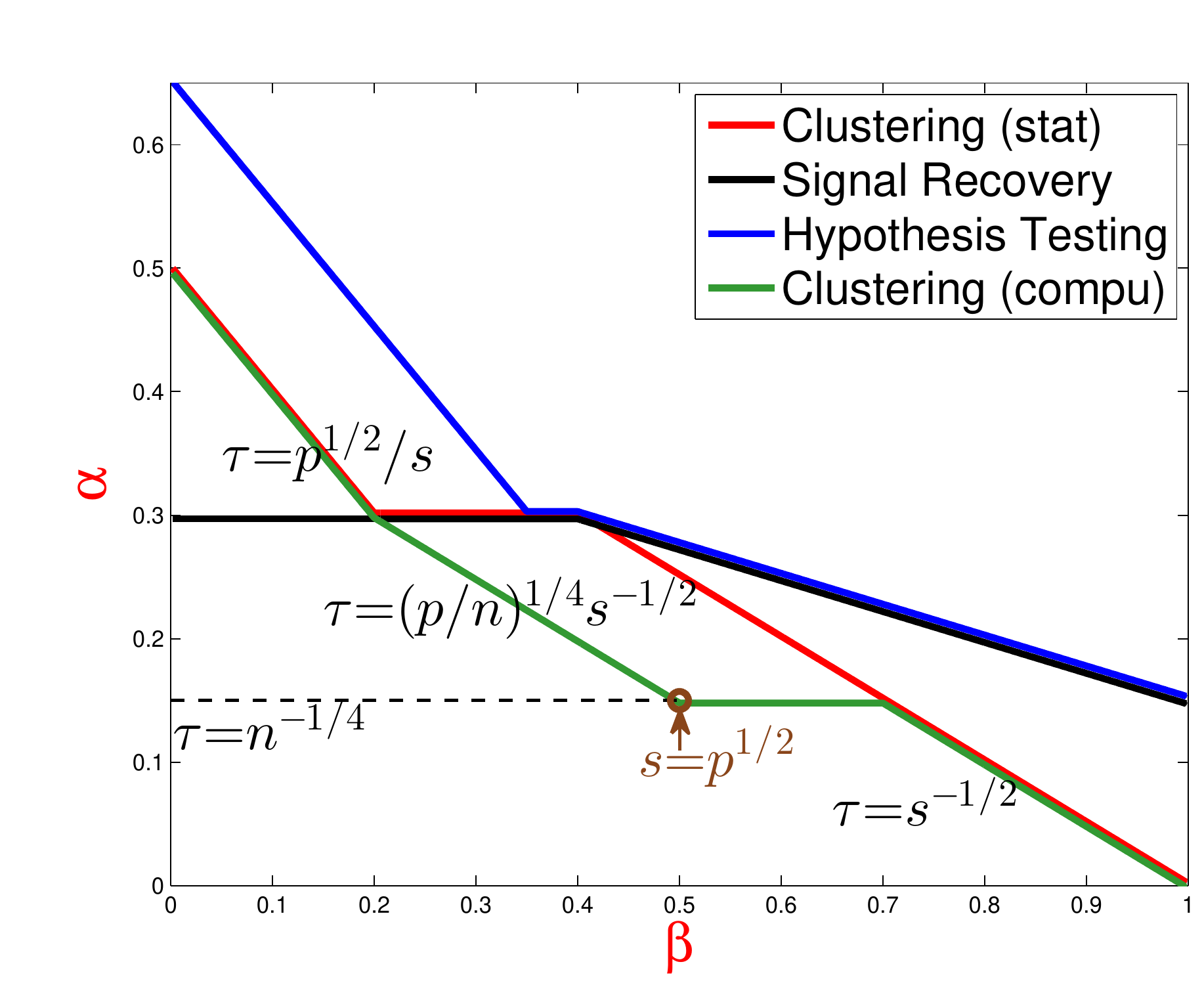} \\ 
\includegraphics[height = 1.83 in]{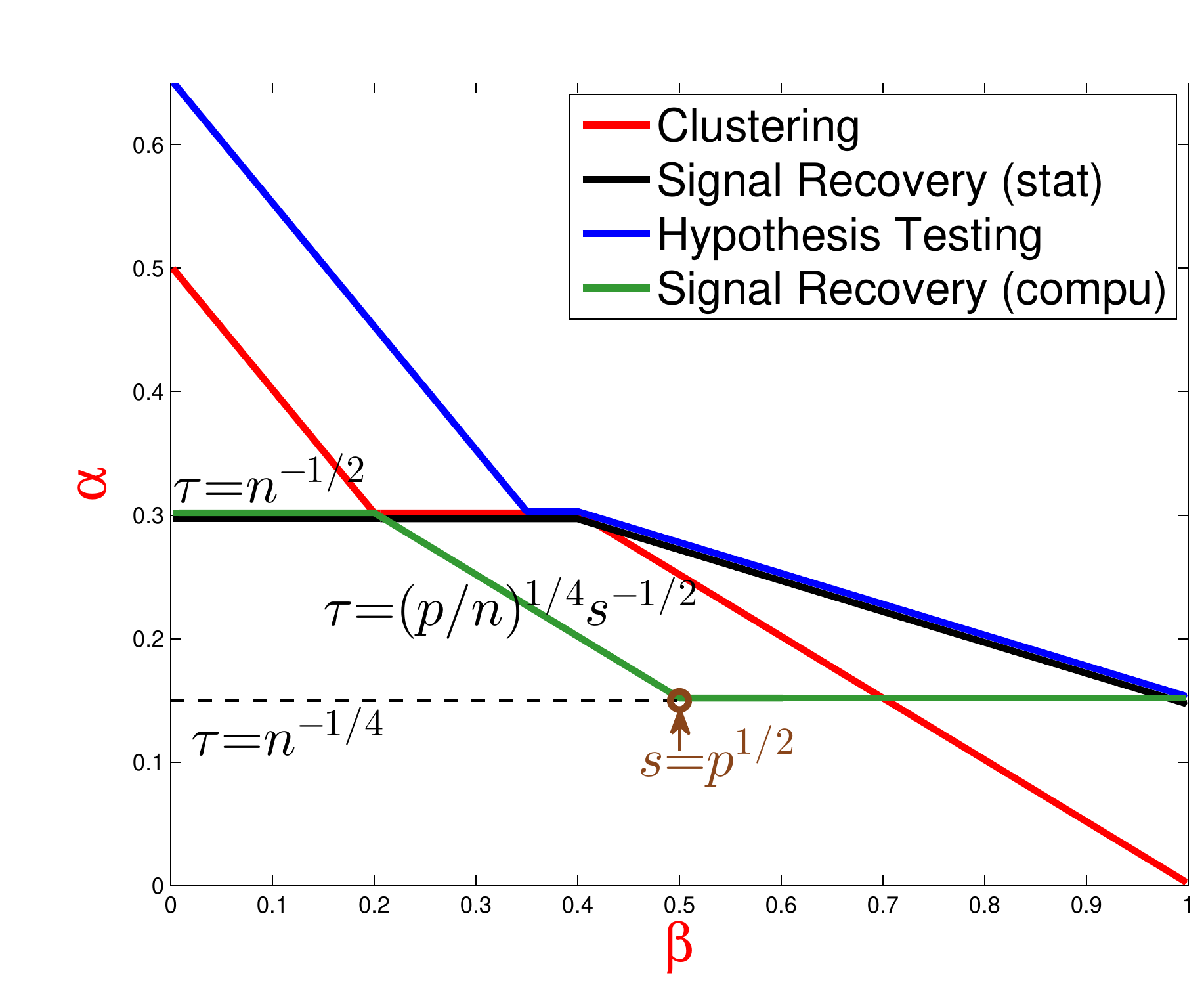}  
\includegraphics[height = 1.80 in]{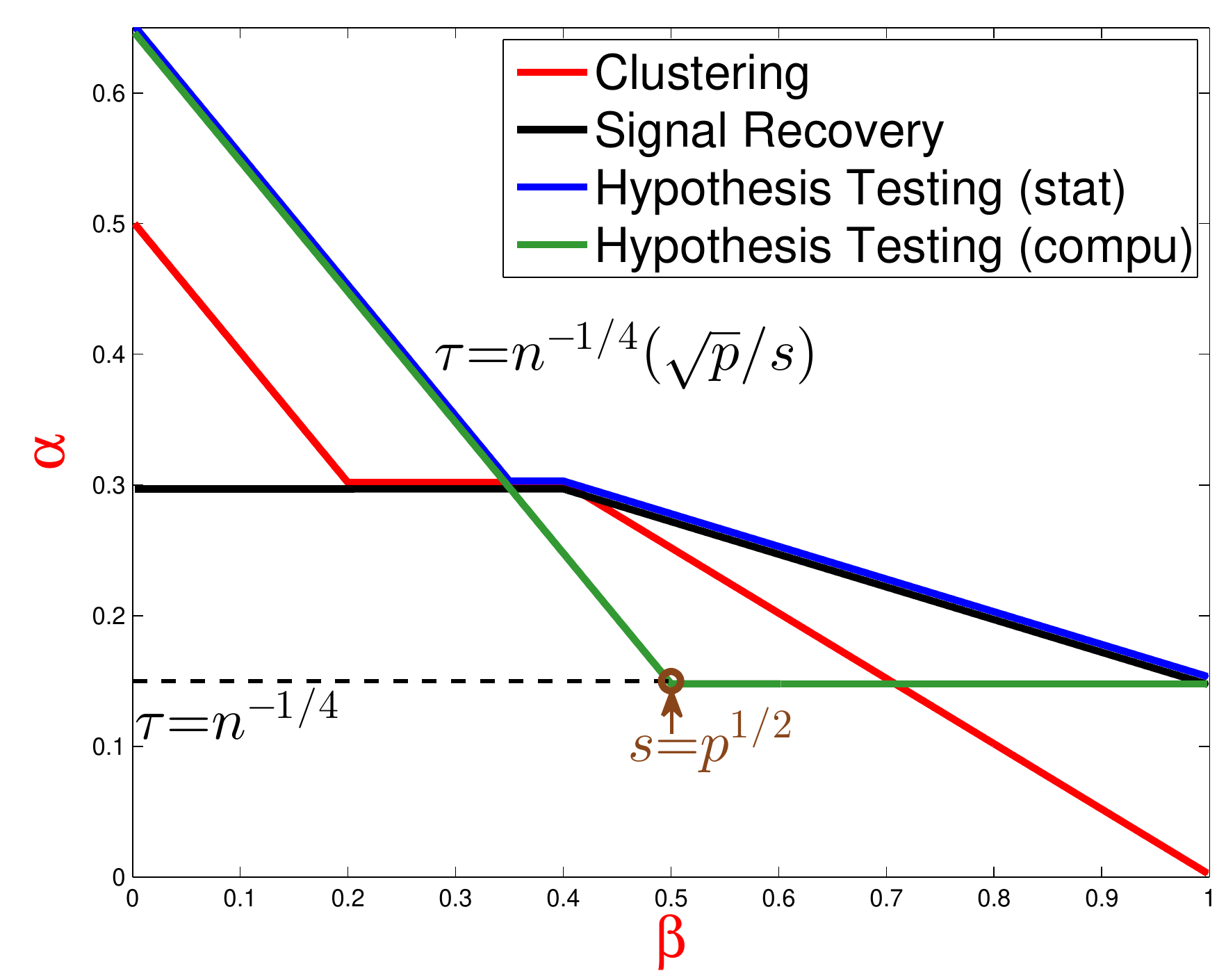} 
\end{center}
\caption{Top left: statistical limits for  
clustering (red), signal recovery (black), and hypothesis testing (blue); $s = p \eps_p$.  
Other three panels: CTUB for clustering (top right), signal recovery (bottom left)  and hypothesis testing (bottom right), respectively (the three statistical limits in the top left panel are also shown for  comparison).}
\label{fig:complimit}
\end{figure}

Figure~\ref{fig:complimit} compares the limits for all three problems: clustering, signal recovery, and hypothesis testing. See details therein.   

{\bf Remark}. Consider an extension of ARW where (\ref{ARW1}) is replaced by a 
more complicated signal configuration: $\mu(j) \stackrel{iid}{\sim} (1 - \eps) \nu_0  + a \eps \nu_{-\tau} + (1 - a) \eps \nu_{\tau}$, where $0 \leq a \leq 1/2$ is a constant ($a = 0$: original ARW).  When $0 <  a < 1/2$, our results on statistical limits and CTUB for all three problems continue to hold, provided with a slight change in the definition of the Hamming distance for signal recovery. 
The case of  $a = 1/2$ is more delicate,  but the changes in statistical limits (compared to the case of $a = 0$)   can be explained with  Figure \ref{fig:complimit} (top left): (a) the black curve (signal recovery) remains the same, (b) the red curve (clustering) remains the same, except for the segment on the left is replaced by $\tau^4 = p/(ns^2)$, (c) for the blue curve (hypothesis testing), the right most segment remains the same, while the other two segments coincide  with those of the red curve. The  CTUBs also change correspondingly.  See Appendix \ref{app:C} for a more detailed discussion.  
%%%%%%%%%%
%%%%%%%%%%
%%%%%%%%%%
%%%%%%%%%%

%%%%%%%%%%%%
%%%%%%%%%%%%
%%%%%%%%%%%%
 
%%%%%%
%%%%%%
%%%%%%
\subsection{Practical relevance and a real data example} \label{subsec:realdata}  
The relatively idealized model we use allows very delicate analysis, but also raises practical concerns. 
In this section, we investigate IF-PCA with a real data example and illustrate that many ideas in previous sections are 
relevant in much broader settings. 

%%%%%%%%%%%%
%%%%%%%%%%%
\begin{table}[ht!]
\centering
\caption{The clustering errors for Leukemia data with different numbers of selected features. Rows highlighted correspond to the threshold choices that yield lowest clustering errors.}
\scalebox{0.88}{
\begin{tabular}{ll||ll||ll}
$\#$\{selected features\}  & Errors  & $\#$\{selected features\}  & Errors  & $\#$\{selected features\}  & Errors \\
\hline
1 & 34  &   1419 & 3 &  2847 & 5  \\ 
347 & 8  & {\bf 1776}  & {\bf 1} & 3204 & 7 \\ 
704 & 6 & {\bf 2133} & {\bf 1}    &  3561 & 11  \\
1062 & 5 & {\bf 2490} & {\bf 1}    \\
\hline
\end{tabular}
} 
\label{table:Diffk}
\end{table}

We use the leukemia data set on gene microarrays. 
This data set was cleaned by Dettling \cite{Dettling}, consisting of $ p = 3571$ measured genes for 
$n = 72$   samples from two classes:  47  from ALL (acute lymphoblastic leukemia), and 25  from AML (acute myeloid leukemia). The data set is available at \url{www.stat.cmu.edu/~jiashun/Research/software/GenomicsData/ALL}. 

To implement IF-PCA, one noteworthy difficulty is the heteroscedasticity across genes in the data set. We apply IF-PCA with small modifications.
In detail, arrange the data matrix as $X = [x_1, \ldots, x_p]$ as before. Let $\bar{x}(j)  = (1/n) \sum_{i = 1}^n  x_j(i)$, $m(x_j) = \mathrm{median}(x_j)$ and $d(j)=\mathrm{median}\{|x_j(1) -m(x_j)|,\cdots, |x_j(n) - m(x_j)|\}$ be the Median Absolute Deviation (MAD). We normalize by $x_j^*(i) = 0.6745 \cdot (x_j(i)  - \bar{x}(j)) / d(j)$, $1 \leq i \leq n, 1 \leq j \leq p$\footnote{The value $0.6745$ is such that $E[(x^*_j(i))^2] = 1$ when $x_j(i) \sim N(0, \sigma^2)$ for any $\sigma > 0$.}. For $q > 0$ to be determined, we select feature $j$ if and only if $(2n)^{-1}|\|x_j^*\|^2-n|> \sqrt{2q\log(p)}$.  
We then obtain the leading left singular vector $(\xi^*)^{(q)}$ of the post-selection data matrix $[x_1^*,\cdots, x_p^*]$ and cluster by applying the standard $k$-means algorithm to the leading eigenvector. 
In the last step, we can also cluster by the sign vector of $(\xi^*)^{(q)}$ and the results are similar. The $k$-means algorithm has a slightly better performance. 

Table~\ref{table:Diffk} displays  the clustering errors for different numbers of selected features (each   corresponds to a choice of $q$). The table suggests that IF-PCA works nicely, with an error rate as low as $1/72$, if $q$ is set appropriately. 

\begin{figure}[t!]
\begin{center}
\includegraphics[width = 1.62 in, height = 1.7 in]{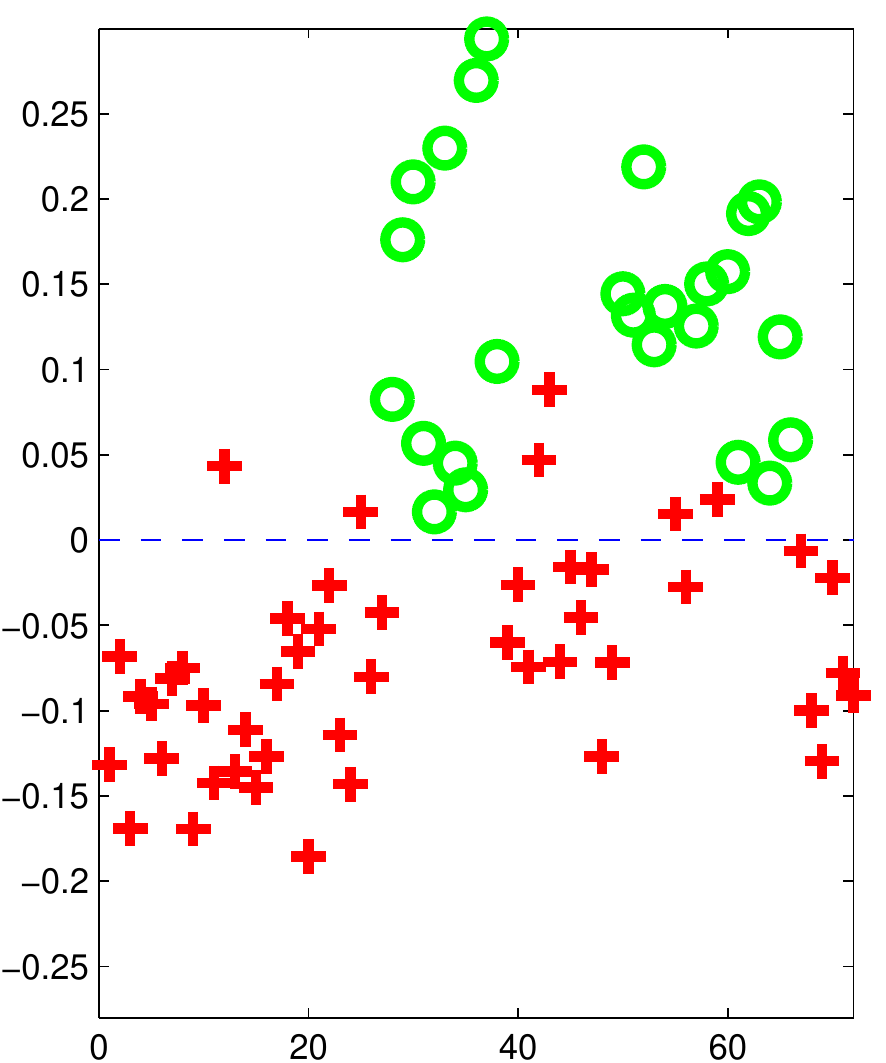}
\includegraphics[width = 1.62 in, height = 1.7 in]{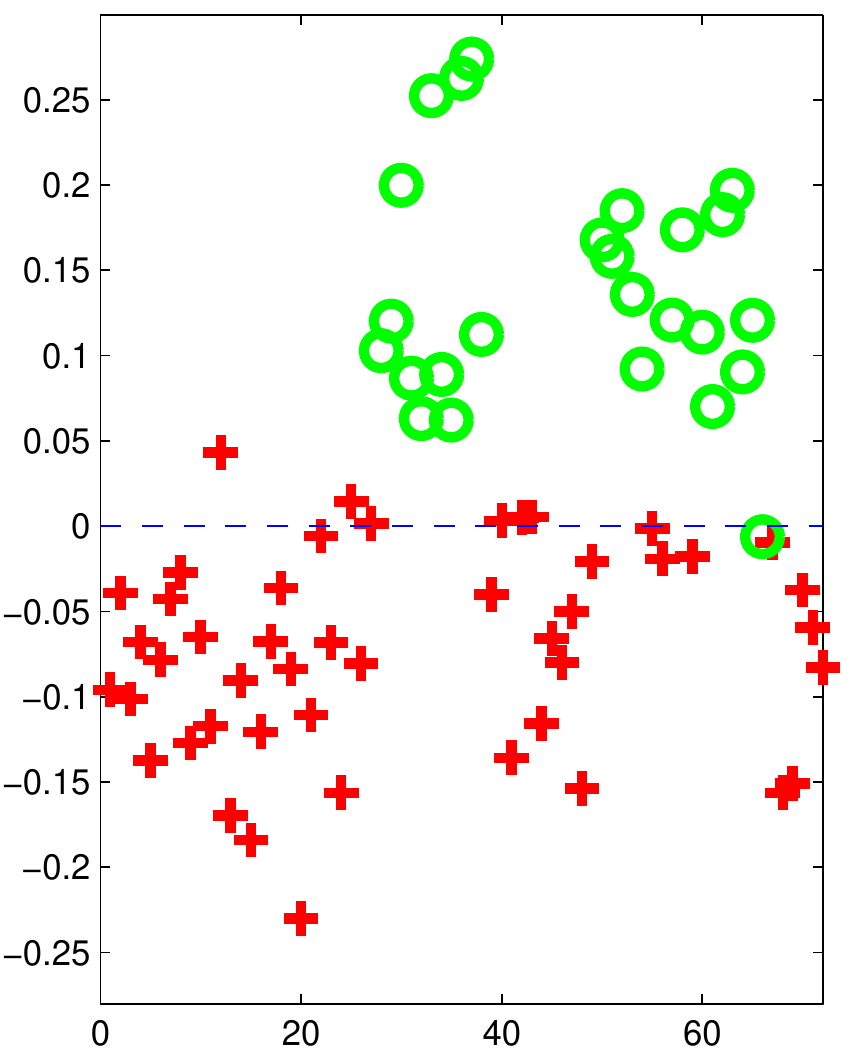}
\includegraphics[width = 1.62 in, height = 1.7 in]{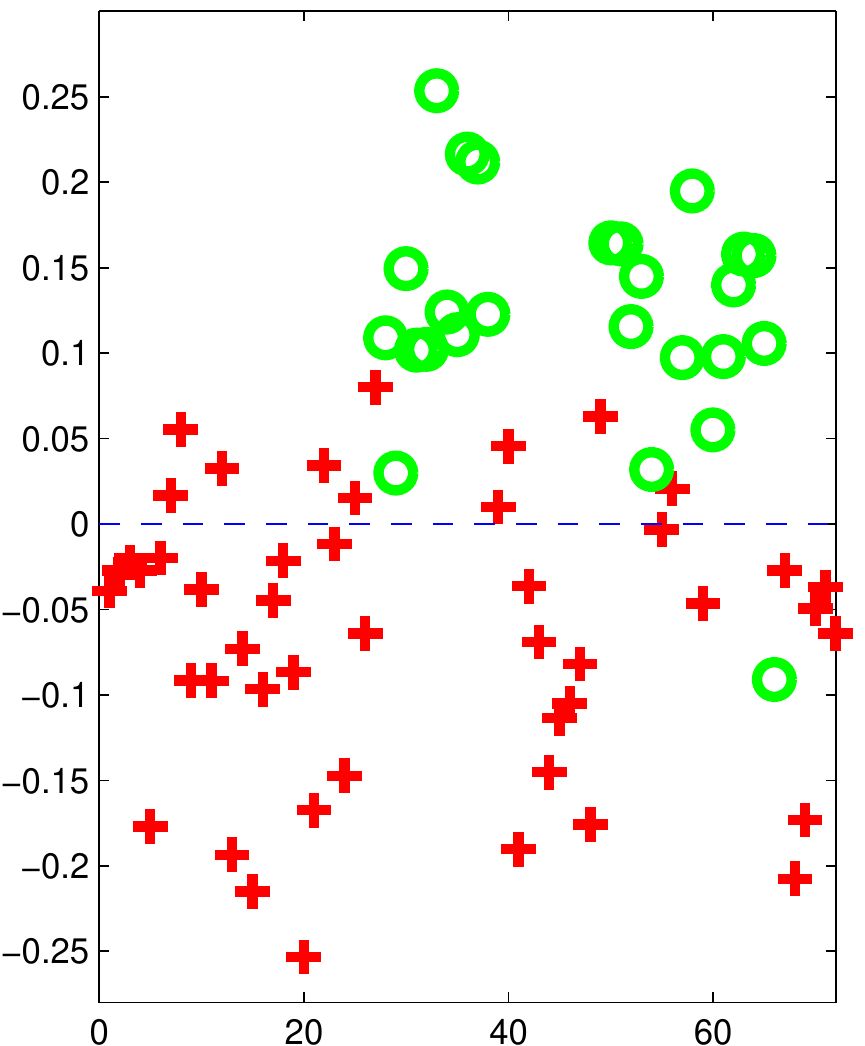}
\end{center}
\caption{Leading left singular vector of the data matrix $X$ with very few features selected by FDR choice (left; 931 features chosen), with ideal number of features selected (middle; 2133 features chosen), and without feature selection (right).
 $y$-axis: entries of the left singular vector, $x$-axis: sample indices.
Plots are based on Leukemia data, where red and green dots represent samples from the two classes ALL and AML, respectively.}
\label{fig:leuk}
\end{figure}

Figure \ref{fig:leuk} compares $(\xi^*)^{(q)}$ for three choices of $q$: (a) the $q$ determined by applying the FDR controlling procedure  \cite{BH95} with the FDR parameter of $.05$ and simulated $P$-values under the null $x_j\sim N(0, I_n)$,  (b) the $q$ associated with the ideal number of selected features (see Table \ref{table:Diffk}),  and (c) the  $q$ corresponding to  classical PCA (any $q$ that allows us to skip the feature selection step works). This suggests that  IF-PCA works well if  $q$ is properly set.  \footnote{A hard problem is how to set $q$ in a data-driven fashion. This is addressed in \cite{JKW}.}

%%%%%%%%%%
%%%%%%%%%%
%%%%%%%%%%
%%%%%%%%%%

We  compare IF-PCA with  classical methods of $k$-means and hierarchical clustering \cite{HTF},   $k$-means++ (a recent revision of the classical $k$-means; \cite{kmeans+}),\footnote{For $k$-means, we use the built-in  Matlab package (parameter `replicates' equals 30). For $k$-means++, we run the program 30 times, and compute the average clustering errors.}
SpectralGem (classical PCA applied to $X^*$; \cite{Lee}),   and sparse $k$-means (a modification of $k$-means with sparse feature weights in the objective; \cite{skmeans}).  
The error rates are in Table \ref{table:Error}, suggesting IF-PCA is effective in this case.

{\begin{table}
\caption{Comparison of clustering errors (Leukemia data). Columns 2--7:  numerator is the number of  clustering errors, and denominator is the number of subjects.}
\centering
\scalebox{0.9}{
\begin{tabular}{ l |cccccc}
\hline
Method   & $k$-means &$k$-means++& Hierarchical  & SpectralGem & Sparse $k$-means & IF-PCA \\
\hline
Error Rate     &  20/72  & 18.5/72  &  20/72 &  21/72  &  20/72 & 1/72   \\
\hline
\end{tabular}}
\label{table:Error}
\end{table}
}

%%%%%%%%%%%
%%%%%%%%%%%
%%%%%%%%%%%
%\subsection{Comparison to literature on sparse PCA and sparse clustering}  \label{subsec:spike} 

%%%%%%%%%
%%%%%%%%%
%%%%%%%%%
\subsection{Comparison to works on the spike model}  
\label{subsec:spike} 
In our model \eqref{model1}-(\ref{model2}), if we replace the Bernoulli model for $\ell_i$ in (\ref{model1})  by a Gaussian model where $\ell_i \stackrel{iid}{\sim} N(0, \sigma^2)$,  then it becomes 
the spike model (Johnstone and Lu \citep{JohnstoneLu}).

In the spike model,  while $\ell_i$'s are also of  interest, the feature vector $\mu$ captures most of the attention:  
most recent works on the spike model (e.g., \citep{amini2009high, Lei, WLL}) have been focused on   signal recovery (and especially, sparse PCA).    The two problems, signal recovery and clustering, are different. There are parameter settings where successful clustering is possible but successful signal recovery is impossible, and there are settings where the opposite is true; see Sections \ref{subsec:IF-PCA} and \ref{subsec:other-limits}. Therefore, a direct  extension of sparse PCA methods to clustering does not always work well.   

Our work is also different from existing works on the spike model  in terms of motivation and validation. 
Our model is motivated by cancer (subject) clustering, where the class labels $\ell_i$'s 
can be conveniently validated in many applications (e.g., see Section \ref{subsec:realdata}). In contrast, it is not easy to find real data sets where the feature vector $\mu$ is known, so it is comparably harder to validate the methods/theory on  signal recovery or sparse PCA. 
Given the growing awareness of reproducibility and replicability \cite{DMRrep},   it becomes  increasingly more important to 
develop methods and theory that can be directly validated by real applications.
In a sense, our model extends the spike model to a new direction, and it helps strengthen (we hope) the ties between the recent theoretical interests on the spike model with real applications.

%We reveal some interesting new insights:  the  problems of clustering and signal recovery are intertwined, 
%and it is preferable to study them simultaneously, not separately. For example, 
%in the less sparse case, clustering is easier, so we can first do clustering and then use results for signal recovery. In the %more sparse case, we should do the opposite; see Section \ref{sec:sig}.   

%Our work is also closely related to the literature of sparse clustering (e.g., \cite{WassermanSparse,CH,skmeans}). 
%Their works and ours share the same spirit of exploring feature sparsity to improve clustering performance, but the %methods and theoretical framework are very different. In particular, we focus on characterizing the precise statistical %limits of clustering. We defer the careful comparison with literature to Section~\ref{sec:discuss}. 

%%%%%%%%%%%%%%%
%%%%%%%%%%%%%%%
%%%%%%%%%%%%%%%
%%%%%%%%%%%%%%%
%%\subsection{Summary}  \label{subsec:summa}   
%%%%%%%%%%%
%%%%%%%%%%%
%%%%%%%%%%%

%%%%%%%%%%%%%%
%%%%%%%%%%%%
%%%%%%%%%%%%
%%%%%%%%%%%% 
\subsection{Content and notations} \label{subsec:content} 
Section~\ref{sec:IF-PCA} studies the phase transition  of IF-PCA, where we prove Theorem~\ref{thm:IF-PCA}. 
Section~\ref{sec:sig} studies the statistical limits for signal recovery, where we prove Theorems \ref{thm:club}, \ref{thm:clu-c},   
\ref{thm:clue}, as well as Theorems~\ref{thm:sigb}-\ref{thm:sigc} (to be introduced). 
Section~\ref{sec:hyp} studies the  statistical  limits for hypothesis testing, where  
we  prove Theorems \ref{thm:hypb}-\ref{thm:hypc} (to be introduced). 
 Section~\ref{sec:proof}  studies  the lower bounds for all three problems and proves Theorems \ref{thm:clu-a} and \ref{thm:clu-d}, as well as Theorems \ref{thm:siga} and \ref{thm:hypa} (to be introduced).  Other proofs  are in  the Appendix.  Section~\ref{sec:discuss}  is for discussion. 

In this paper, $L_p > 0$ denotes a generic multi-$\log(p)$ term; see Section \ref{subsec:colored}.  
%When $x$ is a scalar, $\sgn(x) = -1, 0, 1$ when $x < 0$, $x = 0$, and $x > 0$, respectively. When $x$ is an $n \times 1$ %vector,  $\sgn(x) =(\sgn(x_1),\cdots, \sgn(x_n))'$.  
When $\xi$ is a vector, $\| \xi\|_q$ denotes the vector $L^q$-norm, $0 \leq q \leq \infty$ (the subscript is dropped for simplicity if $q = 2$).   
When $\xi$ is a matrix, $\|\xi\|$ denotes the matrix spectral norm, and $\|\xi\|_F$ denotes the matrix Frobenius norm. 
 For two vectors $\xi, \eta$, $\langle \xi, \eta\rangle$ denotes the inner product of them, and $\ang(\xi, \eta) = |\langle \xi/\|\xi\|, \eta/\|\eta\| \rangle|$. 
For any two probability densities $f$ and $g$,  $\|f-g\|_1$  and $H(f, g)$ are the 
$L^1$-distance and the Hellinger distance, respectively. For any real value $a$, $\lceil a\rceil$ is the smallest integer that is no smaller than $a$. 
We say two positive sequences $a_n\sim b_n$, $a_n \lesssim b_n$ and $a_n \gtrsim b_n$ if $\lim_{n\goto\infty}a_n/b_n= 1$, $\lim\sup_{n \goto \infty} a_n/b_n \leq 1$ and $\lim\inf_{n \goto \infty} a_n/b_n \geq 1$, respectively. 
For two sets $A,B$, $A\Delta B =(A\backslash B)\cup (B\backslash A)$.

%%%%%%%%
%%%%%%%%
%%%%%%%%
\section{Phase transition for IF-PCA}   \label{sec:IF-PCA}    
\setcounter{equation}{0} 

In this section, we prove Theorem~\ref{thm:IF-PCA}.   
Our proofs need very precise characterization of the spectra of the post-selection Gram matrix $X^{(q)} (X^{(q)})'$. 
Specifically, we need both a tight upper bound on the range of the spectra of $X^{(q)} (X^{(q)})'$ (Lemma~\ref{lem:postZ-1}) and a tight 
lower bound for the largest eigenvalue of $X^{(q)} (X^{(q)})'$ (Lemma \ref{lem:postZ-2}).
The main challenges are that, due to feature selection,  
\begin{itemize} 
\item the entries of $X^{(q)}$ are no longer independent,  
\item the conditional distribution of each survived column is unclear.  
\end{itemize} 
For this reason, existing results on RMT do not apply directly and we need to develop new theory on post-selection RMT. 
Our analysis adapts that in Vershynin \cite{Vershynin} and uses the results of covering number in Rogers \cite{Rogers}.

{\bf Remark.} For the spike model, there are results about the spectra of a different post-selection Gram matrix $(X^{(q)})'X^{(q)}$ (e.g., Thereom 2 of \cite{JohnstoneLu}). Since feature selection is column-wise, the leading eigenvectors of $(X^{(q)})'X^{(q)}$ and $X^{(q)}(X^{(q)})'$ have very different behaviors. Moreover, the settings of \cite[Theorem 2]{JohnstoneLu} implicitly force $X^{(q)}$ to have much more rows than columns (which we call the ``skinny" case), but our results do not have such a restriction.

To show the claim, it suffices to show the claim for any fixed realization of $(\ell, \mu)$ in the event 
\[
D_p = \{\mu: \bigl| |S(\mu)| - p \eps_p  \bigr| \leq \sqrt{ 6 p \eps_p\log(p) } \};    
\]
note that $P(D_p^{c})  = O(p^{-3})$ and  the event only has a negligible effect. 
Fixing $0 < q < 1$ and a realization of $(\ell, \mu)$ in $D_p$.  Let $\hat{S}_q^{(if)}(\ell, \mu)$ be the set of all survived features.  In our model,  $X = \ell \mu'  + Z$, and $Z = [z_1, z_2, \ldots, z_p]$. 
Introduce a vector $\mu^{(q)}  = \mu^{(q)}(\ell, \mu) \in R^p$ and a matrix $Z^{(q)}=[z_1^{(q)},\cdots, z_p^{(q)}] \in R^{n,p}$ by 
\[
\mu^{(q)}(j)=\mu(j) \cdot  1\{j\in \hat{S}^{(if)}_q(\ell, \mu) \}, \;\;  z_j^{(q)}=z_j \cdot 1\{j\in \hat{S}^{(if)}_q(\ell, \mu) \}, \;\;  1\leq j\leq p,   
\]
and so the post-selection data matrix $X^{(q)} = X^{(q)}(\ell, \mu)$, viewed as an $n\times p$ matrix with many zero columns, satisfies 
\[
X^{(q)}(\ell, \mu) = \ell \mu^{(q)}(\ell, \mu) + Z^{(q)}(\ell, \mu). 
\] 

Fixing $(\beta,\theta,r) \in (0,1)^3$ and  $q>0$, and assuming $z\sim N(0, I_n)$, introduce 
$m^{(q)}_0(\mu) = (p-|S(\mu)|) \cdot P(\|z\|^2>n+2\sqrt{qn\log(p)})$,
$ m^{(q)}_1(\ell,\mu)= |S(\mu)| \cdot P(\|z + \tau_p^*\ell \|^2>n+2\sqrt{qn\log(p)})$,
and 
$m^{(q)}(\ell, \mu) =  m_0^{(q)}(\ell, \mu) +m_1^{(q)}(\ell, \mu)$. Note that $m_0^{(q)}(\ell,\mu)$ and $m_1^{(q)}(\ell,\mu)$ are the expected numbers of survived useless/useful features, respectively.  
We also need the following counterpart of $m^{(q)}(\ell, \mu)$: 
\begin{align*}
m_*^{(q)}&(\ell, \mu) = (p-|S(\mu)|)\cdot n^{-1} E\big( \|z\|^21\{\|z\|^2>n+2\sqrt{qn\log(p)}\} \big)\cr
& + |S(\mu)| \cdot n^{-1} E\big( \|z\|^21\{\|z+\tau_p^*\ell \|^2>n+2\sqrt{qn\log(p)}\} \big). 
\end{align*}
%where we write $m^{(q)}_*(\ell,\mu)=m^{(q)}_*$ for short.

The dependence on $(\ell, \mu)$ is tedious, so  for notational simplicity, we may drop them without further notices.  

The term $m^{(q)}$ is the expected number of selected features, and plays an important role.  
By tail properties of chi-square distributions (see Section~\ref{subsec:prelim}), with probability 
$1 - O(p^{-3})$, 
\begin{equation} \label{If-pca-asy} 
m_*^{(q)} \sim m^{(q)}  \sim L_p [p^{1 - q} +  p \eps_p p^{-[(\sqrt{q} - \sqrt{r})_+]^2}], 
\end{equation} 
where as before $L_p$ is a generic multi-$\log(p)$ term.  
Recalling $n = p^{\theta}$, define 
\[
\tilde{q}(\beta, \theta, r) =\left\{\begin{array}{lr}
\max\{1 - \theta, (\sqrt{1 - \beta - \theta} + \sqrt{r})^2\}, & \beta<1-\theta,\\
1-\theta, & \beta>1-\theta. 
\end{array}\right. 
\] 
By (\ref{If-pca-asy}) and basic algebra, it is seen that there are two different cases: \begin{itemize}
\item ({\it ``Fat"}). When $q<\tilde{q}(\beta,\theta,r)$,   $m^{(q)}/n\goto\infty$ and $X^{(q)}$ has much more columns than rows.  
\item ({\it ``Skinny"}). When $q> \tilde{q}(\beta,\theta,r)$,   $m^{(q)}/n\goto 0$ and $X^{(q)}$ has much more rows than columns. 
\end{itemize}  
\begin{lemma} \label{lem:postZ-1}
({\it Upper bound for the range of eigenvalues of $Z^{(q)}(Z^{(q)})'$}). Suppose conditions of Theorem~\ref{thm:IF-PCA} hold.  There exists a universal constant $C > 0$ such that  for any fixed $q > 0$, as $p\goto\infty$, conditioning on any realization of $(\ell,\mu)$ from the event $D_p$,  with probability at least $1 - O(p^{-3})$, 
\begin{itemize}
\item (``Fat" case). When $q <  \tilde{q}(\beta,\theta,r)$, all eigenvalues of $Z^{(q)}(Z^{(q)})'$ fall between $m^{(q)}_*\pm [C\sqrt{nm^{(q)}\log(p)}+ o(m_1^{(q)})]$. 
\item (``Skinny" case). When $q  >  \tilde{q}(\beta,\theta,r)$,  all nonzero eigenvalues  of $Z^{(q)}(Z^{(q)})'$ fall between $n   \pm C\sqrt{nm^{(q)}\log(p)}$. 
\end{itemize} 
\end{lemma}
 
{\bf Remark}. Noting that $m^{(q)}$ is the expected number of columns of $Z^{(q)}$, our results are very similar to the well-known results on eigenvalues of RMT in the case where we have an $n \times m^{(q)}$ matrix with $iid$ $N(0,1)$ entries.  
However, we need more sophisticated proofs, as the rows of $Z^{(q)}$ are dependent and the distribution of the columns of $Z^{(q)}$ is unknown and hard to characterize.  

For the ``fat" case, it turns out that Lemma \ref{lem:postZ-1} is insufficient:  we need both 
an improved upper bound on the range (with the $\sqrt{\log(p)}$ factor eliminated) and a
lower bound on the leading eigenvalue.   
%%%%%%%%%
%%%%%%%%%
\begin{lemma} \label{lem:postZ-2} 
({\it Improved bound (``fat" case)}). 
Suppose the conditions of Theorem~\ref{thm:IF-PCA} hold and $r<\rho^*_\theta(\beta)$. There exist constants $c_1>c_2>0$ such that as $p\to\infty$, for any fixed $q>0$, conditioning on any realization of $(\ell,\mu)$ from the event $D_p$, with probability $1-O(n^{-2})$, 
\begin{itemize}
\item All singular values of $Z^{(q)}(Z^{(q)})'$ fall between $m^{(q)}_*\pm c_1\sqrt{nm^{(q)}}$; 
\item $\lambda_{\max}(Z^{(q)}(Z^{(q)})')\geq m_*^{(q)}+c_2\sqrt{nm^{(q)}}$.
\end{itemize}
\end{lemma}

%{\bf Remark}. The condition $r<\rho^*_\theta(\beta)$ is to further eliminate the $o(m_1^{(q)})$ term in Lemma~\ref{lem:postZ-1}. In fact, the above lemma holds for a wider range of $r$, but we assume $r<\rho^*_\theta(\beta)$ for ease of presentation. 
%{\bf Remark}. The condition $r<\rho^*_\theta(\beta)$ is to guarantee that for any $0<q<1$, the number of  selected useful features is comparably small so that ...

We now prove Theorem \ref{thm:IF-PCA}. We show the cases of $r  >  \rho^*_\theta(\beta)$ (Region of Possibility) and $r  < \rho^*_\theta(\beta)$ (Region of Impossibility) separately.

%%%%%%%%%%%%
%%%%%%%%%%%%
%%%%%%%%%%%%
\subsection{Region of Possibility} \label{subsec:IF-PCAUB}  
Consider the case $r > \rho_{\theta}^*(\beta)$.  
Recall that 
\[
 q = q^*(\beta, \theta,r) =  \left\{\begin{array}{ll}
4r, & r < (\beta - \theta/2)/3,\\
\frac{(\beta - \theta/2 + r)^2}{4r}, &(\beta - \theta/2)/3 \leq r < 1.
\end{array}
\right.
 \] 
Let $\xi^*$ be the first left singular vector of $X^{(q)}$ at $q = q^*(\beta, \theta, r)$. 
The goal is to show 
\[
\ang(\ell, \xi^*) \goto 1. 
\]
Write 
\beq  \label{decompose}
X^{(q)}(X^{(q)})' = \|\mu^{(q)}\|^2 \ell\ell' + Z^{(q)}(Z^{(q)})' + A,  
\eeq
where $A=\ell(\mu^{(q)})'(Z^{(q)})' + Z^{(q)}\mu^{(q)}\ell'$ for short. On the right hand side of \eqref{decompose}, the first matrix has a rank $1$,  with $n \|\mu^{(q)}\|^2$ being the only nonzero  eigenvalue and 
 $\ell$ being the associated eigenvector.   
In our model,  the expectation of $\|\mu^{(q)}\|^2$ is equal to $(\tau_p^*)^2 m_1^{(q)}$, 
where by tail properties of chi-square distributions (see Section~\ref{subsec:prelim}), $m_1^{(q)} = L_p  p^{1 - \beta} p^{-[(\sqrt{q} - \sqrt{r})_{+}]^2}$ 
with overwhelming probabilities. It follows that with a probability at least $1 - O(p^{-3})$, 
\[
n\|\mu^{(q)}\|^2 \gtrsim n (\tau_p^*)^2 \cdot m_1^{(q)}  \sim L_p p^{\Delta(q,\beta,\theta,r)},
\]
where 
$\Delta(q,\beta,\theta,r) = 1+\theta/2-\beta-[(\sqrt{q} - \sqrt{r})_{+}]^2$.  
Compare this with (\ref{decompose}).   By perturbation theory in matrices  
\footnote{We use \cite[Proposition 1] {CMW}, a variant of the sine-theta theorem \cite{sin-theta}. By that proposition, if $\hat{\xi}$ and $\xi$ are the respective leading eigenvectors of two symmetric matrices $\hat{G}$ and $G$, where $G$ has a rank $1$, then $\|\hat{\xi}\hat{\xi}'-\xi\xi'\|\leq 2\|G\|^{-1}
\|\hat{G}-G\|$. We also note that for two unit-norm vectors $\hat{\xi}$ and $\xi$, $\cos(\hat{\xi},\xi)\goto 1$ if and only if $\|\hat{\xi}\hat{\xi}'-\xi\xi'\|\goto 0$ by linear algebra.} \cite{CMW, sin-theta}, 
 to show the claim, it suffices to show that there is a scalar $a^*$ (either random or non-random) and a constant $\delta^*>0$ so that \footnote{We have used the fact that adding/subtracting a multiple of the identity matrix does not affect the eigenvectors.}
\[
\| Z^{(q)}(Z^{(q)})' + A - a^* I_n\|  \leq  L_p  p^{\Delta(q,\beta,\theta, r)-\delta^*}. 
\] 
To this end, note that by triangle inequality, 
\[
\| Z^{(q)}(Z^{(q)})' + A- a^* I_n\|  \leq \| Z^{(q)}(Z^{(q)})' - a^* I_n \| + \|A\|. 
\]  
%\vspace{-1.5em}
The following lemma is proved in the Appendix. 
%%%%%%%%%%%%
\begin{lemma}  \label{lem:matrixA}
Suppose conditions of Theorem~\ref{thm:IF-PCA} hold. For any fixed $q>0$, as $p\goto\infty$, conditioning on any realization of $(\ell,\mu)$ from the event $D_p$, with probability $1-O(p^{-3})$, 
$\|\ell(\mu^{(q)})'(Z^{(q)})' + Z^{(q)}\mu^{(q)}\ell'\| \leq Cn\tau_p^* \sqrt{m_1^{(q)}}$. 
 \end{lemma}
The key  to the proof is to control $\|Z^{(q)}\mu^{(q)}\|_\infty$ using the Bernstein inequality \cite{Wellner86} and to study the distribution of  $Z^{(q)}$. See \cite{3Phasesupp} for details. 
 
Now, when $q > \tilde{q}(\beta,\theta,r)$, we are in the ``skinny" case, combining Lemmas \ref{lem:postZ-1} and \ref{lem:matrixA},  we have that with probability at least $1 - O(p^{-3})$,  
\[
\| Z^{(q)}(Z^{(q)})' + A\|   \lesssim  n + C\Big(n\tau_p^*\sqrt{m_1^{(q)}} + \sqrt{n m^{(q)} \log(p)}\Big) \leq L_p p^{\frac{\theta}{2}  +\frac{1}{2} \max\{\theta, \Delta(q,\beta,\theta,r)\}}. 
\]
In the last inequality, we have used (\ref{If-pca-asy}) which indicates that $m_0^{(q)}=L_pp^{1-q}$ and $m_1^{(q)}=L_p p^{1-\beta-[(\sqrt{q}-\sqrt{r})_+]^2}$.  
By the condition of $r > \rho_{\theta}^*(\beta)$, it can be shown that $\Delta(q,\beta,\theta,r)>\theta$, 
and the claim follows by letting $a^* = 0$ and $\delta^* =\frac{\Delta-\theta}{2}$. 
When $q <  \tilde{q}(\beta,\theta,r)$, we are in the ``fat" case.  Combining Lemmas \ref{lem:postZ-1} and \ref{lem:matrixA}, with probability at least $1 - O(p^{-3})$,  
\begin{align*}
\| Z^{(q)}(Z^{(q)})' + A - m_*^{(q)}I_n \| &\leq C \Big( n\tau_p^*\sqrt{m_1^{(q)}} + \sqrt{nm^{(q)}\log(p)} + n^{-1}m_1^{(q)}\Big)\cr
&  \leq L_p p^{\frac{\theta}{2}+ \frac{1}{2}\max\{\theta,\Delta(q,\beta,\theta,r), 1-q\}}+ p^{\Delta(q,\beta,\theta,r)-\frac{3\theta}{2}}. 
\end{align*}
By the condition of $r > \rho^*_{\theta}(\beta)$, it can be shown that $\Delta(q,\beta,\theta,r)>\max\{\theta, \frac{\theta+1-q}{2}\}$,   
 and the claim follows  by letting $a^* = m_*^{(q)}$ and $\delta^* =\min\{ \frac{\Delta-\theta}{2},\Delta-\frac{1-q+\theta}{2}, \frac{3\theta}{2}\}$.

%%%%%%%%%
%%%%%%%%%
%%%%%%%%% 
\subsection{Region of Impossibility} 
Consider the case $r  <  \rho^*_\theta(\beta)$. Fix $0 < q < 1$.  Recall that $\xi^{(q)}$ is first left singular vector of $X^{(q)}$. The goal is to show that 
\begin{equation} \label{ifpcatoshow0} 
\ang(\ell, \xi^{(q)}) \leq c_0 < 1, \qquad \mbox{for any }   0 < q < 1,  
\end{equation} 
where $c_0$ is a universal constant independent of $q$.   
Denote for short 
$H = X^{(q)} (X^{(q)})'$,  $H_0 = Z^{(q)} (Z^{(q)})'$, $\xi=\xi^{(q)}$, and $\tilde{\ell} = \ell/ \|\ell\|$. 
Let the eigenvalues of $H$ be $\lambda_1(H) \geq \lambda_2(H) \geq \ldots \geq  \lambda_n(H)$.    
Write 
\[
\tilde{\ell} = a \xi + \sqrt{1 - a^2} \eta, \qquad \mbox{for a unit-norm vector $\eta$ such that $\eta \perp \xi$}.   
\]
Note that $\xi' H \eta = \lambda_1 \xi' \eta = 0$ and $\eta'H\eta \geq \lambda_n$,   we have 
$\tilde{\ell}' H \tilde{\ell} = a^2  \xi'  H \xi + 2 a \sqrt{1 - a^2} \xi' H \eta + (1 - a^2) \eta' H \eta \geq a^2 \lambda_1  + (1 - a^2) \lambda_n$.  
Rearranging it gives 
$a^2 \leq 1- [\lambda_1(H) -\tilde{\ell}'H\tilde{\ell}]/[\lambda_1(H) -\lambda_n(H)]$.  
Note that 
$\ang(\ell, \xi) = |a|$.  So to show (\ref{ifpcatoshow0}), it suffices to show there  
\begin{equation} \label{ifpcatoshow} 
\frac{\lambda_1(H) - \tilde{\ell}' H \tilde{\ell}}{\lambda_1(H) - \lambda_n(H)} \geq 1-c_0^2,   \qquad  \mbox{for some constant $c_0 \in (0,1)$}. 
\end{equation} 
The following lemma is proved in the Appendix. 
%%%%%%%%%
%%%%%%%%%
\begin{lemma} \label{lem:v'H0v}
Suppose $r<\rho^*_\theta(\beta)$ and the conditions of Theorem~\ref{thm:IF-PCA} hold. As $p\to\infty$, for any fixed $q>0$, conditioning on any realization of $(\ell,\mu)$ from the event $D_p$, for any $v\in \mathcal{S}^{n-1}$, with probability $1-O(p^{-3})$,  $|v' H_0 v -   m_*^{(q)}|\leq C\sqrt{m^{(q)}\log(p)}$, 
and 
%%%%%%%%
\begin{equation} \label{HH0} 
\|H - H_0\|  = \left\{ 
\begin{array}{ll} 
o(n), &\qquad \mbox{$q > \tilde{q}(\beta, r, \theta)$ (``skinny" case)},    \\ 
o(\sqrt{n m^{(q)}}), &\qquad \mbox{$q < \tilde{q}(\beta, r, \theta)$ (``fat" case)}. 
\end{array}
\right. 
\end{equation} 
\end{lemma}

We now show (\ref{ifpcatoshow}).   Similarly, let  $\lambda_1(H_0) \geq \lambda_2(H_0) \geq \ldots \geq  \lambda_n(H_0)$ be the eigenvalues of $H_0$.  We prove  for the cases of $q > \tilde{q}(\beta, r, \theta)$  and $q  < \tilde{q}(\beta, r, \theta)$  separately.  
Consider the first case. This is the ``skinny" case where $m^{(q)} \ll n$. By Lemma \ref{lem:postZ-1} and the first claim of Lemma \ref{lem:v'H0v},  with probability $1 - O(p^{-3})$,  
$\lambda_1(H_0) \sim n$, $\lambda_n(H_0) \geq 0$ and $\tilde{\ell}' H_0 \tilde{\ell} = o(n)$. 
By the second claim of Lemma~\ref{lem:v'H0v}, $\|H-H_0\|=o(n)$.
Combining the above with Weyl's inequality \cite{weyl} (i.e., $\max_{1\leq i\leq n}|\lambda_i(H) - \lambda_i(H_0)| \leq \|H - H_0\|$), we have  
\[
\lambda_1(H)-\lambda_n(H)\leq n+o(n), \qquad \lambda_1(H) - \tilde{\ell}' H \tilde{\ell} \geq n-o(n).  
\]
Inserting these into (\ref{ifpcatoshow}) gives the claim.  

Consider the second case.  This is the ``fat" case and  $m^{(q)}\gg n$.  
By Lemma  \ref{lem:postZ-2} and the first claim of Lemma \ref{lem:v'H0v}, there is
$\lambda_1(H_0) - \lambda_n(H_0) \leq c_2 \sqrt{n m^{(q)}}$ and $ \lambda_1(H_0) - \tilde{\ell}' H_0 \tilde{\ell} \gtrsim c_1 \sqrt{n m^{(q)}}$. 
Similarly, combining these with the second claim of Lemma \ref{lem:v'H0v} and Weyl's inequality, we find that 
\[ 
\lambda_1(H) - \lambda_n(H) \lesssim  c_2 \sqrt{n m^{(q)}}, \qquad  \lambda_1(H) - \tilde{\ell}' H \tilde{\ell} \gtrsim c_1 \sqrt{n m^{(q)}}. 
\]
Inserting these into (\ref{ifpcatoshow})  gives the claim. 

%%%%%%%%%%
%%%%%%%%%%
%%%%%%%%%%
\section{Limits for signal recovery}    \label{sec:sig}  
\setcounter{equation}{0} 
In this section, we discuss limits for signal (support) recovery.  The results are intertwined with those for clustering (namely, Theorems \ref{thm:clu-a}-\ref{thm:clu-c} and Theorem \ref{thm:clu-d}), so we prove all of them together in the later part of the section. 
 
Compare two problems: signal recovery and clustering. One useful insight is that in the less sparse case, 
clustering is comparably easier than signal recovery, so we should estimate $\ell$ first and then use it to estimate $S(\mu)$; in the 
more sparse case, we should do the opposite.  

For the less sparse case, we have introduced two clustering methods, $\hell_*^{(sa)}$ and $\hell_*^{(if)}$, in Section \ref{subsec:4method}. They give rise to two signal recovery 
methods, $\hs_*^{(sa)}$ and $\hs_*^{(if)}$.
In detail, let $y_*^{(sa)} = n^{-1/2} X' \hell_*^{(sa)}$ and $y_*^{(if)} = n^{-1/2} X' \hell_*^{(if)}$, and let   $t_p^* = \sqrt{2 \log(p)}$ be the {\it universal threshold} \citep{DonohoJohnstone}.  Respectively, $\hs_*^{(sa)}$ and $\hs_*^{(if)}$ are defined by 
\[
\hs_{*}^{(sa)} = \{1 \leq j \leq p: |y_{*,j}^{(sa)}| \geq t_p^* \}, \qquad  \hs_{*}^{(if)} =  \{1 \leq j \leq p:   |y_{*,j}^{(if)}| \geq t_p^* \}. 
\]

For the more sparse case, we introduce two methods $\hs_N^{(sa)}$ and $\hs_{q}^{(if)}$; 
they are in fact the ones that give rise to the clustering methods $\hell_N^{(sa)}$ and $\hell_q^{(if)}$ we introduced in 
Section \ref{subsec:4method}. 
In detail,  recalling that $Q(j) = (2n)^{-1/2} (\|x_j\|^2 -n)$ is the column-wise $\chi^2$-statistics, 
\begin{equation} \label{optimization} 
\hs_N^{(sa)} = \margmax_{\{S: S \subset  \{1, 2, \ldots, p\}, |S| = N\}} \{ N^{-1/2}  \| \sum\nolimits_{j  \in S } x_j \|_1  \},  
\end{equation} 
and 
\[
\hs_q^{(if)} = \{1 \leq j  \leq p:  Q(j)  \geq \sqrt{2 q \log(p)}\}. 
\]
 
For any signal (support) recovery procedure $\hat{S}$,  we measure the performance by the normalized size of the difference of 
$\hat{S}$ and the true support 
\begin{equation} \label{hamm2} 
\hamm_p(\hs,  \alpha,\beta,\theta) =   (p \eps_p)^{-1}  E( |\hat{S} \Delta S(\mu)|),     
\end{equation} 
where $A\Delta B=(A\setminus B)\cup(B\setminus A)$ denotes the symmetric difference of two sets and 
the expectation is with respective to the randomness of $(\mu,\ell, Z)$.  
If we think $\hat{S}$ as an estimate of $\mu$, say, $\hat{\mu}$, and $E( |\hat{S} \Delta S(\mu)| )$ is actually the 
 Hamming distance between the two vectors $(\sgn(|\muhat(1)|), \ldots, \sgn(|\muhat(p)|))'$ and $(\sgn(\mu(1)), \ldots, \sgn(\mu(p)))'$. For this reason, we call that in (\ref{hamm2})   the (normalized) Hamming distance. %To simplify the notations,  we omit $(\alpha,\beta,\theta)$ in the brackets when where is no confusion. 

In Section~\ref{subsec:other-limits}, we have introduced the curves $\alpha = \eta^{sig}_\theta(\beta)$ and $\alpha = \tilde{\eta}^{sig}_\theta(\beta)$. The following theorem is proved in Section \ref{sec:proof}.  
%%%%%%%%%%
%%%%%%%%%%
\begin{thm} \label{thm:siga} 
({\it Statistical lower bound for signal recovery}). 
Fix $(\alpha, \beta, \theta) \in (0,1)^3$ and suppose $\alpha > \eta_{\theta}^{sig}(\beta)$.  
Consider the signal recovery problem for Models (\ref{model1})-(\ref{model2}) and (\ref{ARW1})-(\ref{ARW2}).  For any $\hs$ that is an estimate for the  support of $S$,  $\hamm_p(\hs, \alpha, \beta,\theta) \gtrsim 1$ as $p \goto \infty$. 
\end{thm} 
We also have the following theorems, which are proved below. 
%%%%%%%%%%
%%%%%%%%%%
\begin{thm} \label{thm:sigb} 
({\it Statistical upper bound for signal recovery}). 
Fix $(\alpha, \beta, \theta) \in (0,1)^3$ and suppose $\alpha  < \eta_{\theta}^{sig}(\beta)$. 
Consider the signal recovery problem for Models (\ref{model1})-(\ref{model2}) and (\ref{ARW1})-(\ref{ARW2}).    
As $p \goto \infty$,  
\begin{itemize} 
\item $\hamm_p(\hs_*^{(sa)},    \alpha,\beta,\theta) \goto 0$,  if $0 < \beta < (1 - \theta)/2$.  
\item $\hamm_p(\hs_N^{(sa)}, \alpha,\beta,\theta) \goto 0$,    if $(1 - \theta)/2 < \beta < 1$ and $N = \lceil p \eps_p\rceil$.   
\end{itemize} 
\end{thm}    
 
%%%%%%%%%%
%%%%%%%%%%
\begin{thm} \label{thm:sigc} 
({\it CTUB for signal recovery}). 
 Fix $(\alpha, \beta, \theta) \in (0,1)^3$ and suppose $\alpha  < \tilde{\eta}_{\theta}^{sig}(\beta)$. 
Consider the signal recovery problem for Models (\ref{model1})-(\ref{model2}) and (\ref{ARW1})-(\ref{ARW2}).    
As $p \goto \infty$,  
\begin{itemize} 
\item $\hamm_p(\hs_*^{(sa)},    \alpha,\beta,\theta) \goto 0$,  if $0 < \beta < (1 - \theta)/2$.  
\item $\hamm_p(\hs_*^{(if)}, \alpha,\beta,\theta) \goto 0$,    if $(1 - \theta)/2 < \beta < 1/2$. 
\item $\hamm_p(\hs_q^{(if)}, \alpha,\beta,\theta) \goto 0$,    if $(1 - \theta)/2 < \beta < 1/2$ and $q \geq  3$. 
\end{itemize} 
\end{thm} 

%%%%%%%%%%%
%%%%%%%%%%%
%%%%%%%%%%%
\subsection{Proofs of Theorems \ref{thm:club}-\ref{thm:clu-c}, \ref{thm:clue} and \ref{thm:sigb}-\ref{thm:sigc}}  
\label{subsec:UBpf1} 
 
We need two lemmas. 
The first one is on classical PCA, and it is needed for studying $\hell_*^{(if)}$ and $\hat{S}_*^{(if)}$.  
The second one is a large-deviation inequality for folded normal random variables and it is needed for 
studying the optimization problem in (\ref{optimization}). 
%%%%%%%%%%%%%%%
%%%%%%%%%%%%%%%
%%%%%%%%%%%%%%%
\begin{lemma}  \label{lemma:PCA} 
Fix $(\alpha, \beta, \theta) \in (0,1)^3$ such that $(1 - \theta)/2 < \beta < 1/2$ and $\alpha < \tilde{\eta}_{\theta}^{clu}(\beta)$. In Models (\ref{model1})-(\ref{model2}) and (\ref{ARW1})-(\ref{ARW2}), let $\lambda$ be the 
first eigenvalue of $XX'$ and 
$\xi$ be the corresponding eigenvector. 
There is a generic constant $\delta  = \delta(\alpha, \beta, \theta) > 0$ such that 
with probability $1 - O(p^{-3})$, 
\[
\min\{ \| \sqrt{n} \xi + \ell\|_{\infty},  \|\sqrt{n} \xi - \ell\|_{\infty} \}  < p^{-\delta}.
\]   
The claim continues to hold if we replace the model \eqref{model2} by (\ref{newARW}) for $A = I_n$ and  $B$ such that $\max\{\|B\|,  \|B^{-1}\| \} \leq L_p$. 
\end{lemma} 
%%%%%%%%%%%
%%%%%%%%%%%
%%%%%%%%%%%
%%%%%%%%%%%
\begin{lemma} \label{lemma:FN} 
({\it Large-deviation on Folded Normals}).   As $n \goto \infty$, for any $h > 0$ and $0  \leq x \leq \sqrt{n} / \log(n)$, and $n$ independent samples 
$z_i$ from $N(0,1)$, 
\[
P \bigl(  \bigl| \sum\nolimits_{i = 1}^n  \bigl( |z_i + h| - E[|z_i +h|] \bigr)  \bigr|  \geq \sqrt{n} x \bigr)  \leq 
2 \mathrm{exp}\bigl(-  (1 + o(1))  x^2/2  \bigr),   
\]
where $o(1) \goto 0$, uniformly for all $h > 0$ and $0 <  x \leq \sqrt{n} / \log(n)$. 
\end{lemma} 
%%%%%%%%
%%%%%%%%
%%%%%%%%
We now show all theorems about upper bound.  Since Theorems \ref{thm:club}-\ref{thm:clu-c} are special cases of Theorem \ref{thm:clue} with $B = I_p$,  it suffices to show Theorems \ref{thm:clue} and \ref{thm:sigb}-\ref{thm:sigc}. As there are four methods involved, it is more convenient to prove in a way by grouping  the items  
associated with each method together.  Fixing $(\alpha, \beta, \theta) \in (0,1)^3$ 
and viewing all statements in Theorems \ref{thm:clue} and Theorems \ref{thm:sigb}-\ref{thm:sigc}, what we need to show can be re-organized as follows (for the statements regarding $\hat{\ell}$, we need to prove that they hold for a general B where $\max\{\|B\|, \|B^{-1}\|\} \leq L_p$).   
\begin{itemize} 
\item {\it (a). Simple Aggregation}.  Consider the case  $0 < \beta < (1 - \theta)/2$. In this range, $\eta_{\theta}^{sig}(\beta) < \eta_{\theta}^{clu}(\beta)$. All we need to show is that if $\alpha < \eta_{\theta}^{clu}(\beta)$, then $\hell_*^{(sa)} = \ell$ with probability at least $1 - O(p^{-3})$, and that  if additionally $\alpha < \eta_{\theta}^{sig}(\beta)$, then $\hamm_p(\hs_*^{(sa)}, \alpha, \beta, \theta) \goto 0$. 
\item {\it (b). Sparse Aggregation}.  Consider the case $(1 - \theta)/2  < \beta < 1$. In this case, 
$\eta_{\theta}^{clu}(\beta) \leq \eta_{\theta}^{sig}(\beta)$.  Letting $N = \lceil p \eps_p\rceil$, all we need to show is that $\hamm_p(\hs_{N}^{(sa)}, \alpha, \beta, \theta) \goto 0$ if $\alpha < \eta_{\theta}^{sig}(\beta)$ and  $\hamm_p(\hell_N^{(sa)}, \alpha, \beta, \theta) \goto 0$  if additionally $\alpha < \eta_{\theta}^{clu}(\beta)$. 
\item {\it (c). Classical PCA}.  Consider the case $(1 - \theta)/2 < \beta < 1/2$ where only computationally tractable bounds are concerned and $\tilde{\eta}^{clu}_\theta(\beta)=\tilde{\eta}_{\theta}^{sig}(\beta)$. All we need to show is that if $\alpha < \tilde{\eta}_{\theta}^{clu}(\beta)$, then $\hell_*^{(if)} = \pm \ell$ with probability at least $1 - O(p^{-3})$ and that $\hamm_p(\hs_*^{(if)}, \alpha, \beta, \theta) \goto 0$. 
\item {\it (d). IF-PCA}.  Consider the case $1/2 < \beta < 1$ where only computationally tractable bounds are concerned and $\tilde{\eta}_{\theta}^{clu}(\beta) \leq \tilde{\eta}_{\theta}^{sig}(\beta)$. All we need to show is that if $\alpha < \tilde{\eta}_{\theta}^{sig}(\beta)$, then $\hat{S}_{q}^{(if)}  = S(\mu)$ with probability at least $1 - O(p^{-3})$; and if additionally $\alpha<\tilde{\eta}_\theta^{clu}(\beta)$, then $\hamm_p(\hell_q^{(if)}, \alpha, \beta, \theta) \goto 0$. 
\end{itemize}  
Consider (a). Note that $\hat{\ell}_*^{(sa)}=\mathrm{sgn}(\sum_{j=1}^p x_j)$ and $\sum_{j=1}^p x_j\sim N(\|\mu\|_0 \tau \ell, p I_n)$. By \eqref{(b)proof-0}, $\|\mu\|_0 \tau = p^{1 - \beta - \alpha} (1 + o(1))$. Hence, $\alpha<\eta_\theta^{clu}(\beta)$ implies $\|\mu\|_0 \tau \gg \sqrt{p}$, and it follows that $\hat{\ell}_*^{(sa)} = \ell$ with overwhelming probability. Once $\hat{\ell}_*^{(sa)}=\ell$, $y_*^{(sa)}=n^{-1/2}X'\ell\sim N(\sqrt{n}\mu, I_p)$. Noting that $\alpha<\eta^{sig}_{\theta}(\beta)$ implies $\sqrt{n}\tau\gg 1$, we have $\hamm_p(\hs_*^{(sa)}) \goto 0$ with overwhelming probability. Consider (c). The first claim is a direct result of Lemma~\ref{lemma:PCA}, and the second claim can be proved similarly as in (a). Consider (d). Recall that the column-wise test statistic $Q(j)$ is approximately distributed as $N(0,1)$ for useless features and $N(\sqrt{n/2}\tau^2, 1)$ for useful features. So $\tau \gg n^{-1/4}$ will assure successful signal recovery, which translates to $\alpha<\tilde{\eta}_\theta^{sig}(\beta)$. Once $\hat{S}_q^{(if)}=S(\mu)$, we restrict our attention to $X^{S(\mu)}$, the sub-matrix of $X$ restricted to the columns in $S(\mu)$, and the claim of Lemma~\ref{lemma:PCA} continues to hold by adapting the proof there (see the Appendix for details). So $\hamm_p(\hell^{(if)}_q)\goto 0$ with overwhelming probability.  
It remains to prove (b). 
 
%%%%%%%%%%%
%%%%%%%%%%%
%%%%%%%%%%%
We now show (b). Define $\muhat_N^{(sa)}$ such that $\hat{\mu}_N^{(sa)}(j)=\tau_p\cdot 1\{j\in \hat{S}_N^{(sa)}\}$. Write $\hat{S}_N^{(sa)}=\hat{S}$, $\hat{\mu}_N^{(sa)}=\muhat$, $\hell^{(sa)}_N=\hell$ and $s_p=p\eps_p$. With probability $1-O(p^{-3})$, 
\beq \label{(b)proof-0}
|\|\mu\|_0-s_p|\leq C\sqrt{s_p\log(p)}. 
\eeq
Since any event of probability $O(p^{-3})$ has a negligible effect to the Hamming distances, we always condition on a fixed realization $(\ell,\mu)$ that satisfy \eqref{(b)proof-0}; so the probabilities below are with respective to the randomness of $Z$. 
To show (b), all we need to show are 
\begin{itemize} 
\item (b1).  $\hamm_p(\hat{\ell}, \alpha, \beta, \theta)\goto 0$,  if $\alpha<\eta_\theta^{clu}(\beta)$. In this item, the matrix $B$ may be any matrix that satisfies $\max\{\|B\|, \|B^{-1}\|\}\leq L_p$.   
\item (b2).  $\hamm_p(\hat{S}, \alpha, \beta, \theta) \goto 0$, if  $\alpha < \eta_{\theta}^{sig}(\beta)$. In this item, $B = I_p$.    
\end{itemize}

Consider (b1) first. 
It suffices to show
\beq \label{(b)proof-9}
n^{-1}\langle\hell,\ell\rangle \goto 1. 
\eeq
For any realized $\mu$, we construct $\tilde{\mu}$ as follows:   
\begin{itemize} 
\item If $\|\mu\|_0 > N$, replace $\|\mu\|_0 - N$ nonzero entries by $0$. 
\item If $\|\mu\|_0 < N$, replace $N - \|\mu\|_0$  zero entries by $\tau_p$.
\end{itemize} 
Let $\tilde{S}$ be the support of $\tilde{\mu}$. Write $X \tilde{\mu} = \|B \tilde{\mu}\| \cdot  [Z (B \tilde{\mu} / \|B \tilde{\mu}\|) + \langle\mu, \tilde{\mu}/\|B \tilde{\mu}\|\rangle \ell ]$, where $Z (B\tilde{\mu} / \|B \tilde{\mu}\|) \sim N(0, I_n)$ and $\langle\mu, \tilde{\mu}/\|B \tilde{\mu}\| \rangle \gtrsim \|B\|^{-1}\tau_p\sqrt{s_p}= L_p  p^{(1 - \beta - 2\alpha)/2}$, with $(1 - \beta - 2\alpha) > 0$ in our range of interest. According to Mills' ratio \cite{Wellner86}, with probability $1 - O(p^{-3})$, the absolute value of standard normal variable is bounded by $\sqrt{6 \log(p)}$, which is less than $L_p p^{(1 - \beta - 2\alpha)/2}$  when $p \goto \infty$. It follows that with probability at least $1   - O(p^{-3})$, $\sgn(X \tilde{\mu})=\ell$. Furthermore,  $\ell' X \tilde{\mu} = \|X \tilde{\mu}\|_1=\tau_p\|\sum_{j\in \tilde{S}}x_j\|_1$.  
Since that $\hat{\ell}' X \hat{\mu} =  \|X\hat{\mu}\|_1=\tau_p\|\sum_{j\in \hat{S}} x_j\|_1$ and that $\hat{S}$ solves the  optimization  problem \eqref{optimization},  
\beq  \label{(b)proof-8}
\hell' X\muhat\geq \ell' X\tilde{\mu}. 
\eeq
Write $\hell' X\muhat=\langle\hell, \ell\rangle \langle \mu, \muhat\rangle + \hell'ZB\muhat$. We aim to obtain  an upper bound for $|\hell'ZB\muhat|$ (an upper bound for $|\ell'ZB\tilde{\mu}|$ can be obtained similarly). Denote by $(ZB)^{\hat{S}}$ the sub-matrix of $ZB$ containing columns in $\hat{S}$. Then $|\hell'ZB\muhat|\leq \sqrt{n}\|(ZB)^{\hat{S}}\|\|\mu\|\leq \sqrt{ns_p}\tau_p \|(ZB)^{\hat{S}}\|$, where $\|(ZB)^{\hat{S}}\|\leq \|B\|\|Z^{\hat{S}}\|\leq L_p \max_{|S|=N}\|Z^S\|$. By classical RMT \cite{Vershynin}, $\max_{|S|=N}\|Z^S\|\leq L_p\max\{\sqrt{n}, \sqrt{s_p}\}$ with probability at least $1-O(p^{-3})$. 
Inserting them into \eqref{(b)proof-8} gives
\begin{equation} \label{(b1)1} 
\langle \hell, \ell\rangle \langle \mu, \muhat\rangle \geq n \langle\mu, \tilde{\mu}\rangle -  L_p \sqrt{ns_p}\tau_p(\sqrt{n}+\sqrt{s_p}). 
\end{equation} 
First, $\langle \mu,\hat{\mu} \rangle\leq \max\{\|\mu\|_0, N\} \tau_p^2\sim s_p\tau_p^2$. Second, by \eqref{(b)proof-0} and the definition of $\tilde{\mu}$, $\langle \mu, \tilde{\mu}\rangle =s_p\tau_p^{2} (1+o(1))$. Inserting these into (\ref{(b1)1}) gives 
$n^{-1}\langle \hell,\ell \rangle \geq 1 - L_p (\sqrt{n}+\sqrt{s_p})/(\tau_p\sqrt{ns_p})$.  
When $\alpha<\eta_\theta^{clu}(\beta)$, the second term on the right hand side is $\leq p^{-\delta}$ for some $\delta=\delta(\alpha,\beta,\theta)>0$, and \eqref{(b)proof-9} follows. 

We now consider (b2). Let $\tilde{\mu}$ and $\tilde{S}$ be the same as above. 
Due to \eqref{(b)proof-0}, $|\tilde{S}\cap S(\mu)|\geq |S(\mu)|(1+o(1))$. It 
suffices to show that 
\beq  \label{(b)proof-1}
|\hat{S}\cap S(\mu)| \geq | \tilde{S}\cap S(\mu)|-o(s_p). 
\eeq
Since $|\tilde{S}|=|\hat{S}|=N$ and that $\hat{S}$ solves the optimization \eqref{optimization},  
\beq  \label{(b)proof-2}
G(\hat{S})\equiv N^{-1/2}\sum_{i=1}^n |\sum_{j\in \hat{S}} X_i(j)|\geq N^{-1/2}\sum_{i=1}^n |\sum_{j\in \tilde{S}} X_i(j)|\equiv G(\tilde{S}).
\eeq
For any $S\subset\{1,\cdots,p\}$ such that $|S|=N$, we define $w_i(S)=N^{-1/2} \sum_{j\in S}Z_i(j)$ and  $h(S)=N^{-1/2}|S\cap S(\mu)|\tau_p$.  
It follows that
$G(S) \overset{(d)}{=}\sum_{i=1}^n |w_i(S)+h(S)|$, where $w_i(S)\overset{iid}{\sim} N(0,1)$, $1\leq i\leq n$. 
For any $h>0$, we define the function $u(h)=E_{X\sim N(0,1)}(|X+h|)$. Let $E_p$ be the event that $\{\max_{S\subset\{1,\cdots,p\}, |S|= N}|G(S)-u(h(S))|\leq \sqrt{6N\log(p)/n}\}$. 
By Lemma~\ref{lemma:FN} and the fact that 
there are no more than $p^N$ such $S$,  $P(E_p^c)=O(p^{-3})$; so those realizations $Z$ in $E_p^c$ has a negligible effect. Combining it with \eqref{(b)proof-2} gives
\beq \label{(b)proof-4}
u(h(\hat{S})) \geq u(h(\tilde{S})) - L_p\sqrt{s_p \log(p)/n}.
\eeq
The following lemma is proved in the Appendix.  
%%%%%%%%%%
%%%%%%%%%%
\begin{lemma} \label{lem:u-prop}
%$u(h)$ is strictly convex and monotone increasing for $h\in(0,\infty)$, $u(h)\to \sqrt{2/\pi}$ as $u\to 0$ and $u(h)/h \to 1$ as $u\to \infty$. Moreover, 
There exists a constant $C>0$ such that for any $0<h_1<h_2$, $u(h_2)-u(h_1)\geq C \min\{(h_2-h_1), (h_2-h_1)^2\}$. 
\end{lemma}
%%%%%%%%%%
\noindent
Since $h(S)\leq h(\tilde{S})$ for any $S$ with $|S|=N$, by Lemma~\ref{lem:u-prop},
\beq \label{(b)proof-6}
u(h(\tilde{S})) \geq u(h(\hat{S})) - C \min\{ h(\tilde{S})-h(\hat{S}), \; [h(\tilde{S})-h(\hat{S})]^2 \} 
\eeq
We combine \eqref{(b)proof-4}-\eqref{(b)proof-6}. It yields that
\beq \label{(b)proof-7}
0\leq  \frac{h(\tilde{S})-h(\hat{S})}{\sqrt{s_p}\tau_p} \leq L_p \tau_p^{-1}\max\left\{(\log(p)/n)^{1/2},\;\; (\log(p)/ns_p)^{1/4}\right\}. 
\eeq
The assumption $\alpha<\eta^{sig}_\theta(\beta)$ implies $\tau_p\leq p^{-\delta}\min\{n^{-1/2}, (ns_p)^{-1/4}\}$. So the right hand side of \eqref{(b)proof-7} is $o(1)$. Then \eqref{(b)proof-1} follows. 

%%%%%%%%%%%
%%%%%%%%%%%
%%%%%%%%%%%
\section{Limits for hypothesis testing} \label{sec:hyp}   
\setcounter{equation}{0} 
The goal for (global) hypothesis testing is to test a null hypothesis 
\begin{equation} \label{null} 
H_0^{(p)}: \qquad X_i \stackrel{iid}{\sim} N(0, I_p), \qquad 1 \leq i \leq n, 
\end{equation} 
against a specific alternative in the complement of the null, 
\begin{equation} \label{alt} 
H_1^{(p)}: \;\;  \mbox{$X_i$'s are generated from Models (\ref{model1})-(\ref{model2}) and (\ref{ARW1})-(\ref{ARW2})}. 
\end{equation} 
We consider three different tests. 

The first test $\hat{T}_*^{(sa)}$ is connected to the idea of simple aggregation. Recall that $\bar{x}$ is the average of all columns. 
The idea is to test whether $E[\bar{x}] = 0$ or not using the classical $\chi^2$. This test rejects $H_0^{(p)}$ if and only if 
\[
(2n)^{-1/2} [p  \|\bar{x}\|^2  - n] \geq 2\sqrt{2 \log(p)}.  
\]
The second test $\hat{T}_N^{(sa)}$ is connected to sparse aggregation. Let $\hs_N^{(sa)} $ be as in \eqref{optimization}. This test rejects $H_0^{(p)}$ if and only if 
\[
 N^{-1/2} \| \sum_{j \in \hs^{(sa)}_N}x_j \|_1    \geq     \sqrt{2/\pi} n +  \sqrt{2n(N + 2) \log(p)}. 
\]
The third test $\hat{T}^{(hc)}$ is connected to the Higher Criticism  in Donoho and Jin \cite{DJ04}. 
Recalling that $Q(j)=(2n)^{-1/2}(\|x_j\|^2-n)$ are the column-wise $\chi^2$-tests,  the idea is to test whether some of the $Q(j)$'s have non-zero means. 
\begin{itemize}
\item For $1 \leq j\leq p$, obtain a $P$-value $\pi_j = P\{ (2n)^{-1/2}[\chi^2_n(0)-n] \geq Q(j)\}$.
\item Sort  the $P$-values in the ascending order:
$\pi_{(1)} < \pi_{(2)} < \ldots < \pi_{(p)}$.
\item Compute the Higher Criticism statistic $HC_p^*   =  \max_{\{ 1 \leq i \leq p/2 \}}  HC_{p, i}$, where 
$HC_{p, i} \equiv  \sqrt{p}  [(i/p)  - \pi_{(i)}]/ [\pi_{(i)} (1 - \pi_{(i)})]^{1/2}$. 
\end{itemize}
The test rejects $H_0^{(p)}$ if and only if
$HC_p^* \geq  2\sqrt{2\log\log(p)}$.

The test $\hat{T}_{N}^{(sa)}$ is similar to a test in \cite{Ery}, which is designed for the case that there is (unknown) dependence among features and so the test is more complicated than ours. The other two tests are newly proposed.

For any testing procedure $\hat{T}$ that tests $H_1^{(p)}$ against $H_0^{(p)}$,  we measure the performance by the sum of Type I and Type II errors: 
\begin{equation} \label{df} 
Err(\hat{T},\alpha,\beta,\theta) = P_{H_0^{(p)}}(\mbox{$\hat{T}$ rejects $H_0^{(p)}$})  + P_{H_1^{(p)}}(\mbox{$\hat{T}$ accepts $H_0^{(p)}$}),
\end{equation} 
where the probabilities are with respective to the randomness of $(\ell,\mu,Z)$. 

In Section~\ref{subsec:other-limits}, we have introduced two curves $\eta_\theta^{hyp}(\beta)$ and $\tilde{\eta}_\theta^{hyp}(\beta)$. 
The following theorem is proved in Section \ref{sec:proof}. 
%%%%%%%%
%%%%%%%%
%%%%%%%%
%%%%%%%%
\begin{thm} \label{thm:hypa} 
({\it Statistical lower bound for hypothesis testing}). 
Fix $(\alpha, \beta, \theta) \in (0,1)^3$ with $\alpha > \eta_{\theta}^{hyp}(\beta)$.  Consider the testing problem (\ref{null})-(\ref{alt}) for 
 Models (\ref{model1})-(\ref{model2}) and (\ref{ARW1})-(\ref{ARW2}).  For any  test $\hat{T}$, $Err(\hat{T},\alpha,\beta,\theta)  \gtrsim 1$ as $p \goto \infty$. 
\end{thm} 
Consider the upper bound. By the definitions (see (\ref{hypphasefunction})), when $\alpha < \eta_{\theta}^{hyp}(\beta)$, we have either $\alpha < \eta_{\theta}^{hyp, 1}(\beta)$ or $\alpha < \rho_{\theta}^{hyp, 2}(\beta)$, or both. 
%%%%%%%%%%
%%%%%%%%%%
\begin{thm} \label{thm:hypb} 
({\it Statistical upper bound for hypothesis testing}). 
Fix $(\alpha, \beta, \theta) \in (0,1)^3$ such that $\alpha <  \eta_{\theta}^{hyp}(\beta)$.  Consider the testing problem (\ref{null})-(\ref{alt}) for 
 Models (\ref{model1})-(\ref{model2}) and (\ref{ARW1})-(\ref{ARW2}).
As $p \goto \infty$, 
\begin{itemize} 
\item $Err(\hat{T}_*^{(sa)},\alpha, \beta,\theta)   \goto 0$ if $\alpha < \eta_{\theta}^{hyp,1}(\beta)$.  
\item $Err(\hat{T}_N^{(sa)}, \alpha,\beta,\theta) \goto 0$ if  $\alpha < \eta_{\theta}^{hyp,2}(\beta)$  and we take $N  = \lceil p \eps_p\rceil$.  
\end{itemize} 
\end{thm} 

%%%%%%%%%%
%%%%%%%%%%
\begin{thm} \label{thm:hypc} 
({\it CTUB for hypothesis testing}).  
Fix $(\alpha, \beta, \theta) \in (0,1)^3$ such that $\alpha <  \tilde{\eta}_{\theta}^{hyp}(\beta)$.  Consider the testing problem (\ref{null})-(\ref{alt}) for 
 Models (\ref{model1})-(\ref{model2}) and (\ref{ARW1})-(\ref{ARW2}).
As $p \goto \infty$, 
\begin{itemize} 
\item $Err(\hat{T}_*^{(sa)}, \alpha,\beta,\theta)   \goto 0$ if $0 < \beta < 1/2$.   
\item $Err(\hat{T}^{(hc)}, \alpha,\beta,\theta) \goto 0$ if  $1/2 < \beta < 1$.  
\end{itemize} 
\end{thm} 

\begin{comment}
\begin{figure}[t!]
\begin{center}
\includegraphics[height = 1.9 in]{hyp.pdf}
\includegraphics[height = 1.9 in]{clusterlimits_sign.pdf}
\end{center}
\caption{Left: statistical limits for hypothesis testing problem when $\theta = 4/9$ (red), $2/3$ (blue), and $8/9$ (black). Right: the statistical limits (red) and the CTUB (green) for clustering when the signals have different signs.} 
\label{fig:hypphase}
\end{figure}
\end{comment}
%%%%%%%%%%%%%%
%%%%%%%%%%%%%%
%%%%%%%%%%%%%%
%%%%%%%%%%%%%%
\subsection{Proofs of Theorems \ref{thm:hypb}-\ref{thm:hypc}}  \label{subsec:UBpf2}  
Similarly, as three tests are involved,   it is more convenient to prove the results  in a way by grouping items associated with each test separately. Fixing $(\alpha, \beta, \theta) \in (0, 1)^3$ and viewing the two theorems, the following is what we need to show. 
\begin{itemize} 
\item ({\it Simple Aggregation}).  When $\alpha < \eta_{\theta}^{hyp,1}(\beta)$, $Err(\hat{T}_*^{(sa)}, \alpha,\beta,\theta) \goto 0$. 
\item ({\it Sparse Aggregation}).  When $\alpha < \eta_{\theta}^{hyp,2}(\beta)$, $Err(\hat{T}_N^{(sa)}, \alpha,\beta,\theta) \goto 0$ if we take $N = \lceil p \eps_p \rceil$. 
\item ({\it HC}).  When $1/2 < \beta < 1$ and $\alpha < \theta/4$, $Err(\hat{T}^{(hc)}, \alpha,\beta,\theta) \goto 0$. 
\end{itemize} 

In the above, (c) is an easy extension of \cite{DJ04}, so we omit its proof. Below, we prove (a) and (b). 
Consider (a). $\hat{T}_*^{(sa)}$ is defined through $\bar{x}$, where $\bar{x}\sim N(p^{-1}\|\mu\|_0\tau, p^{-1}I_n)$. So the claim follows directly from the tail probability of chi-square distributions.
Consider (b). Under $H_0^{(p)}$, for each fixed $S$ with $|S|=N$, we can write $N^{-1/2}\|\sum_{j\in S}x_j\|_1=\sum_{i=1}^n |w_i|$, where $w_i$'s are {\it iid} standard normal variables. Since $E(|w_i|)=\sqrt{2/\pi}$, by Lemma~\ref{lemma:FN}, $\hat{T}_{N}^{(sa)}\leq\sqrt{(2/\pi)} n + \sqrt{2n(N + 2) \log(p)}$ with probability $1 - O(p^{-2})$. 
On the other hand, it is seen that 
\[
\hat{T}_{N}^{(sa)} = \max\nolimits_{\ell\in\{\pm 1\}^n, \mu\in\{0,1\}^p, \|\mu\|_0 = N} \ell' X (\mu/\|\mu\|). 
\] 
Under $H_1^{(p)}$, let $\tilde{\mu}$ be defined in the same way as Section~\ref{subsec:UBpf1} and so $\hat{T}_{N}^{(sa)} \geq \ell' X\tilde{\mu}/\|\tilde{\mu}\|=n\langle \mu,\tilde{\mu}\rangle/\|\tilde{\mu}\| + \ell'Z\tilde{\mu}/\|\tilde{\mu}\|$. Since $|S(\mu)| \sim p\eps_p$ with probability $1 - O(p^{-2})$, $\langle \mu,\tilde{\mu} \rangle/\|\tilde{\mu}\|\geq \|\mu\|(1+o(1))$. Moreover, $\|\ell'Z\tilde{\mu}\|\leq C\sqrt{n}\|\tilde{\mu}\|(\sqrt{n}+\sqrt{N})$ with probability $1 - O(p^{-3})$, by classical RMT \cite{Vershynin}. Combining the above gives
$\hat{T}_N^{(sa)} \gtrsim n\|\mu\| - C\sqrt{n}(\sqrt{n}+\sqrt{p\eps_p}) \geq n\|\mu\|/2$, 
where the last inequality is because $\alpha<\eta_\theta^{hyp,2}(\beta)$ implies $\tau_p\gg \max\{n^{-1/2}, s_p^{-1/2}\}$. 
Therefore, $\hat{T}_N^{(sa)}\gtrsim n \tau_p \sqrt{N}/2 \gg \max\{n, \sqrt{n N \log(p)}\}$, and the claim follows. 

%Consider (c). In this range, feature $j$ is selected with probability $O(p^{-2})$ when $\mu(j) = 0$ and with probability $1 - O(p^{-2})$ when $\mu(j) \neq 0$, so $\hat{S} = S(\mu)$ with probability at least $1 - O(p^{-1})$.  
%Now, over the event $\{\hat{S} = S(\mu)\}$ (noting that $|S(\mu)|\sim p\eps_p$),  by elementary RMT (e.g., Lemma~\ref{lem:vershynin}),  under $H_0^{(p)}$, there is $\lambda^{(q)} \leq (\sqrt{n} + \sqrt{|S(\mu)|} + 2\sqrt{\log(p)})^2$ with probability $1 - O(p^{-2})$; under $H_1^{(p)}$, $\lambda^{(q)} \geq n \|\mu\|^2-L_p\max\{|S(\mu)|, n\} \gtrsim np\eps_p\tau_p^2 -L_p\max\{|S(\mu)|, n\} \gg (\sqrt{n} + \sqrt{|S(\mu)|} + 2\sqrt{\log(p)})^2$ , and  the third claim follows. 
%%%%%%%%%%
%%%%%%%%%%
%%%%%%%%%%
%%%%%%%%%%
%%%%%%%%%%
%%%%%%%%%
%%%%%%%%%
%%%%%%%%%
%%%%%%%%
%%%%%%%%
%%%%%%%%
\section{Proofs of Theorems  \ref{thm:clu-a}, \ref{thm:clu-d},  \ref{thm:siga}, and  \ref{thm:hypa} (lower bounds)}   \label{sec:proof}  
\setcounter{equation}{0} 
%%%%%%%%%%%%
%%%%%%%%%%%%
%%%%%%%%%%%%
\subsection{Proof of Theorem~\ref{thm:clu-a}}     
\label{subsec:cluLB} 
For each $1 \leq i  \leq n$, 
consider the testing of two hypotheses, $H_{-1}^{(i)}:  \ell_i = -1$ versus $H_1^{(i)}:     \ell_i = 1$. 
Let $f_{\pm}^{(i)}$ be the joint density of $X$ under $H_{\pm 1}^{(i)}$, respectively. Since $\ell_i  = \pm 1$ with equal probabilities,  it follows from the connection between $L^1$-distance and the sum of Type I and Type II testing errors  \cite{Aad} that for any clustering procedure $\hat{\ell}$,  
$P(\hat{\ell}_i  \neq \ell_i)   \geq  1 - \|f_{-}^{(i)} - f_{+}^{(i)}\|_1$.  
Comparing this with the desired claim, it suffices to show that for all $1 \leq i \leq n$,  
\begin{equation} \label{cluLB1} 
\|f_{-}^{(i)} - f_{+}^{(i)}\|_1   = o(1),   \qquad  \mbox{where $o(1) \goto 0$ and does not depend on $i$}. 
\end{equation} 

We now show (\ref{cluLB1}) for every fixed $1 \leq i \leq n$. For short, we drop the superscript ``$(i)$"  in $f_{\pm}^{(i)}$ and $H_{\pm 1}^{(i)}$.   Recall that $X = \ell \mu' + Z$.   Denote $\tilde{\ell} = \ell - \ell_i e_i$, where $e_i$ is the $i$-th standard basis vector of $R^n$; note that $\tilde{\ell}_i  = 0$.   
By basic calculus and Fubini's theorem, 
\begin{align*} 
\|f_{-} - f_{+}\|_1 & =  E \bigl[\bigl| \int \sinh(X_i' \mu) e^{-\|\mu\|^2/2}  e^{ \tilde{\ell}' X \mu - (n-1) \|\mu\|^2/2} d F(\mu) d F(\tilde{\ell}) \bigr|  \bigr]     \\
&  \leq  E \bigl[\int \bigl| \int \sinh(X_i' \mu) e^{-\|\mu\|^2/2}  e^{ \tilde{\ell}' X \mu - (n-1) \|\mu\|^2/2} d F(\mu)  \bigr|  d F(\tilde{\ell})  \bigr]     \\
& =   \int   E \bigl[\bigl| \int \sinh(X_i' \mu) e^{-\|\mu\|^2/2}  e^{ \tilde{\ell}' X \mu - (n-1) \|\mu\|^2/2} d F(\mu)  \bigr| \bigr]    d F(\tilde{\ell}),  
\end{align*} 
where $E$ denotes the expectation under the law of $X = Z$. 
Seemingly, to show (\ref{cluLB1}), it suffices to show that for every realization of $\tilde{\ell}$, 
\begin{equation} \label{cluLB2}
E \bigl[\bigl| \int \sinh(X_i' \mu) e^{-\|\mu\|^2/2}  e^{ \tilde{\ell}' X \mu - (n-1) \|\mu\|^2/2} d F(\mu)  \bigr| \bigr]  = o(1);  
\end{equation}  
note that the left hand side does not depend on $i$ and $\tilde{\ell}$. 
We now show (\ref{cluLB2}) for the cases of $\beta > (1 - \theta)$ and $\beta < (1 - \theta)$, separately.

Consider the case $\beta < (1 - \theta)$ first. Introduce $V = (n-1)^{-1/2} X' \tilde{\ell}$; note that $V \sim N((n-1)^{1/2} \mu, I_p)$.  Let $g_{-}^{(i)}$,  $g_{+}^{(i)}$, and $g_0^{(i)}$ be the joint densities of $(X_i, V)$ for the cases of 
$X_i  = - \mu + z$, $X_i = \mu + z$, and $X_i  = z$, where  $z \sim N(0, I_p)$ and is independent of $\mu$   
(in all three cases, $V = (n-1)^{1/2} \mu + \tilde{z}$ where $\tilde{z}$ is independent of $(\mu, z)$).  
By the triangle inequality and symmetry, 
$\|g_{-}^{(i)} - g^{(i)}_{+}\|_1 \leq  \|g_{-}^{(i)} - g^{(i)}_0\|_1 + \|g_{+}^{(i)} - g^{(i)}_0\|_1  =  2  \|g_{+}^{(i)} - g_0^{(i)}\|_1$.  

We recognize that the left hand side of (\ref{cluLB2}) is nothing else but  $\|g_{-}^{(i)} - g_{+}^{(i)}\|_1$. Combining these, to show (\ref{cluLB2}),  it is sufficient to show 
\begin{equation} \label{cluLB3}
 \|g_{+}^{(i)} - g_0^{(i)}\|_1 = o(1). 
\end{equation} 

Now,  denote by $A(f, g)$ the Hellinger affinity for any two densities $f$ and $g$.   
Denote $h_p(V(j)) = \eps_p e^{\sqrt{n-1} \tau_p V(j) - (n- 1) \tau_p^2/2} /[1 - \eps_p + \eps_p e^{\sqrt{n-1} \tau_p V(j) - (n-1)\tau_p^2 /2}]$.  By definitions and direct calculations,  $A(g_{+}^{(i)},    g_{0}^{(i)})$ 
equals to
\begin{align*} \label{cluLB4} 
 \Pi_{j = 1}^p  E\big\{ \big[1 + h_p(V(j)) (e^{\tau_p X_i(j) - \tau_p^2/2} -1)\big]^{1/2}\big\} 
=  \bigl(E\big\{ \big[1 + h_p(V(1)) (e^{\tau_p X_i(1) - \tau_p^2/2} -1)\big]^{1/2}\big\} \bigr)^p. 
\end{align*} 
Write for short $u = X_i(1)$ and $w = V(1)$.  According to \cite[Page 221]{Aad}, for any probability densities $f$ and $g$, $ \|f-g\|_1\leq 2\sqrt{2-2A(f,g)}$. Combining this with the expression of $A(g_+^{(i)}, g_0^{(i)})$, to show (\ref{cluLB3}), it suffices to show  
\begin{equation} \label{cluLB5} 
E[ \bigl(1 + h_p(w)   [e^{\tau_p u - \tau_p^2/2} -1]\bigr)^{1/2}] = 1 + o(p^{-1}). 
\end{equation}
Note that for any $x > -1$,  $|\sqrt{1 + x} - 1 - x/2|  \leq C x^2$,  
\begin{eqnarray} \label{cluLB6a} 
& \bigl| E[ \bigl(1 + h_p(w)   [e^{\tau_p u - \tau_p^2/2} -1]\bigr)^{1/2}] - E[ 1 +   \frac{h_p(w)}{2}  (e^{\tau_p u - \tau_p^2/2} - 1)] \bigr| \nonumber\\ 
& \leq  C E[  h_p^2(w)  (e^{\tau_p u - \tau_p^2/2} -1)^2]. 
\end{eqnarray} 
On one hand, due to the independence between $w$ and $u$ and the fact that $E[e^{\tau_p u - \tau_p^2/2}] = 1$,   
we have $E\bigl[h_p(w) [e^{\tau_p u - \tau_p^2/2} - 1] \bigr]  = 0$ and  
$E[ h_p^2(w)  (e^{\tau_p u - \tau_p^2/2} -1)^2  ]  =  E [h_p^2(w)]  E  [(e^{\tau_p u - \tau_p^2/2} -1)^2]$. 
On the other hand, since $h_p(w) \leq \eps_p e^{\sqrt{n-1}  \tau_p w - (n-1) \tau_p^2/2}$, by direct calculations there is 
$E[h_p^2(w)] \leq  \eps_p^2 e^{(n-1) \tau_p^2}$,  and $E  [(e^{\tau_p u - \tau_p^2/2} -1)^2] =  e^{\tau_p^2} -1$.  
Inserting these into (\ref{cluLB6a}) and invoking $\eps_p = p^{-\beta}$, $\tau_p = p^{-\alpha}$, and $n = p^{\theta}$, 
%%%%%%%%%%
%%%%%%%%%%
%%%%%%%%%%
\begin{equation*} \label{cluLB6b} 
\bigl| E[ \bigl(1 + h_p(w)   [e^{\tau_p u - \tau_p^2/2} -1]\bigr)^{1/2}] - 1  \bigr| \leq C \eps_p^2 (e^{\tau_p^2} -1)  e^{(n-1) \tau_p^2} \leq C p^{-2 \beta - 2 \alpha} e^{ p^{\theta  -  2 \alpha}}. 
\end{equation*} 
By the assumptions of $\alpha > \eta_{\theta}^{clu}(\beta)$ and $\beta < (1 - \theta)$, we have  $2(\beta + \alpha) > 1$ and $\theta < 2\alpha$, and  (\ref{cluLB5}) follows. 

We now consider the case of  $\beta > (1 - \theta)$.  
In this case, similarly, by basic algebra and Fubini's theorem, the left hand side of (\ref{cluLB2}) is no greater than 
\begin{align} 
   & E [ \int |\sinh(X_i' \mu)|  e^{-\|\mu\|^2/2}  e^{ \tilde{\ell}' X \mu - (n-1) \|\mu\|^2/2}   d F(\mu)\bigr] \nonumber   \\
= &  \int E \bigl[ | \sinh(X_i' \mu)|   e^{-\|\mu\|^2/2} e^{\tilde{\ell}' X \mu - (n-1)\|\mu\|^2/2} \bigr] d F(\mu) \nonumber   \\
 = & \int E|\sinh(X_i' \mu)| e^{- \|\mu\|^2/2}  d F(\mu),   \label{cluLB2A} 
\end{align} 
where in the last step we have used the independence between $X_i$ and $\{X_k: k \neq i, 1 \leq k  \leq n\}$, and that 
$E [e^{\tilde{\ell}' X \mu - (n-1)\|\mu\|^2/2}]  = 1$.   
Finally, let $A_p$ be the event of $\{\mu:   \|\mu\|_0 / (p \eps_p)  \leq 2\}$, and write 
\begin{equation} \label{cluLB30}
\int E|\sinh(X_i' \mu)| e^{- \|\mu\|^2/2}  d F(\mu)  = I  + II, 
\end{equation} 
where 
$I = \int \bigl( E|\sinh(X_i' \mu)| e^{- \|\mu\|^2/2} \cdot 1_{A_p} \bigr)  d F(\mu)$, and   $II =  \int \bigl(E|\sinh(X_i' \mu)| e^{- \|\mu\|^2/2} \cdot  1_{A_p^c} \bigr)  d F(\mu)$. 

By Cauchy-Schwarz inequality,  $(E|\sinh(X_i' \mu)|)^2  \leq  E[(\sinh(X_i' \mu))^2]   =   (e^{2 \|\mu\|^2} -1)/2 $  for any realized $\mu$ in $A_p$. Combining this with basic algebra,  it follows that 
$I \leq \int ( \sqrt{\sinh(\|\mu\|^2)} \cdot 1_{A_p})  d F(\mu) \leq  \sqrt{\sinh(2 p \eps_p    \tau_p^2)}$, 
where in the last step, we have used the fact that over the event $A_p$, $\|\mu\|^2 \leq 2 p \eps_p \tau_p^2$.   
By our assumption of $\tau_p = p^{-\alpha}$, $\eps_p = p^{-\beta}$ and  $\alpha >  \eta_{\theta}^{clu}(\beta) =  (1 - \beta)/2$,   $p \eps_p \tau_p^2  = o(1)$. Combining these gives 
\begin{equation} \label{cluLB3a} 
|I| \leq o(1).
\end{equation} 
At the same time,   since $|\sinh(x)| \leq \cosh(x)$ for any $x$, 
\begin{equation} \label{cluLB3b} 
 II \leq \int \bigl( E \cosh(X_i' \mu) e^{-\|\mu\|^2/2}  \cdot 1_{A_p^c} \bigr) d F(\mu ) = P(A_p^c);  
\end{equation} 
note that $P(A_p^c) = o(1)$. We insert (\ref{cluLB3a})-(\ref{cluLB3b}) into (\ref{cluLB30}), and find that $\int E|\sinh(X_i' \mu)| e^{- \|\mu\|^2/2}  d F(\mu)  = o(1)$. Then (\ref{cluLB2}) follows from (\ref{cluLB2A}). \qed

%%%%%%%%%%%%%%
%%%%%%%%%%%%%%
\subsection{Proof of Theorem \ref{thm:clu-d}} 
Recall that $Z$ has {\it iid} entries from $N(0,1)$. By elementary statistics\footnote{Note that $ZB$ has the same distribution as $\tilde{c}_p \tilde{Z} + \tilde{W}$, where $\tilde{Z}$ has $iid$ normal entries, $\tilde{c}_p$ is half of the minimum eigenvalue of $BB'$, and the columns of $\tilde{W}$ follow $N(0, BB' - \tilde{c}_p I_p)$ distribution. Similar analysis for $A(\tilde{c}_p \tilde{Z} + \tilde{W})$ gives the result.} and conditions on $A$ and $B$, there is a non-stochastic term $c_p$ such that (a) $c_p^{-1} \leq L_p$, (b) there is a random matrix $W \in R^{n,p}$ such that $c_p Z + W$ has the same distribution of $A Z B$ ($W$ is independent of $(\ell, \mu, Z)$). 
Compare two experiments 
\[
\mbox{Experiment 1}.  \;\;      X = \ell \mu'  + c_p Z, \qquad   \mbox{Experiment 2}.  \;\;       X  = \ell \mu' + c_p Z + W.  
\]
Fixing  $1 \leq i \leq n$, consider the testing of two hypotheses,
$H_{-1}^{(i)}:    \ell_i = -1$ versus $H_1^{(i)}:  \ell_i = 1$. 
Let $f_{\pm}^{(i)}$ be the joint density of $X$ under  $H_{\pm}^{(i)}$, respectively, for Experiment 1, and let 
 $g_{\pm}^{(i)}$ be the joint density of $X$ under $H_{\pm}^{(i)}$, respectively, for Experiment 2. 
By Neyman-Pearson's fundamental lemma on  testing \citep{Aad},   for any clustering procedure $\hell$,  tight lower bounds for  $P(\hell_i \neq \ell_i)$ (expected Hamming error at location $i$) associated with the two experiments are  
$1 - \|f_{+}^{(i)} - f_{-}^{(i)} \|_1$ and  $1 - \|g_{+}^{(i)} - g_{-}^{(i)}\|_1$,  
respectively,   where $\|f-g\|_1$ denotes the $L^1$-distance  between two densities $f$ and $g$.  
Le Cam's idea can be solidified  as follows: 
\begin{thm} \label{thm:LeCam} 
({\it Monotonicity of $L^1$-distance}). 
$\|g_{+}^{(i)} - g_{-}^{(i)}\|_1 \leq  \|f_{+}^{(i)} -f_{-}^{(i)}\|_1$. 
\end{thm}
Using this, Theorem \ref{thm:clu-d} follows directly from the proof of Theorem \ref{thm:clu-a}. 

It remains to show Theorem \ref{thm:LeCam}. 
Without loss of generality, we assume $i = 1$, and drop the superscripts  in $g_{\pm}^{(i)}$ and $f_{\pm}^{(i)}$ for simplicity. 
Let $a \in R^{n-1}$ be the vector such that $a_i \stackrel{iid}{\sim}  2 \mathrm{Bernoulli}(1/2) -1$. 
For any realization of $a$, let $\ell_{\pm} = \ell_{\pm}(a) \in R^n$ be the vectors of $(\pm 1,  a')'$, respectively.  
Let $F(a)$, $F(\mu)$, and $F(w)$ be the CDF of $a$ and $\mu$ respectively, and let $h(z)$ be the (joint) density of the matrix $Z$. 
It follows that 
$g_{\pm}(x) = \int h(x - \ell_{\pm}(a) \mu' - w) d F(a) d F(\mu) d F(w)$, $x  \in R^{n,p}$,  and $\| g_{+}  - g_{-}\|_1$ equals to  
\[
\int  |\int [h(x - \ell_{+}(a) \mu' - w) - h(x - \ell_{-}(a) \mu' - w)]  d F(a) dF(\mu) d F(w)| dx. 
\]   
Using Fubini's theorem, this is no greater than 
$\int G(w) d F(w)$, where $G(w) =  |\int h(x - \ell_{+}(a) \mu' - w) - h(x - \ell_{-}(a) \mu' - w) d F(a) dF(\mu) | dx$. 
Note that for any fixed $w \in R^{n,p}$, $A(w)$ does not depend on $w$ and equals to $\|f_{+} - f_{-}\|_1$,  
and the claim follows. \qed

%%%%%%%%%
%%%%%%%%%
%%%%%%%%%  
\subsection{Proof of Theorem \ref{thm:siga}}    
For each $1 \leq j  \leq p$, 
consider the testing of two hypotheses, $H_{0}^{(j)}:   \mu(j)  = 0$ versus $H_1^{(j)}:      \mu(j) = \tau_p$.
Let $f_{0}^{(j)}$ and $f_1^{(j)}$ be the joint density of $X$ under $H_0^{(j)}$ and $H_1^{(j)}$, respectively. Since $P(\mu(j) = \tau_p) = \eps_p$,     it follows from the connection between $L^1$-distance and the sum of Type I and Type II testing errors  \cite{Aad} that for any clustering procedure $\hat{\mu}$,  
\begin{align*} 
P(\sgn(\hat{\mu}(j))   \neq \sgn(\mu(j)))  = &   (1 - \eps_p) P(\hat{\mu}(j) \neq 0 | \mu(j)  = 0) + \eps_p P(\hat{\mu}(j)  = 0 | \mu(j)  = \tau_p)   \\
 \geq  &  (1/2) [1 - \|(1 - \eps_p) f_{0}^{(j)} - \eps_p f_{1}^{(j)}\|_1] \\ 
 \geq  &   \eps_p [1 -  (1/2)  \| f_0^{(j)} - f_1^{(j)}\|_1], 
\end{align*} 
where in the last step we have used $\|(1 - \eps_p) f_0^{(j)} - \eps_p f_1^{(j)} \|_1  = \| (1 - 2 \eps_p) f_0^{(j)} + \eps_p (f_0^{(j)} - f_1^{(j)})\|_1  \leq (1 - 2 \eps_p) + \eps_p \| f_0^{(j)} - f_1^{(j)}\|_1$. 
Comparing this with the desired claim, it suffices to show that for all $1 \leq j  \leq p$,  
\begin{equation} \label{sigLB1} 
\| f_0^{(j)} -  f_1^{(j)}\|_1   =   o(1),   \;\;  \mbox{where $o(1) \goto 0$ and does not depend on $j$}. 
\end{equation} 

We now show (\ref{sigLB1}) for every fixed $1 \leq j  \leq p$. We first consider the case $\beta<1-\theta$. For short, we drop the superscript ``$(j)$"  in $f_0^{(j)}$ and $f_1^{(j)}$.  Recall that $X = \ell \mu' + Z =  [x_1, x_2, \ldots, x_p]$ and let  $\tilde{\mu} = \mu - \mu(j)  e_j$, where $e_j$ is the $j$-th standard basis vector of $R^p$; note that $\tilde{\mu}(j)   = 0$.   
Let $E$ denote the expectation under the law of $X = Z$.  By basic calculus and Fubini's theorem, 
\begin{align} 
\| f_0 -   f_1\|_1 & =  E[\bigl| \int [1   -   e^{\tau_p  \langle\ell, x_j\rangle - n \tau_p^2/2}]    e^{ \ell'  X \tilde{\mu} -  n \|\tilde{\mu}\|^2/2} d F(\tilde{\mu}) 
d F(\ell) \bigr|  \bigr]    \nonumber \\
&  \leq   \int  E \bigl[ |  1  -   e^{\tau_p  \langle\ell, x_j\rangle - n \tau_p^2/2} |  e^{ \ell'  X \tilde{\mu} -  n \|\tilde{\mu}\|^2/2}  \bigr]      d F(\tilde{\mu})    d F(\ell)   \nonumber \\
& = \int  E \bigl[ | 1   -   e^{\tau_p  \langle\ell, x_j\rangle - n \tau_p^2/2} | \bigr]  d F(\ell)  , 
\label{sigLB2}
\end{align} 
where in the last step, we have used the fact that  $x_j$ and $X \tilde{\mu}$ are independent and that $E[e^{ \ell'  X \tilde{\mu} -  n \|\tilde{\mu}\|^2/2}] = 1$.  Additionally, note that $E \bigl[ |1 -  e^{\tau_p  \langle \ell, x_j \rangle - n \tau_p^2/2} | \bigr]$ does not depend on $\ell$. Denote $z = n^{-1/2} \langle \ell, x_j\rangle$; note that $z \sim N(0,1)$.  
Inserting these into (\ref{sigLB2}) gives 
\begin{equation}\label{sigLB3}
\|  f_0 -   f_1\|_1  = E_0 \bigl[ | 1  -   e^{\sqrt{n} \tau_p  z - n \tau_p^2/2} | \bigr],  
\end{equation}
where $E_0$ denotes the expectation under the law of $z \sim N(0,1)$.  
By the conditions of $\alpha> \eta_{\theta}^{sig}(\beta)$ and $\beta < (1 - \theta)$,  we have $\alpha > \theta/2$, and 
$n \tau_p^2 = p^{\theta -2\alpha} = o(1)$. In this simple setting,  it is seen that $E_0 \bigl[ | 1  -   e^{\sqrt{n} \tau_p  z - n \tau_p^2/2} | \bigr]  = o(1)$.  
Combining (\ref{sigLB1})-(\ref{sigLB3}) gives the claim. 

We now consider the case $\beta > (1 - \theta)$.   In this case,  
$\eta_{\theta}^{sig}(\beta) = \eta_{\theta}^{hyp}(\beta)$, so intuitively, the claim follows by the argument that ``as long as 
it is impossible to have (global) hypothesis testing, it is impossible to identify the signals". Still, for mathematical rigor, 
it is desirable to  provide a proof using the $L^1$-distance. 
Similarly to that in the proof on the lower bound for global testing, 
write $\mu=\tilde{\mu} + \mu(j)e_j$ and let $d_p = (6 p \eps_p \log(p))^{1/2}$, 
$A_s$ be the event $\{\|\tilde{\mu}\|_0 = s\}$ and $F_s$ be the conditional distribution of $\tilde{\mu}$ given the event of $A_s$,  $1 \leq s \leq p$. 
Define 
$
a_s = \int  e^{\tau_p \langle \ell, x_j \rangle - n \tau_p^2/2} e^{\ell' X \tilde{\mu} - n \|\tilde{\mu}\|^2/2}  d F_s(\tilde{\mu}) d F(\ell)$ and  
$\tilde{a}_s = \int e^{\ell' X \tilde{\mu} - n \|\tilde{\mu}\|^2/2}  d F_s(\tilde{\mu}) d F(\ell)$. 
It suffices to show that for all $s$ such that $|s - p \eps_p| \leq d_p$ that 
\begin{equation}\label{sigLB4}
E[(a_s - \tilde{a}_s)^2] = o(1). 
\end{equation}

Let $\nu$ be an independent duplicate of $\mu$. 
By similar arguments and noting that $\mu' \nu = \tilde{\mu}' \tilde{\nu} + \tau_p^2$ and $\tilde{\mu}'\nu=\tilde{\mu}'\tilde{\nu}$, we have 
$
E[a_s^2] = \int [\cosh(\tilde{\mu}' \tilde{\nu} + \tau_p^2)]^n  d F_s(\tilde{\mu}) d F_s(\tilde{\nu})$,   
$E[\tilde{a}_2^2] = \int [\cosh(\tilde{\mu}' \tilde{\nu})]^n  d F_s(\tilde{\mu}) d F_s(\tilde{\nu})$,
and the cross term
$
E[\tilde{a}_s a_s] = \int [\cosh(\tilde{\mu}' \tilde{\nu})]^n  d F_s(\tilde{\mu}) d F_s(\tilde{\nu})$.
Combining these terms and noting that $\cosh(x+y) = \cosh(x) [1 + \tanh(x) \tanh(y)]$,  there is
\[
E[(a_s - \tilde{a}_s)^2] = \int [\cosh(\tilde{\mu}' \tilde{\nu})]^n \{[1 +  \tanh(\tau_p^2) \tanh(\tilde{\mu}' \tilde{\nu})]^n -1\} d F_s(\tilde{\mu}) d F_s(\tilde{\nu}). 
\]
Now, over the event $\{(\tilde{\mu}, \tilde{\nu}):  \|\tilde{\mu}\|_0 = \|\tilde{\nu}\|_0 = s\}$,  where $s \sim p \eps_p$, 
we have 
$|\tilde{\mu}' \tilde{\nu}| \leq  s \tau_p^2  \lesssim p \eps_p \tau_p^2 \leq p^{(1 - \beta - \theta)/2}$; note that by the assumption of $r > \eta_{\theta}^{hyp}(\beta)$ and $\beta > (1 - \theta)$, the exponent $(1 - \beta - \theta)/2 < 0$. 
As a result,  
it is seen that $
\tanh(\tilde{\mu}' \tilde{\nu}) \tanh(\tau_p^2) \lesssim \tilde{\mu}' \tilde{\nu}  \tau_p^2  \lesssim p \eps_p \tau_p^4$, 
where $ p \eps_p \tau_p^4 = o(n^{-1})$ by the assumption of $\alpha > \eta_{\theta}^{sig}(\beta)$. Inserting this into (\ref{sigLB4}) gives 
\begin{equation}\label{sigLB5}
E[(a_s - \tilde{a}_s)^2] = o(1)  \cdot  \int [\cosh(\tilde{\mu}' \tilde{\nu})]^n  d F_s(\tilde{\mu}) d F_s(\tilde{\nu}). 
\end{equation}
According to (\ref{gtLB3})-(\ref{gtLB3a}) in Section~\ref{subsec:pfhypa}, the second term on the right hand side of \eqref{sigLB5} is $1 + o(1)$. This gives the claim. \qed

%%%%%%%%%%%%
%%%%%%%%%%%%
%%%%%%%%%%%%
\subsection{Proof of Theorem \ref{thm:hypa}}\label{subsec:pfhypa}  
Recall that $X = \ell \mu' + Z$. Let $f_0(X)$ and $f_1(X)$   be the joint density of $Z$ and $X$, respectively.   
It is sufficient to show that as $p \goto \infty$,  under the conditions of Theorem \ref{thm:hypa},  
\begin{equation} \label{gtLB1} 
\|f_1 - f_0\|_1 \goto 0. 
\end{equation}   
Recall that $\|\mu\|_0$ and $\|\mu\|$ denote the $L^0$-norm and the $L^2$-norm 
of $\mu$ respectively.  For $1 \leq s \leq p$, let $A_s$ be the event $A_s = \{ \|\mu\|_0 = s\}$, %and define $\hat{\eps}_p = s/p$.  
$F(\ell)$ and $F(\mu)$  be the distributions of $\ell$ and $\mu$, respectively, and 
let $F_s(\mu)$ be the conditional distribution of $\mu$ given the event of $A_s$.  
Introduce a constant $d_p = (6p \eps_p \log(p))^{1/2}$, a set $D_p = \{s: |s - p\eps_p| < d_p\}$, and functions $a_s(X) = \int e^{\ell' X \mu - n \|\mu\|^2/2 }   d F_s(\mu) d F(\ell)$, $1 \leq s \leq p$. 
Let $E$ be the expectation under the law of $X = Z$.  It is seen that $f_1(X) / f_0(X) = \int e^{\ell' X \mu - n \|\mu\|^2 /2} d F(\mu) d F(\ell) = \sum_{s = 1}^p P(A_s) a_s(X)$, and so $\|f_1 - f_0\|_1$ equals to
%%%%%%%%%%%
%%%%%%%%%%%
%%%%%%%%%%%
\begin{equation} \label{gtLB2} 
E \bigl|\sum\nolimits_{s = 1}^p P(A_s)  \bigl( a_s(X) - 1 \bigr) \bigr| \leq  \sum\nolimits_{D_p} P(A_s) E\bigl[|a_s(X) -1| \bigr]   + rem,  
\end{equation} 
where $rem  =  \sum_{D^c_p} P(A_s) E[|a_s(X) - 1|]$.  
Since $E[|a_s -1|] \leq E[a_s]  + 1 = 2$,  $rem \leq \sum_{D_p^c} 2 P(A_s)   \leq 2 P\bigl( \|\mu\|_0 \in D_p^c)$. Note that  $\|
\mu\|_0 \sim \mathrm{Binomial}(p, \eps_p)$, where $p \eps_p = p^{1 - \beta}$ with $0 < \beta < 1$,  it follows from basic statistics that 
$rem  = o(1)$.  At the same time, by Cauchy-Schwarz inequality,   
$(E[|a_s(X) - 1|])^2 \leq  E[(a_s(X) - 1)^2] = E[a^2_s(X)] - 1$.   
Combining these with (\ref{gtLB2}), to show (\ref{gtLB1}), it suffices to show  that 
\begin{equation} \label{gtLB3} 
E[(a_s^2(X)]  \leq 1 +    o(1), \qquad \forall s \in D_p,  
\end{equation} 
where $o(1) \goto 0$ uniformly for all such $s$ as $p \goto \infty$.  

We now show (\ref{gtLB3}).  Fix an $s \in D_p$. Let $\nu \in R^p$ be an independent copy of $\mu$, and let $F_s(\nu)$ be the distribution of $(\nu | \{\|\nu\|_0 = s\})$. 
Using basic statistics and the independence of $X_i$, 
\begin{eqnarray*}
a_s^2(X)   =  \int e^{-n\|\mu\|^2/2 - n\|\nu\|^2/2} \Pi_{i = 1}^n [\cosh(\mu' X_i) \cosh(\nu' X_i)]  d F_s(\mu) d F_s(\nu).   
\end{eqnarray*}
First, by the independence of $X_i$ and basic statistics,  $E[a_s^2(X)]$ equals to
\begin{equation} \label{gtLB3a}  
\int [\cosh(\mu' \nu)]^n dF_s(\mu) d F_s(\nu) = 
 \sum_{k = 0}^n   \int {n \choose k} \frac{e^{(2k - n) \mu' \nu}}{2^n} dF_s(\mu) d F_s(\nu). 
\end{equation} 
Recalling that any nonzero entry of $\mu$ or $\nu$ is $\tau_p$, it is seen that over the event $\{\|\mu\|_0 = \|\nu\|_0 = s\}$,  $\tau_p^{-2} \langle \mu, \nu\rangle$ is distributed as a hyper-geometric distribution $H(p, s, s)$.  Write $\hat{\eps}_p=s/p$. As $s \in D_p$, $\hat{\eps}_p \sim \eps_p$. Following \cite{Aldous}, there is a $\sigma$-algebra ${\cal B}$   and a random variable $b \sim \mathrm{Binomial}(s, \hat{\eps}_p)$ such that $\tau_p^{-2} \langle \mu, \nu\rangle$ has the same distribution as that of $E[b | {\cal B}]$. Using Jensen's inequality, $e^{(2k-n)\mu'\nu}\leq E[e^{(2k-n)\tau_p^2b}|{\cal B}]$, for $0 \leq k \leq n$. It follows that 
\begin{equation} \label{gtLB3b} 
E\int e^{(2k - n) \mu' \nu} dF_s(\mu) d F_s(\nu) \leq E[e^{ (2k - n) \tau_p^2 b}] = (1 - \hat{\eps}_p + \hat{\eps}_p e^{(2k - n) \tau_p^2})^s. 
\end{equation} 
Inserting (\ref{gtLB3b}) into (\ref{gtLB3a}) and rearranging,  
\begin{equation} \label{gtLB4} 
E[a_s^2(X)] \leq  2^{-n} \sum_{k = 0}^n {n \choose k}[1 - \hat{\eps}_p + \hat{\eps}_p e^{(2 k - n) \tau_p^2}]^s. 
\end{equation} 

We now analyze the right hand side of (\ref{gtLB4}).  
Denote $S$ by $\{1, 2, \ldots, n\}$. We split $S$ as the union of three disjoint subsets $S = S_1 \cup S_2 \cup S_3$,  where $S_1 = \{k \in S:  |2k - n|   <   \sqrt{n} \log(n) \}$, $ S_3 = \{k \in S:   |2 k - n| >  n \wedge  \sqrt{2 \log(n) n p \eps_p} \}$. 

Also, let $\tilde{\tau}_p  = p^{- \eta_{\theta}^{hyp}(\beta)}$. By our assumption of $\alpha > \eta_{\theta}(\beta)$, there is a constant $\delta = \delta(\theta, \alpha) > 0$ such that 
$\tau_p^2 = p^{-\delta} \tilde{\tau}_p^2$. 
We also claim that when $\alpha >  \eta_{\theta}^{hyp}(\beta)$,  
$\tau_p^2 |2k - n| = o(1)$ for any $k \in S_1 \cup S_2$. 
In fact,  by definitions and direct calculations, we have $\eta_{\theta}^{hyp}(\beta) > \theta/2$ when $\beta < \max\{1 - \theta, (2 - \theta)/4\}$ and $\eta_{\theta}^{hyp}(\beta) = (1 + \theta - \beta)/4$ otherwise. In the first case, recalling $n = p^{\theta}$,  the claim follows since $\tau_p^2 |2k - n| \leq \tau_p^2 n = p^{\theta - 2\alpha}$ and $\alpha > \theta/2$.  
In the second case, noting that $\tau_p^2 = p^{-\delta} \tilde{\tau}_p^2 = p^{-\delta}  (n p\eps_p)^{-1/2}$, it follows $|2k - n| \tau_p^2 \leq \sqrt{2 \log(n) n p \eps_p} \cdot (p^{-\delta} (n p \eps_p)^{-1/2}) = o(1)$ for all $k \in S_1 \cup S_2$, and the claim follows.  
Now,  since for any $x \in (-1, 1)$ and $y\in R$,   $1 - \hat{\eps}_p + \hat{\eps}_p e^{x}  \leq 1 + 2 \hat{\eps}_p |x| \leq e^{2 \hat{\eps}_p |x|}$, and $1 - \hat{\eps}_p + \hat{\eps}_p e^{y} \leq 1 - \hat{\eps}_p + \hat{\eps}_p e^{|y|} \leq e^{|y|}$,
\begin{equation} \label{gtLB5a}
[(1 - \hat{\eps}_p + \hat{\eps}_p e^{\tau_p^2 (2k - n)}]^s \leq 
\left\{
\begin{array}{ll}
 [1 + 2 \hat{\eps}_p \tau_p^2 |2k - n|]^s \leq e^{ 2 p\hat{\eps}_p^2  \tau_p^2 |2k - n|}, & \, k \in S_1 \cup S_2,  \\ 
 (e^{\tau_p^2 |2 k -n|})^s  = e^{s  \tau_p^2 |2k - n|},  & \, k \in S_3.  
\end{array}
\right.  
\end{equation} 
If we take $Y \sim \mathrm{Binomial}(n, 1/2)$, then $P(Y = k) = 2^{-n} {n \choose k}$. At the same time, by de Moivre-Laplace Theorem and Hoeffding inequality \cite{Wellner86}, 
\begin{equation} \label{gtLB5b} 
P(Y = k)  \left\{ 
\begin{array}{ll}  
\sim (\pi n /2)^{-1/2} e^{-(2 k - n)^2 / (2n)},  &\quad k \in S_1, \\
\leq e^{- (2k - n)^2/(2n)},  &\quad  k \in S_2 \cup S_3.
\end{array} 
\right. 
\end{equation} 
Combining (\ref{gtLB5a})-(\ref{gtLB5b}),  we have the following. First,  the summation over $k \in S_1$ is smaller than that that ($\phi$ is the probability density of $N(0,1)$)
\begin{equation} \label{gtLB6a} 
(2/\sqrt{n})  \sum_{k \in S_1} e^{2 p \hat{\eps}_p^2 \tau_p^2 |2k - n|} \phi(\frac{2k - n}{\sqrt{n}}) \sim \int_{-\log(n)} ^{\log(n)}  e^{2 p \hat{\eps}_p^2 \tau_p^2 \sqrt{n} x } \phi(x) dx.  
\end{equation} 
By the assumption of $\alpha > \eta_{\theta}^{hyp}(\beta)$ and basic algebra, we have $\alpha > (2 + \theta - 4 \beta)/4$. It follows that $p \hat{\eps}_p^2 \tau_p^2 \sqrt{n} \sim 2 p^{-\delta} p \eps_p^2 \tilde{\tau}_p^2 \sqrt{n} = 2 \cdot p^{-\delta} \cdot  p^{(2 + \theta - 4 \beta)/2 - 2 \eta_{\theta}^{hyp}(\beta)}$, where the exponent is negative.  It follows that the right hand side of (\ref{gtLB6a}) is $1 + o(1)$.  
Second, let $II$ be the summation over $k \in S_2$, then
\begin{equation} \label{gtLB6b} 
II  \leq \sum_{k \in S_2} e^{2 p \hat{\eps}^2_p \tau_p^2 |2 k - n|} e^{- (2k -n)^2/(2n)} \leq  \sum_{k \in S_2} e^{- (2k - n)^2/(2n)},  
\end{equation} 
where the second inequality is because $2 p \hat{\eps}^2_p \tau_p^2 \leq 2p^{-\delta}/\sqrt{n}\leq |2k-n|/(4n)$. The right hand side does not exceed $n e^{-  \log^2(n)} = o(1)$ since $|2k - n| \geq \sqrt{n}\log(n)$.  
Last, we consider the summation over $k \in S_3$. We only consider the case of $\beta > (1 - \theta)$ since 
only in this case $S_3$  is  non-empty. Note that in this case, 
$n \wedge \sqrt{2 \log(n) n p\eps_p} = \sqrt{2 \log(n) n p \eps_p}$ and that 
for any $k \in S_3$, $s \tau_p^2 \lesssim p^{-\delta} p \eps_p \tilde{\tau}_p^2  \leq  |2k - n|/(4n)$,  
\begin{equation} \label{gtLB6c}
III \leq \sum_{k  \in S_3} e^{s \tau_p^2 |2 k -n|  - (2k - n)^2/(2n)}  \leq  \sum_{k \in S_3} e^{- (2k - n)^2/(4n)}, 
\end{equation} 
which $\leq n e^{- \log(n) p \eps_p/2}   = o(1)$.  Combining (\ref{gtLB6a})-(\ref{gtLB6c}) with (\ref{gtLB4}) gives the claim. \qed   

%%%%%%%%%%%
%%%%%%%%%%%
%%%%%%%%%%%
\section{Discussions}  \label{sec:discuss}
We have studied the statistical limits for three interconnected problems: clustering, signal recovery, and hypothesis testing. For each  problem,  
in the two-dimensional phase space calibrating the signal sparsity and strength, we identify the exact separating boundary for the Region of Possibility and  Region of Impossibility.    We have also derived a computationally tractable upper bound (CTUB),  part of which is tight, and the other part is conjectured to be tight.  
Our study on the limits are extended  to the case where the parameters fall exactly on the separating boundaries and  the case of colored noise. 

We propose several different methods, including IF-PCA.  IF-PCA  is  a two-fold dimension reduction algorithm:  
we first reduce dimensions from (say) $10^4$ to a few hundreds by screening, and then further reduce it to just a few by PCA. Each of the two steps can be useful in other high-dimensional settings.   
 Compared to popular penalization approaches, our approach has advantages for it is highly extendable and computationally inexpensive.

The work is closely related to Jin and Wang \cite{IFPCA}   but is also very different. 
The focus of \cite{IFPCA} is to investigate the performance of IF-PCA with real data examples and to study the 
consistency theory. The primary focus here,  however, is on the statistical limits 
for three problems including clustering. The paper is also closely related to the very interesting paper by Arias-Castro and Verzelen \cite{Ery}. However, two papers are different in important ways. 
\begin{itemize} 
\item  The focus of our paper is on clustering, while the focus of their paper is on hypothesis testing (without careful discussion on clustering). 
\item  Both papers addressed  signal recovery, but there are 
important  differences: we provided the statistical lower bound but they did not;  
the CTUB they derived is not as sharp as ours. See Figure \ref{fig:complimit}. 
\item Both papers studied hypothesis testing, but since the models  are different, the 
separating boundaries (and so the proofs) are also different.  See Sections~\ref{subsec:other-limits} and ~\ref{sec:hyp} (also Figure~\ref{fig:complimit}) for details. 
\item Both papers studied the case with colored noise, besides the different focuses (clustering v.s. hypothesis testing), their setting in the colored case is also different from ours. In their setting, coloration makes a substantial difference to statistical limits.
\end{itemize} 
For these reasons,  the methods and theory (especially that on IF-PCA) in our paper are very different from those in \cite{Ery}.
With that being said, we must note that since two papers have overlapping interest, it is not surprising that  certain part of this  paper overlaps\footnote{Compare the critical signal strength required for successful  hypothesis testing/signal recovery in our paper with those in \cite{Ery}, we note some discrepancies in terms of some multi-logarithmic factors.  This is due to that we choose a simpler calibration than that in \cite{Ery}:   all the parameters $(n,\eps,\tau)$ are expressed as a (constant) power of $p$ and multi-logarithmic factors are neglected. Such a calibration makes the presentation more succinct.} with that in \cite{Ery} (e.g., some parts of the separating boundaries and some of the ideas and methods). 

The paper is related to recent ideas in spectral clustering (e.g.,   
Azizyan {\it et al} \cite{WassermanSparse}, Chan and Hall \cite{CH}; see also \cite{PanShen2007, raftery2006,Sun2012,skmeans}).  In particular, 
the high level idea of IF-PCA (i.e., combining feature selection with classical methods)  is not new and can be found in    \cite{WassermanSparse, CH},  but the methods and theory are  different.   
Azizyan {\it et al}  \cite{WassermanSparse} study the clustering 
problem in a closely related setting, but they use a different loss function and so the separating boundaries are also different.  Chan and Hall \cite{CH} use   a very different screening idea (motivated by real data analysis) and do not study phase transitions.

%their bound that $\eta^{clu, ASW}_\theta(\beta) \max\{(1 - \beta)/2, (1 + \theta - \beta)/4\}$ is not sharp. 
%The paper also has close connection with the work of Azizyan, Singh and Wasserman \cite{WassermanSparse}. Their paper studied the statistical lower bound and CTUB for clustering problem when the signal is moderately sparse and strong, and they found there is a gap between the two bounds. In our paper, we derived exact statistical limits for the whole range, and we derived more tight CTUB with four clustering methods different from theirs. 

Our work is closely related to recent interest in the spike model (e.g., \citep{amini2009high, Lei, WLL}).  
In particular, {\it mathematically},  Model (\ref{model2}) is similar to the spike model \cite{JohnstoneLu}, 
and theoretical results on one can shed light  on those for the other.  However, two models are also different from a {\it scientific perspective}: (a) two models are motivated by different application  problems, (b) the primary interest of Model (\ref{model2}) is on the class labels $\ell_i$, which are sometimes easy to validate in real applications, and (c) the primary interest of spike model is on the feature vector $\mu$, which is relatively hard to validate in real applications. The focus and scope of our study are very different from many recent works on the spike model, and most part of the bounds (especially those for clustering and IF-PCA) we derive are new.  

This paper is also related to  the recent interest on computationally tractable lower bounds and sparse PCA \citep{BR13learning, CMW}, but it is also very different in terms of our focus on clustering and statistical limits. 
%%%%%
It is also related to the lower bound for hypothesis testing problem \citep{LugosiTest} and the  sub-matrix detection problem \citep{MW}, but the model is different. Recovering of $\ell$ and $\mu$ can also be interpreted as recovering a low-rank matrix from the data matrix, which is closely related to the low rank matrix recovery studies \citep{candes2009exact}.
%%%%%
In terms of the phase transitions, the paper is closely related to \cite{DJ04} on signal detection, \cite{DJ08} on classification, and \cite{Ke} on variable selection, but is also very different for the primary focus here is on clustering.

For simplicity, we focus on the ARW model, where we have several assumptions such as $\ell_i = \pm 1$ equally likely, 
the signals have  the same sign and  equal strength, etc.  
Many of these assumptions  can be largely relaxed. For example, 
 Theorems \ref{thm:clu-a}-\ref{thm:clu-c} continue to hold if we replace the model  
$\mu(j) \stackrel{iid}{\sim}  (1 - \eps_p) \nu_0 + \eps_p \nu_{\tau_p}$ by that of 
$\mu(j) \stackrel{iid}{\sim}  (1 - \eps_p) \nu_0 + \eps_p G_p$, where $G_p$ is a distribution supported in the interval   
$[a_p \tau_p, b_p \tau_p]$   
with $0 <  \max\{ a_p^{-1}, b_p\} \leq  L_p$ (a multi-$\log(p)$ term).   Also, in 
 Section~\ref{subsec:other-limits}, we have  discussed the case where we replace $\mu(j) \stackrel{iid}{\sim}  (1 - \eps_p) \nu_0 + \eps_p \nu_{\tau_p}$ in Model (\ref{ARW1})  by that of $\mu(j) \stackrel{iid}{\sim}  (1 - \eps_p) \nu_0 + a\eps_p \nu_{-\tau_p} + (1-a)\eps_p\nu_{\tau_p}$ for a constant $0\leq a\leq 1/2$. Theorems \ref{thm:clu-a}-\ref{thm:clu-c} continue to hold if $a\neq 1/2$. If $a=1/2$, the left part of the boundaries will change and the aggregation methods need to be modified. We discuss this case in detail in Section~\ref{app:C}.  It requires a lot of time and effort to fully investigate how broad the main theorems hold, 
so we leave it to the future.

The paper motivates an array of interesting problems in post-selection Random Matrix Theory that could be future research topics.  For the perspective of spectral clustering, it is of great interest to precisely characterize the limiting behavior of the singular values  (bulk and the edge singular values) and leading singular vectors of  the post-selection  data matrix. These problems are  technically very challenging, and we leave them to the future. 

Our paper supports the philosophy in Donoho \cite[Section 10]{50years} that simple and homely methods 
are just as good as more charismatic methods in Machine Learning  for analyzing (real) high dimensional data.   
%%%%%%%%%%%%%

\appendix  
\setcounter{equation}{0}

%%%%%%%%%%%
\section{An extension of the ARW Model} \label{app:D}

We consider an extension of the Asymptotic Rare and Weak (ARW) model in Section~\ref{subsec:ARW}, where Models \eqref{model1}-\eqref{model2} and the calibration \eqref{ARW2} continue to hold but \eqref{ARW1} is replaced by a more sophisticated signal configuration: 
\beq  \label{ARW4}
\mu(j) \stackrel{iid}{\sim} (1 - \eps) \nu_0  + a \cdot  \eps \cdot  \nu_{-\tau} + (1 - a) \cdot  \eps \cdot \nu_{\tau}, \qquad 1\leq j\leq p,
\eeq
where $0\leq a\leq 1/2$ is a constant. This extended model includes the original ARW as a special case with $a=0$. In this extension, we allow the nonzero coordinates of the feature vector $\mu$ to have positive and negative signs.   Due to such a change,  we need to  slightly modify the definition of the (normalized) Hamming distance for signal recovery: $\hamm_p(\hat{\mu}, \alpha,\beta,\theta) = (p \eps_p)^{-1} \sum_{j = 1}^p P(\sgn(\mu(j))\neq \sgn(\hat{\mu}(j)))$. The loss functions for clustering and hypothesis testing remain the same.

When $0 <  a < 1/2$, with high probability, the majority of the nonzero coordinates of $\mu$ are positive, and the performance of the four methods in Section~\ref{subsec:4method} is not affected. Furthermore, the statistical limits and CTUB for all three problems continue to hold.  For brevity, we omit the details. 

The case of $a=1/2$ is more delicate. In this case, the two aggregation methods turn out to be ineffective. In light of this,   we introduce a variant of the Sparse Aggregation, where we cluster the $n$ subjects  by 
\beq \label{SAclu}
\hat{\ell}_N^{(sa)} = \sgn(X\hat{\mu}_N^{(sa)}).  
\eeq
Here,  
\beq \label{SAsig}
\hat{\mu}_N^{(sa)} = \mathrm{argmax}_{\mu\in\{-1,0,1\}^p: \|\mu\|_0=N} \| X\mu\|_1. 
\eeq
Also, we use $\hat{\mu}_N^{(sa)}$ to estimate the sign of $\mu$ (i.e., for signal recovery), and use the test statistic
\beq \label{SAhyp}
\hat{T}_N^{(sa)} = N^{-1/2}\|X\hat{\mu}_N^{(sa)}\|_1
\eeq 
for hypothesis testing. Note that if we force $\mu(j)\in\{0,1\}$ in \eqref{SAsig}, then it reduces to the original Sparse Aggregation.

%For Simple Aggregation, we have not found a variant that both achieves the statistical limit and is computationally tractable. However, in the less sparse case, the classical PCA turns out to be optimal. We extend the classical PCA to methods for signal recovery and hypothesis. Let $\lambda_1(X)$, $\xi_1(X)$  and $\eta_1(X)$ be the leading singular value of $X$ and the associated left and right singular vectors, respectively. 

{\bf Remark}. We have not found a variant of Simple Aggregation that both achieves the statistical limit and is computationally tractable. However, in the less sparse case, the classical PCA turns out to be already optimal.
\footnote{Classical PCA for hypothesis testing is to reject the null hypothesis when the leading singular value of $X$ is larger than $\sqrt{p} + \sqrt{n} + \log(p)$; for signal recovery is as the description in Section~\ref{sec:sig}. 
	%Classical PCA as a method for signal recovery and hypothesis testing are described in Section~\ref{sec:sig} and Section~\ref{sec:hyp}, respectively. 
	However, for signal recovery, since we need to estimate not only the support but also the sign of $\mu$, we slightly modify it to $\hat{\mu}_*^{(if)}(j)=\sgn(\hat{y}(j))\cdot 1\{|\hat{y}(j)|>2\sqrt{\log(p)}\}$, where $\hat{y}=X\hat{\ell}_*^{(if)}$ and $\hat{\ell}_*^{(if)}$ denotes the class label vector estimated by classical PCA.}

%%%%%%%%%%
%%%%%%%%%%
%%%%%%%%%%
%%%%%%%%%%
\begin{figure}[t!]
	\centering
	\includegraphics[height = 1.85 in]{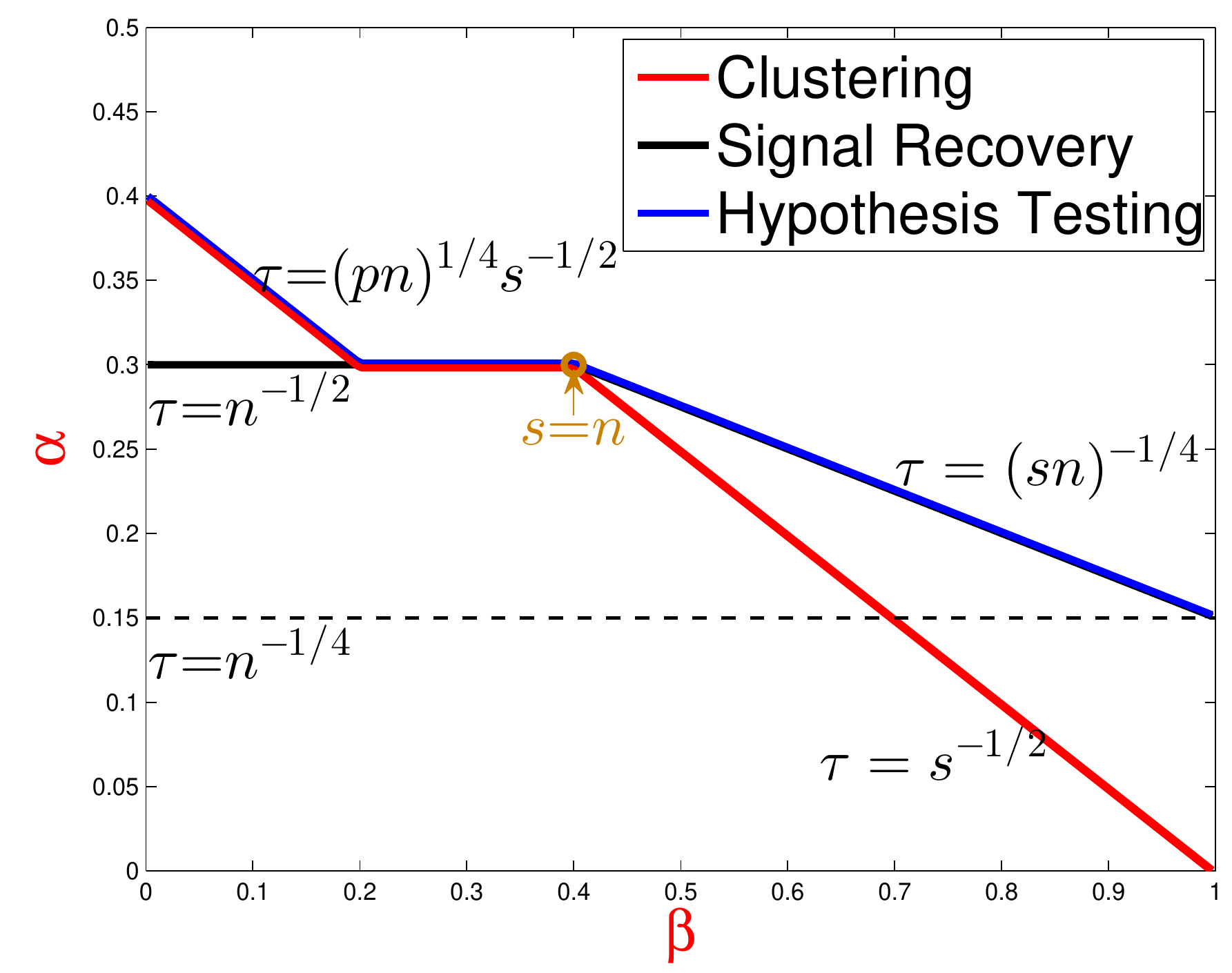}
	\includegraphics[height = 1.83 in]{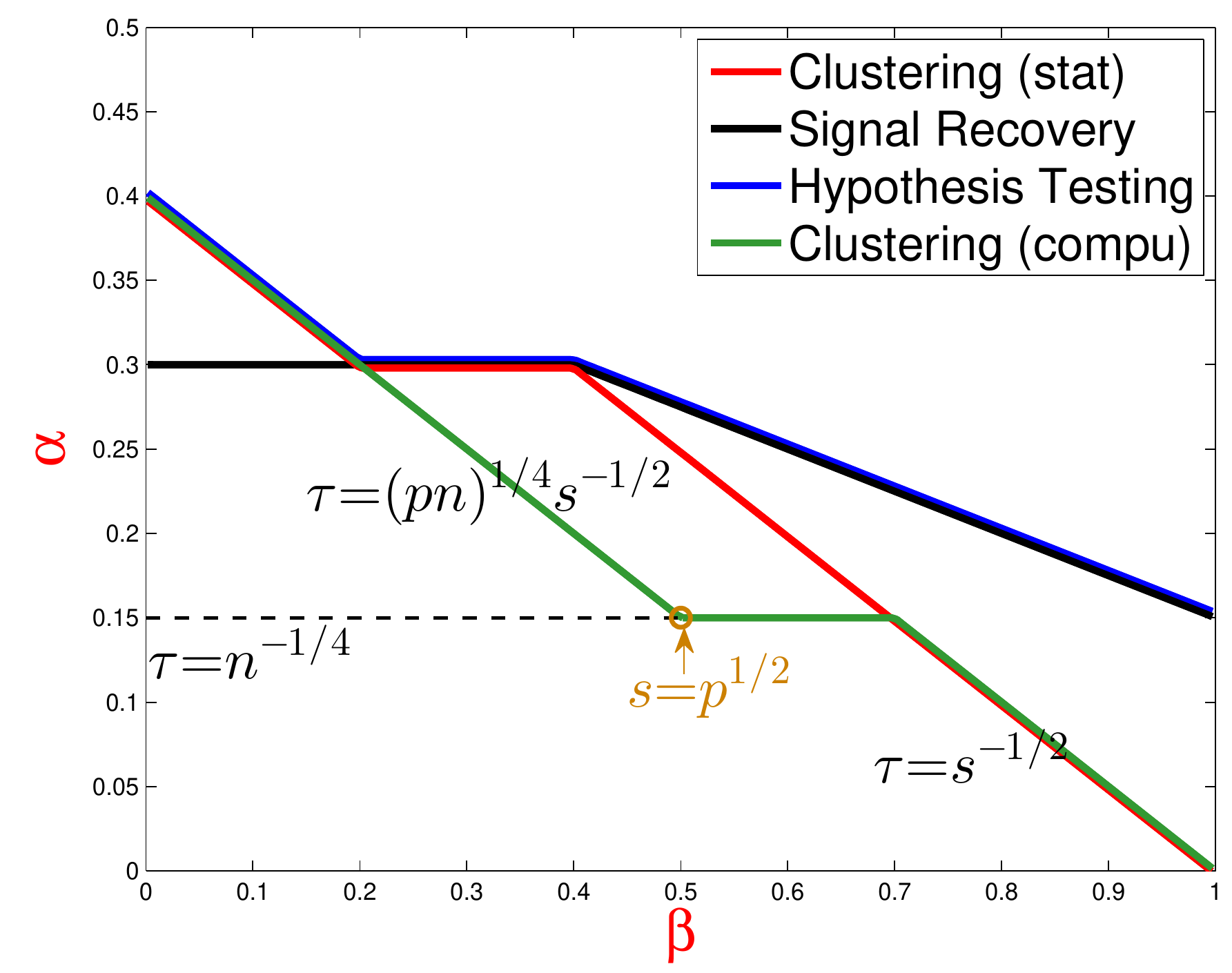} \\ 
	\includegraphics[height = 1.83 in]{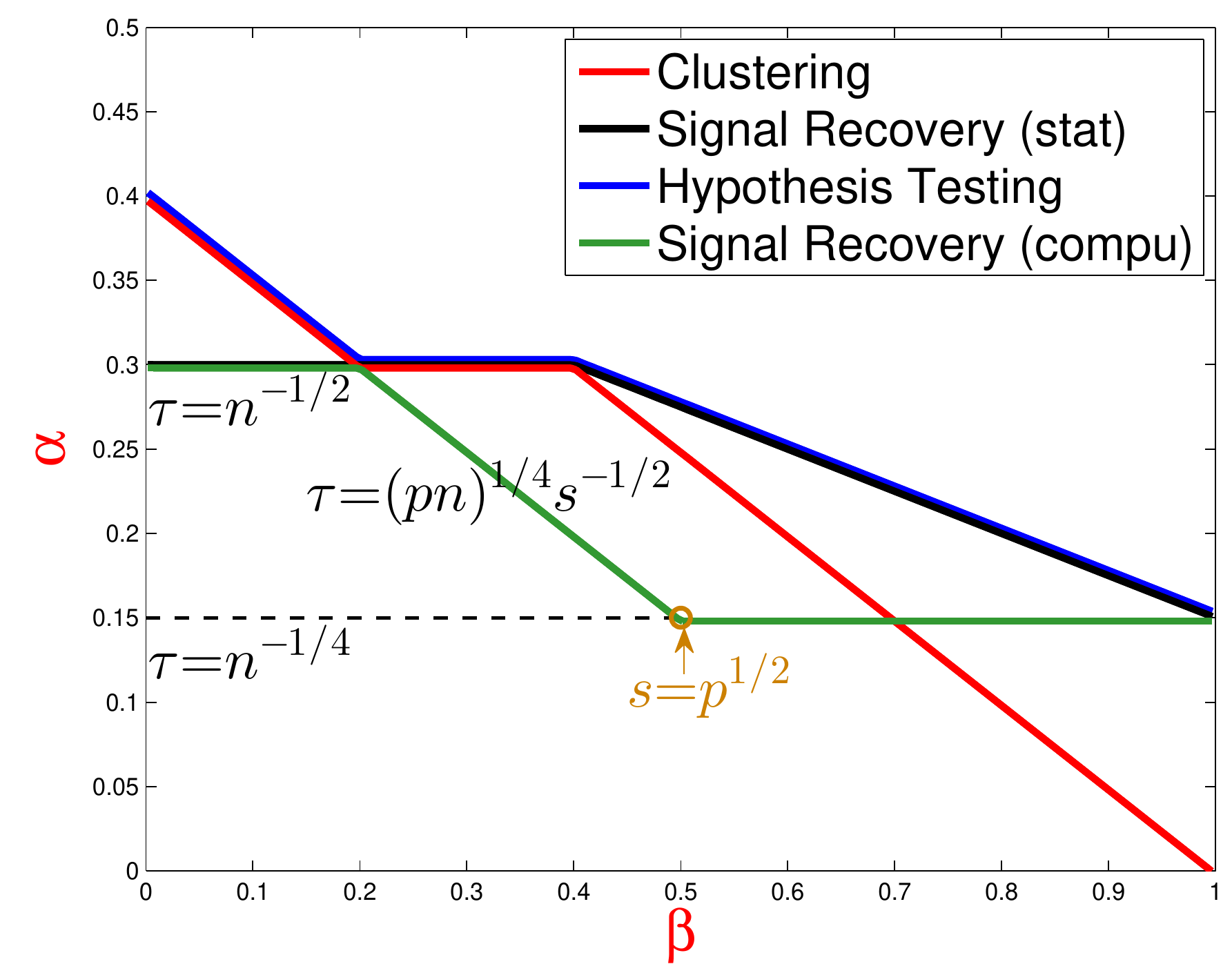}  
	\includegraphics[height = 1.83 in]{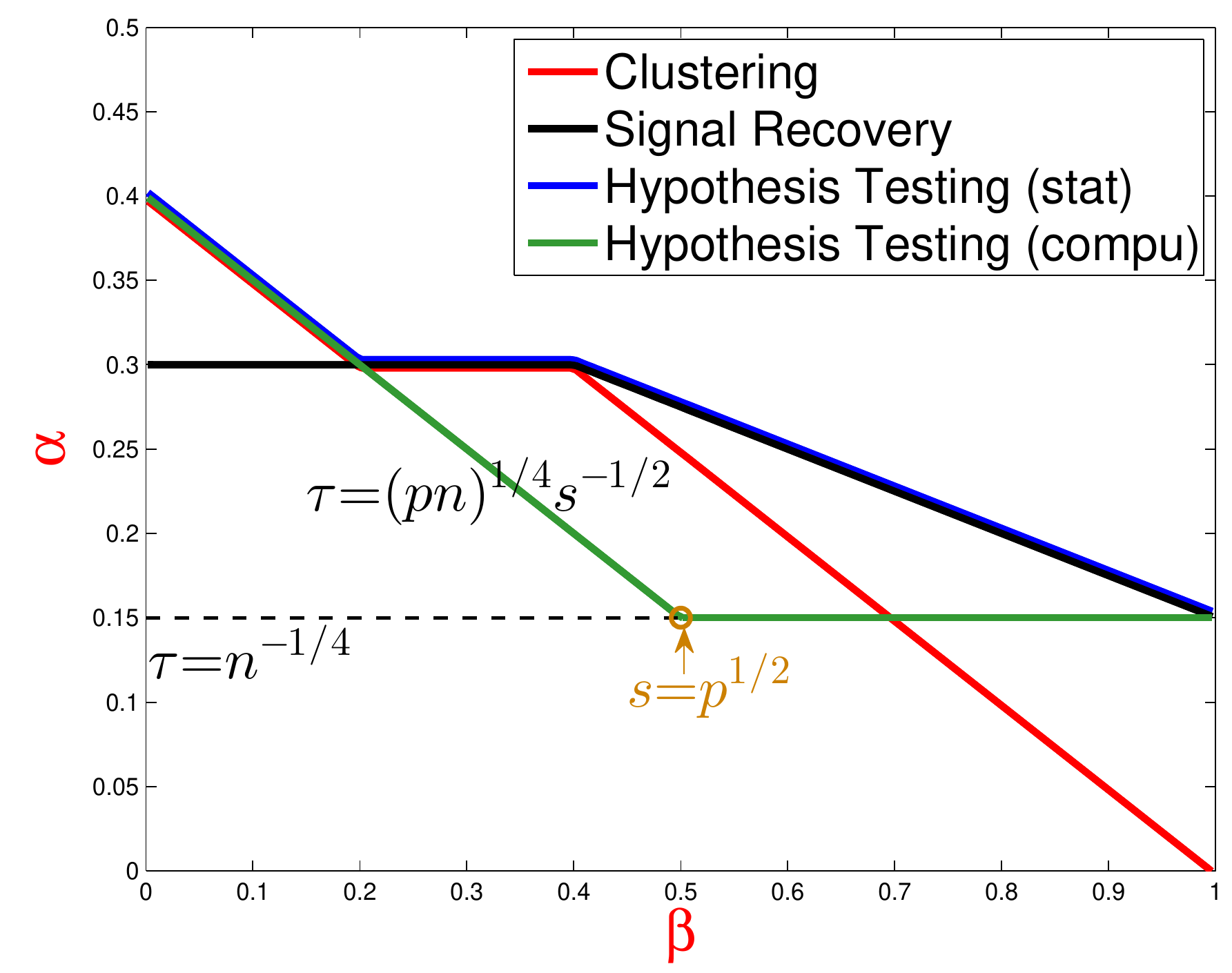} 
	\caption{Top left: statistical limits for  
		clustering (red), signal recovery (black), and hypothesis testing (blue); $s = p \eps_p$.  
		Other three panels: CTUB (green) for clustering (top right), signal recovery (bottom left)  and hypothesis testing (bottom right), respectively.}
	\label{fig:complimit-sign}
\end{figure}

We now present the statistical limits and CTUB for all three problems. They are different from the ones we present in the main paper \cite{JKW}. First, we look at the statistical limits. 
\begin{align*}
& \eta^{clu}_{\theta}(\beta)  = 
\left\{
\begin{array}{ll}
(1+\theta - 2 \beta)/4,	 	& \qquad   \beta < (1-\theta)/2,     \\
\theta/2, 				& \qquad  (1 - \theta)/2 < \beta < (1 - \theta),     \\
(1 - \beta)/2,	&  \qquad  \beta > (1 - \theta).  
\end{array}
\right.\cr
& 
\eta^{sig}_{\theta}(\beta) = 
\left\{
\begin{array}{ll}
\theta/2,		 				& \qquad  \beta < (1 - \theta),  \\
(1 + \theta - \beta)/4,	                         & \qquad  \beta > (1 - \theta).   
\end{array}
\right.\cr
& \eta^{hyp}_{\theta}(\beta)  = 
\left\{
\begin{array}{ll}
(1+\theta - 2 \beta)/4,	 	& \qquad   \beta < (1-\theta)/2,     \\
\theta/2, 				& \qquad  (1 - \theta)/2 < \beta < (1 - \theta),     \\
(1 + \theta - \beta)/4,	&  \qquad  \beta > (1 - \theta).  
\end{array}\right. 
\end{align*}
Figure~\ref{fig:complimit-sign} (top left panel) displays the statistical limits for three problems. Comparing  it with Figure~\ref{fig:complimit} (top left panel), we find that : (a) the black curve (signal recovery) remains the same, (b) the red curve (clustering) remains the same, except for the segment on the left is replaced by $\tau^4 = p/(ns^2)$, (c) for the blue curve (hypothesis testing), the right most segment remains the same, while the other two segments coincide  with those of the red curve. 

{\it Achievability}. The statistical limit of clustering is achieved by the classical PCA (the left segment) and the variant \eqref{SAclu} of Sparse Aggregation (the right two segments). For signal recovery, the right two segments are achieved by the modified Sparse Aggregation \eqref{SAsig}, and the left segment is achieved by  classical PCA. For hypothesis testing, the left segment is achieved by classical PCA and the right two segments are achieved by the modified Sparse Aggregation \eqref{SAhyp}. 

Next, we present a CTUB for each of the three problems: 
\begin{align*}
& \tilde{\eta}^{clu}_{\theta}(\beta)  = 
\left\{
\begin{array}{ll}
(1+\theta - 2 \beta)/4,	 	& \qquad   \beta < 1/2,     \\
\theta/4, 				& \qquad  1/2 < \beta < 1 - \theta/2,     \\
(1 - \beta)/2,	&  \qquad  \beta > (1 - \theta).  
\end{array}
\right.\cr
& 
\tilde{\eta}^{sig}_{\theta}(\beta) = 
\left\{
\begin{array}{ll}
\theta/2,		 				& \qquad  \beta < (1 - \theta)/2,  \\
(1 + \theta - 2 \beta)/4,   &\qquad   (1 - \theta)/2  <  \beta < 1/2,  \\ 
\theta/4,	& \qquad  \beta >  1/2.   
\end{array}
\right.\cr
& \tilde{\eta}^{hyp}_{\theta}(\beta)  = 
\left\{
\begin{array}{ll}
(1+\theta - 2 \beta)/4,	 	& \qquad   \beta < 1/2,     \\
\theta/4, 				& \qquad  \beta> 1/2.  
\end{array}\right. 
\end{align*}
See Figure~\ref{fig:complimit-sign} (top right and the two bottom panels). 

{\it Methods associated with CTUB}. The CTUB for clustering is associated with the methods of classical PCA (left segment) and IF-PCA (right two segments). The CTUB for signal recovery is associated with the methods of classical PCA (left segment) and IF-PCA (right segment). The CTUB for hypothesis testing is associated with the methods of classical PCA (left segment) and IF-PCA (right segment). %\footnote{See Section~\ref{sec:hyp} for the Higher Criticism test. IF-PCA also achieves this part of boundary.}  

{\bf Remark.} We now make a connection to the recent literature on the Gaussian mixture learning (e.g.\cite{WassermanSparse,chaudhuri2009}). In our framework, we calibrate with $(\eps, \tau)$. In the latter, we calibrate with $\|\mu\|$ and $\|\mu\|_0$.  For brevity,  we only discuss  the problem of hypothesis testing. 
The statistical limits for hypothesis testing can be (roughly) re-stated as follows: 
\begin{itemize} 
	\item $\sqrt{np} \ll s \ll p$:  $s \tau^2 = \sqrt{p/n}$.  
	\item $n \ll s \ll  \sqrt{np}$: $\tau = n^{-1/2}$.   \item $s \ll n $:  $\tau = (sn)^{-1/4}$.  
\end{itemize} 
Note that the first item corresponds to the {\it non-sparse} cases in the Gaussian mixture learning literature, where   $\|\mu\|^2 = s \tau^2 = \sqrt{p/n}$; the results match with those in, e.g.,   \cite{WassermanSparse,chaudhuri2009}.   The second one is part of  the {\it sparse} case in the Gaussian mixture learning literature, where $1 \ll  \|\mu\|^2=p/n \ll \sqrt{p/n}$  and $n \ll \|\mu\|_0 \ll \sqrt{np}$. The last one is also part of the sparse case, where  $\|\mu\|^2  = \sqrt{s/n}$ and $s \ll n$.

%%%%%%%%%%%
\section{Proof of Lemmas in Section 2}\label{app:A}
In this sectoin, we prove the post-selection random matrix theory results in Section~\ref{sec:IF-PCA}, specifically Lemmas~\ref{lem:postZ-1}--\ref{lem:v'H0v}. 

\subsection{Preliminary lemmas for Section~\ref{sec:IF-PCA}}  \label{subsec:prelim}
Lemma~\ref{lem:Bernstein} states the well-known Bernstein inequality \cite{Wellner86}.  Lemma~\ref{lem:vershynin} is a result from classical Random Matrix Theory \cite[Page 21]{Vershynin}. Lemma~\ref{lem:chi-square} states some properties about columns of the matrix $Z^{(q)}$; it is proved in Section~\ref{proof:chi-square}. 
%%%%%%%%%%%
%%%%%%%%%%%
\begin{lemma}  \label{lem:Bernstein}
	Let $X_1,\cdots, X_N$ be independent random variables with $E[X_k]=0$ and $\mathrm{var}(X_k)\leq v_k$, for $1\leq k\leq N$. Suppose $E(|X_k|^m)\leq v_km!c^{m-2}/2$ for all $m\geq 2$, where $c>0$ is a constant. Then for all $\lambda>0$,
	\[
	P\big( \big|\sum_{k=1}^N X_k\big| \geq \lambda\sqrt{N}\big) \leq 
	\exp\left( - \frac{\lambda^2/2}{\sum_{k=1}^n v_k/N + c\lambda/\sqrt{N}} \right).
	\]
\end{lemma}
%%%%%%%%%%%%

%%%%%%%%%%%%
%%%%%%%%%%%%
\begin{lemma} \label{lem:vershynin}
	Let $A$ be an $N \times n$ matrix whose entries are independent standard normal random variables. Then for every $x \geq 0$, with probability at least $1 - 2\exp(- x^2/2)$, 
	\[
	\sqrt{N} - \sqrt{n} - x \leq s_{\min}(A) \leq s_{\max}(A) \leq \sqrt{N} + \sqrt{n} + x,
	\]
	where $s_{\min}(A)$ and $s_{\max}(A)$ are the respective minimum and maximum singular values of $A$. 
\end{lemma}
%%%%%%%%%%%%%

Fix $q>0$. With $e_1=(1,0,\cdots,0)'$ and $z\sim N(0, I_p)$, we introduce a few notations: 
\begin{align*}
& \pi_0^{(q)} = P( \|z\|^2>n+2\sqrt{qn\log(p)} ),\cr
& \pi_1^{(q)} = P( \|z + \sqrt{n}\tau_p^*e_1 \|^2>n+2\sqrt{qn\log(p)}),\cr
& a_p^{(q)} = E\big[ (z(1))^2\cdot 1\{ \|z \|^2> n+2\sqrt{qn\log(p)}\}\big], \cr
& b_p^{(q)} = E\big[ (z(2))^2\cdot 1\{ \|z + \sqrt{n}\tau_p^* e_1\|^2> n+2\sqrt{qn\log(p)}\}\big],\cr
& c_p^{(q)} = E\big[ (z(1))^2\cdot 1\{ \|z + \sqrt{n}\tau_p^* e_1\|^2> n+2\sqrt{qn\log(p)}\}\big].
\end{align*}
For notation simplicity, we omit all the superscripts. In the following lemma, $m^{(q)}(\ell,\mu)$, $m_*^{(q)}(\ell,\mu)$ and the event $D_p$ are defined in Section~\ref{sec:IF-PCA}.   
%%%%%%%%%%%%
%%%%%%%%%%%%
\begin{lemma}  \label{lem:chi-square}
	Let $S(\mu)$ denote the support of $\mu$, $\kappa_m=E(|z(1)|^m)$ and $\kappa_{2m}(n)=E(\|z\|^{2m})$, where $z\sim N(0, I_n)$. Below, all the probabilities are conditioning on $(\ell,\mu)$, and the $o(1)$ terms are uniform for all realizations of $(\ell,\mu)$ in the event $D_p$.   
	\begin{itemize}
		\item[(a)] Fix $j\notin S(\mu)$. For any $v\in \mathcal{S}^{n-1}$ and any integer $m\geq 1$
		\[
		\begin{array}{l}
		E[(v'z_j^{(q)})^2]=a_p,\\
		E(|v'z_j^{(q)}|^m) \leq \kappa_m \pi_0(1+o(1)),\\
		E(\|z_j^{(q)}\|^2)=na_p,\\
		E(\|z_j^{(q)}\|^{2m})=\kappa_{2m}(n) \pi_0(1+o(1)),\\
		a_p = \pi_0(1+L_pn^{-1/2}).    
		\end{array}
		\]
		\item[(b)] Fix $j\in S(\mu)$. For any $v\in \mathcal{S}^{n-1}$ and any integer $m\geq 1$
		\[
		\begin{array}{l}
		E[(v'z_j^{(q)})^2]=b_p + (c_p - b_p)\tfrac{(v'\ell)^2}{\|\ell\|^2}\\
		E(|v'z_j^{(q)}|^m) \leq \kappa_m \pi_1(1+o(1)),\\
		E(\|z_j^{(q)}\|^2)= nb_p + (c_p - b_p),\\ 
		E(\|z_j^{(q)}\|^{2m})\leq 2^m\kappa_{2m}(n) \pi_1(1+o(1)),\\
		b_p=\pi_1(1+L_pn^{-1/4}), \;\; c_p=\pi_1(1+L_pn^{-1/4}).     
		\end{array}
		\]
		\item[(c)] $m^{(q)}(\ell,\mu)=(p-|S(\mu)|)\pi_0 + |S(\mu)|\pi_1$, \\
		$m_*^{(q)}(\ell,\mu)=(p-|S(\mu)|)a_p + |S(\mu)|[b_p+n^{-1}(c_p-b_p)]$. 
	\end{itemize}
\end{lemma}

\noindent
{\bf Remark}. Lemma~\ref{lem:chi-square} allows us to characterize the quantities $m^{(q)}$ and $m^{(q)}_*$. First, by (a)-(b), $a_p\sim \pi_0$ and $b_p\sim c_p\sim \pi_1$. Combining them with (c) gives that $m_*^{(q)}\sim m^{(q)}$. Second, we look at $m^{(q)}$. By Lemma~\ref{lem:chi-square-0} and Mills' ratio \cite{Wellner86}, $\pi_0\sim \bar{\Phi}(\sqrt{2q\log(p)})=L_pp^{-q}$. Similarly, $\pi_1\sim L_p p^{-[(\sqrt{r}-\sqrt{q})_+]^2}$. Plugging them into (c) gives 
\[
m^{(q)}\sim L_p p^{1-q} + p\epsilon_p \cdot L_p p^{-[(\sqrt{r}-\sqrt{q})_+]^2}. 
\]
This is the equation \eqref{If-pca-asy} in the main text of \cite{3Phase}.

%%%%%%%%%%%%%%
%%%%%%%%%%%%%%
\subsection{Proof of Lemma~\ref{lem:postZ-1}}
Fix $q>0$ and write $H_0=Z^{(q)}(Z^{(q)})'$ and $\hat{S}=\hat{S}^{(if)}_q$ for short. Fix a realization $(\ell,\mu)$. With probability at least $1-O(p^{-3})$, 
\beq \label{thm-postZ-1}
||\hat{S}| - m^{(q)}| \leq \sqrt{6m^{(q)}\log(p)} . 
\eeq 

First, we consider $q> \tilde{q}(\beta,\theta,r)$, so that $m^{(q)}\leq n p^{-\delta}$ for some $\delta>0$.  
Let $k = \lceil m^{(q)} + \sqrt{6m^{(q)}\log(p)} \rceil$. Under \eqref{thm-postZ-1}, 
\begin{align*}
& \lambda_{\max}(H_0) \leq \max_{T\subset\{1,\cdots, p\}, |T|\leq k} \lambda_{\max}((ZZ')^{T,T}),\cr
& \lambda_{\min}^+(H_0) \geq \min_{T\subset\{1,\cdots, p\}, |T|\leq k} \lambda_{\min}^+((ZZ')^{T,T}),
\end{align*}
where for a matrix $A$, $\lambda_{\min}^+(A)$ denotes the minimum non-zero eigenvalue and $A^{T,T}$ is the submatrix restricted to rows and columns in $T$. 
For each fixed $T$, we can write $(Z'Z)^{T,T}=Z_T(Z_T)'$, where $Z_T=(z_j, j\in T)$ is an $n\times |T|$ matrix with {\it iid} entries of $N(0,1)$. Using Lemma~\ref{lem:vershynin},  for each $T$, with probability at least $1-O(p^{-(k+3)})$, all non-zero eigenvalues of $(ZZ')^{T,T}$ fall into 
\begin{align*}
\Big[ \big(\sqrt{n} + \sqrt{|T|}  + \sqrt{6k\log(p)}\big)^2, \; \big(\sqrt{n} - \sqrt{|T|} - \sqrt{6k\log(p)}\big)^2 \Big]= n\pm C\sqrt{nk\log(p)}.   
\end{align*}
Note that the number of subsets $T$ such that $|T|\leq k$ is no more than $p^k$. Combining the above results, we find that with probability at least $1-O(p^{-3})$, all non-zero eigenvalues of $H_0$ fall into 
\beq  \label{thm-postZ-2}
n \pm C\sqrt{nm^{(q)}\log(p)}. 
\eeq
The claim then follows. 

Next, we consider $q< \tilde{q}(\beta,\theta,r)$, so that $m^{(q)}\geq n p^{\delta}$ for some $\delta>0$.  Write for short 
\[
\omega_p = \sqrt{nm^{(q)}}+o(1) |S(\mu)|\pi_1, 
\]
where $\pi_1$ is as in Lemma~\ref{lem:chi-square} and $m_1^{(q)}=|S(\mu)|\pi_1$ by definition. It suffices to show that with probability at least $1-O(p^{-3})$, 
\beq \label{thm-postZ-3}
\|H_0 - m^{(q)}_* I_n\| \leq C\omega_p. 
\eeq
%Then the claim follows from Lemma~\ref{lem:chi-square} where we see that $m_1^{(q)}=|S(\mu)|\pi_1$. 
%By elementary statistics (Mills' ratio), $n^{-5/4}\|\mu\|_0\pi_1=L_pp^{-5\theta/4+1-\beta-[(\sqrt{q}-\sqrt{r})_+]^2}$ and $\sqrt{nm^{(q)}}\geq n=p^{\theta}$. It is seen that $\omega_p\sim \sqrt{nm^{(q)}}$ under our assumptions. Therefore,  the claim follows directly once \eqref{thm-postZ-3} is proved. 

We now show \eqref{thm-postZ-3}. Fix $\alpha>0$. A subset $\mathcal{M}_\alpha$ of the unit sphere $\mathcal{S}^{n-1}$ is called an $\alpha$-net if for any $v \in \mathcal{S}^{n-1}$, there exits $u \in \mathcal{M}_\alpha$ such that $\|u-v\| \leq \alpha$. 
The following lemma states some well-known results and its proof can be found in 
\cite[Page 8]{Vershynin}.
%%%%%%%%%%%%%%%%
\begin{lemma}\label{lem:alpha-net}
	Fix $\alpha\in (0,1/2)$. For any $\mathcal{M}_\alpha$, an $\alpha$-net of $\mathcal{S}^{n-1}$, and  any symmetric matrix $A\in R^{n,n}$,
	$\|A\| \leq (1 - 2 \alpha)^{-1} \sup_{u \in \mathcal{M}_\alpha} \{ |u' A u| \}$. Moreover, there exists an $\alpha$-net $\mathcal{M}^*_\alpha$ of $\mathcal{S}^{n-1}$ such that  $|\mathcal{M}^*_\alpha |\leq (1+2/\alpha)^n$. 
\end{lemma}
%%%%%%%%%%%%%%%%%
\noindent 
By Lemma~\ref{lem:alpha-net} with $\alpha=1/4$ , there exists a subset $\mathcal{M}^*$, such that $|\mathcal{M}^*|\leq 9^n$ and $\sup_{v\in \mathcal{M}^*} v'Av\geq \|A\|/2$ for any $n\times n$ matrix $A$. Therefore, 
to show the claim, it suffices to show that for each fixed $v\in \mathcal{M}^*$, with probability $\geq 1-O(9^{-n}p^{-3})$, 
\beq  \label{thm-postZ-4}
|v'(H_0 - m^{(q)}_* I_n)v| \leq C \omega_p.  
\eeq 

We now show \eqref{thm-postZ-4}. 
Fix $v$ and define 
\[
W_j=(v'z_j^{(q)})^2 - a_p, \;\; \mbox{for $j\notin S(\mu)$};  \qquad W_j=(v'z_j^{(q)})^2 - b_p,\;\;\mbox{for $j\in S(\mu)$},
\]
where $a_p$ and $b_p$ are defined in Section~\ref{subsec:prelim}. By (c) of Lemma~\ref{lem:chi-square}, $m_*^{(q)}=(p-|S(\mu)|)a_p+|S(\mu)|b_p+n^{-1}(c_p-b_p)|S(\mu)|$. Since $|c_p-b_p|=o(\pi_1)$,   we can rewrite 
\beq
v'(H_0-m^{(q)}_* I_n)v = \sum_{j=1}^p W_j + o(n^{-1}|S(\mu)|\pi_1). 
\eeq
Here $W_j$'s are independent of each other. Applying Lemma~\ref{lem:chi-square}, we get the following results. For $j\notin S(\mu)$, $E(W_j)=0$, $\mathrm{var}(W_j)\leq 3\pi_0(1+o(1))$ and $E(|W_j|^m)\leq \kappa_{2m} \pi_0 (1 + o(1))$. For $j\in S(\mu)$, $|E(W_j)| \leq |b_p-c_p|=\pi_1\cdot o(1)$, $\mathrm{var}(W_j)\leq 3\pi_1(1+o(1))$ and $E(|W_j|^m)\leq \kappa_{2m} \pi_1(1 + o(1))$. So we have
\[
|\sum_{j=1}^p E(W_j)| =o(1)|S(\mu)|\pi_1, \qquad \sum_{j=1}^p \mathrm{var}(W_j)\lesssim 3m^{(q)}. 
\]
We apply Lemma~\ref{lem:Bernstein} with $\lambda=\sqrt{9p^{-1}m^{(q)}(n\log(9)+2\log(p)+\log(2))}$. To check the moment conditions, we note that $\kappa_{2m}=E_{Y\sim N(0,1)}(|Y|^{2m})\leq 2^m m!$ for all $m\geq 1$. 
Furthermore, since $m^{(q)}/n\goto \infty$, we have $\sum_j\mathrm{var}(W_j)/p\sim 3m^{(q)}/p\gg \lambda/\sqrt{p}$. 
It follows that with probability $\geq 1-O(9^{-n}p^{-3})$,
\[
|\sum_{j=1}^p W_j|\lesssim  3\sqrt{\log(9)} \sqrt{nm^{(q)}} + o(1)|S(\mu)|\pi_1. 
\]
This gives \eqref{thm-postZ-4}, and the proof is now complete.\qed

%%%%%%%%%%%%%%
%%%%%%%%%%%%%%
\subsection{Proof of Lemma~\ref{lem:postZ-2}}
We have shown the first claim in \eqref{thm-postZ-3}, noting that $\omega_p\sim \sqrt{nm^{(q)}}$ when $r<\rho^*_\theta(\beta)$ (see also \eqref{rate}). 

We now show the second claim. Write for short $H_0=Z^{(q)}(Z^{q})'$. The key is the following lemma, which is proved in Section~\ref{app:C}. 
%%%%%%%%%%%
%%%%%%%%%%%
\begin{lemma}  \label{lem:Frobenius}
	Under conditions of Lemma~\ref{lem:postZ-2}, as $p\to\infty$, conditioning on any realization of $(\ell,\mu)$ on the event $D_p$, with probability at least $1-O(n^{-2})$, 
	\begin{align*}
	& |n^{-1}\mathrm{tr}(H_0)-m_*^{(q)}| \leq C\sqrt{m^{(q)}\log(p)}, \cr
	%\leq p^{-\delta} \sqrt{nm^{(q)}}, \quad \mbox{for some } \delta>0, \cr
	& \|H_0\|^2_F \geq n^{-1}[\mathrm{tr}(H_0)]^2 + Cn^2 m^{(q)}. 
	\end{align*}
\end{lemma}

%%%%%%%%%%%

%For this claim to hold, we need an extra assumption: $n^{-7/4}\|\mu\|_0\pi_p^{(1)}(\ell)/[N_p^*(\mu)]^{1/2}= L_pn^{-\alpha}$ for some $\alpha>0$. It translates to
%\beq  \label{thm-postZ-10}
%L_pn^{-5/4}\|\mu\|_0\pi_p^{(1)} = o(1)\cdot \sqrt{nN_p^*(\mu)}, \qquad \delta_p= \sqrt{nN_p^*(\mu)}(1+o(1)). 
%\eeq
%
%To prove \eqref{thm-postZ-5}, the key is to show 
%\beq  \label{thm-postZ-6}
%n^{-1}|\mathrm{tr}(H_0)-N_p^*(\mu)| = o(\delta_p). 
%\eeq
%and
%\beq  \label{thm-postZ-7}
%\|H_0\|^2_F \geq n^{-1}[\mathrm{tr}(H_0)]^2 + C\delta_p^2.  
%\eeq
%Suppose \eqref{thm-postZ-6}-\eqref{thm-postZ-7} are true. 

Let $k$ be the largest integer that is no larger than $m^{(q)}/2$. Since
$k\gg n$, for each fixed $k\times k$ submatrix of $ZZ'$, its rank is $n$ with probability $1$. Using \eqref{thm-postZ-1}, the rank of $H_0$ is $n$ with probability at least $1-O(p^{-3})$. 
Let $\lambda_1\geq \lambda_2\geq \cdots \geq \lambda_{n}>0$ be the eigenvalues of $H_0$ and write $\bar{\lambda}=n^{-1}\sum_{i=1}^n\lambda_i$. For $\delta_p\equiv \sqrt{nm^{(q)}}$, \eqref{thm-postZ-3} and Lemma~\ref{lem:Frobenius} imply
\beq  \label{thm-postZ-6}
|\lambda_1-\lambda_n|\leq A_1\delta_p, \qquad \bar{\lambda}=m_*^{(q)}+o(\delta_p), \qquad \sum_{i=1}^n \lambda_i^2 \geq n\bar{\lambda}^2 + A_2n\delta_p^2, 
\eeq
for some constants $A_1,A_2>0$. On one hand, 
\[
\sum_{i=1}^n (\lambda_i-\bar{\lambda})^2 \geq A_2n\delta_p^2. 
\]
On the other hand, $\lambda_i-\bar{\lambda}\leq \lambda_1-\bar{\lambda}$ for $i$ satisfying $\lambda_i\geq \bar{\lambda}$; moreover, $\bar{\lambda}-\lambda_i\leq A_1\delta_p$ for $i$ such that $\lambda_i<\bar{\lambda}$. It follows that
\begin{align*} 
\sum_{i=1}^n (\lambda_i-\bar{\lambda})^2 &\leq  (\lambda_1-\bar{\lambda})\sum_{i:\lambda_i\geq \bar{\lambda}}(\lambda_i-\bar{\lambda}) + A_1\delta_p \sum_{i:\lambda_i< \bar{\lambda}} (\bar{\lambda}-\lambda_i)\cr
& = [(\lambda_1-\bar{\lambda}) + A_1\delta_p] \sum_{i:\lambda_i\geq \bar{\lambda}}(\lambda_i-\bar{\lambda})\cr
&\leq n [(\lambda_1-\bar{\lambda}) + A_1\delta_p] (\lambda_1-\bar{\lambda}). 
\end{align*}
Now, if we write $x=\lambda_1-\bar{\lambda}$, then $x(x+A_1\delta_p)\geq A_2\delta_p^2$. It follows that
\beq
\lambda_1-\bar{\lambda} \geq \frac{\sqrt{A_1^2+4A_2}-A_1}{2}\delta_p. 
\eeq
Combining it with the second equation in \eqref{thm-postZ-6}, we obtain that $\lambda_1\geq m^{(q)}_*+C\sqrt{nm^{(q)}}$. \qed

\subsection{Proof of Lemma~\ref{lem:matrixA}}
Write for short $A=\ell(Z^{(q)}\mu^{(q)})' + (Z^{(q)}\mu^{(q)})\ell'$. Since 
\[
\|A\|\leq 2\|\ell\|\|Z^{(q)}\mu^{(q)}\|\leq 2n\|Z^{(q)}\mu^{(q)}\|_\infty,
\]
it suffices to show that with probability $1-O(p^{-3})$, 
\beq \label{lem-matrixA-1}
\|Z^{(q)}\mu^{(q)}\|_\infty \leq C\tau_p^*\sqrt{m_1^{(q)}}.
\eeq

Note that
\beq \label{lem-matrixA-2}
\|Z^{(q)}\mu^{(q)}\|_\infty = \max_{1\leq i\leq n} |\sum_{j\in S(\mu)}\mu(j)z_j^{(q)}(i)|
=\tau_p^* \max_{1\leq i\leq n} |\sum_{j\in S(\mu)}z_j^{(q)}(i)|.
\eeq
Fix $i$ and write $V_j=z_j^{(q)}(i)$ for short. Then $V_j$'s are independent and $E(V_j)=0$ by symmetry. We apply Lemma~\ref{lem:chi-square} with $v=e_1$ and find that $\mathrm{var}(V_j)\leq b_p+2|c_p-b_p|=\pi_1(1+o(1))$. By Lemma~\ref{lem:Bernstein} (the moment conditions can be verified using Lemma~\ref{lem:chi-square}), $|\sum_{j\in S(\mu)}V_j|\leq 2\sqrt{|S(\mu)|\pi_1}$ with probability $1-O(p^{-4})$. It follows that with probability $1-O(p^{-3})$, 
\beq \label{lem-matrixA-3}
\max_{1\leq i\leq n} |\sum_{j\in S(\mu)}z_j^{(q)}(i)| \leq  C \sqrt{|S(\mu)|\pi_1}=C\sqrt{m_1^{(q)}}. 
\eeq
Combining \eqref{lem-matrixA-2}-\eqref{lem-matrixA-3} gives \eqref{lem-matrixA-1}. \qed

%%%%%%%%%%%
%%%%%%%%%%%
\subsection{Proof of Lemma~\ref{lem:v'H0v}}
Introduce
\begin{align*} \label{DefineDelta}
& \Delta^{\dag}(q, \beta, r, \theta)\cr
= & \left\{\begin{array}{ll}
\beta - \frac{1}{2}\min\{q, \beta - \frac{\theta}{2}\}, & q \leq r,    \\
\beta + (\sqrt{q} - \sqrt{r})^2 - \frac{1}{2} \min\{q, \beta - \frac{\theta}{2} + (\sqrt{q} - \sqrt{r})^2\}, & q > r.  
\end{array}
\right. 
\end{align*} 
By elementary algebra, 
\[
r<\rho^*_\theta(\beta) \qquad \Longleftrightarrow \qquad \min_{q>0}\Delta^{\dag}(q,\beta,r,\theta)>1/2,
\]
Moreover, using Mills' ratio, 
\[
\min\left\{ 
\frac{n(\tau_p^*)^2 |S(\mu)|\pi_1}{\sqrt{ nm^{(q)}}}, \;\; 
\tau_p^*\sqrt{|S(\mu)|\pi_1} \right\}  
= L_pp^{1/2-\Delta^{\dag}(q,\beta,r,\theta)}. 
\]
It follows that for some $\delta>0$, 
\beq \label{rate}
r<\rho^*_\theta(\beta) \qquad \Longleftrightarrow \qquad |S(\mu)|\pi_1\leq p^{-\delta}
\max \{\sqrt{m^{(q)}}, \sqrt{n} \}. 
\eeq

Consider the first claim. The proof is similar to that of \eqref{thm-postZ-4}, except that we take $\lambda= C\sqrt{p^{-1}m^{(q)}\log(p)}$ when applying Lemma~\ref{lem:Bernstein}. 
It follows that for any $v\in \mathcal{S}^{n-1}$, with probability at least $1-O(p^{-3})$,
\[
|v'(H_0-m^{(q)}_* I_n)v| \leq C\sqrt{m^{(q)}\log(p)} + o(|S(\mu)|\pi_1). 
\]
By \eqref{rate}, the second term above is negligible and the claim follows. 

Consider the second claim. $H-H_0=\|\mu^{(q)}\|^2\ell\ell' + A$, where with probability at least $1-O(p^{-3})$, $\|A\|\leq Cn\tau_p^*\sqrt{m_1^{(q)}}$ by Lemma~\ref{lem:matrixA}. Furthermore, by elementary statistics, $\|\mu^{(q)}\|^2\leq Cm_1^{(q)}(\tau_p^*)^2$. Note that $m_1^{(q)}=|S(\mu)|\pi_1$ due to the spherical symmetry of $N(0, I_n)$.  Together, we see that
\[
\|H-H_0\|\leq L_p(\sqrt{n} |S(\mu)|\pi_1 + n^{3/4} \sqrt{|S(\mu)|\pi_1}). 
\]
If $q>\tilde{q}(\beta,r,\theta)$, then $m^{(q)}=o(n)$ and \eqref{rate} implies $\|H-H_0\|\leq p^{-\delta}n$. 
If $q<\tilde{q}(\beta,r,\theta)$, then $n=o(m^{(q)})$ and \eqref{rate} implies  $\|H-H_0\|\leq p^{-\delta}\sqrt{nm^{(q)}}$. \qed

%%%%%%%%
%%%%%%%%
%%%%%%%%
\section{Proof of Lemmas in Section 3}  \label{app:B}
In this section, we prove Lemmas~\ref{lemma:PCA}--\ref{lem:u-prop}.

\subsection{Proof of Lemma \ref{lemma:PCA}}  

We first show the claim for $B=I_p$ and then generalize it to any $B$ satisfying $\max\{\|B\|, \|B^{-1}\|\}\leq L_p$. 

Fix $B=I_p$. We use $\delta > 0$ to denote a generic constant which only depends on $(\alpha, \beta, \theta)$ but may change from occurrence to occurrence.  
In our model, $X = \ell \mu' + Z$. Let  $H_0 = ZZ' - p I_n$. It is seen 
\begin{equation} \label{XX'} 
XX' - p I_n   = [\|\mu\|^2 \ell \ell' +  \ell \mu' Z' + Z \mu \ell']  + ZZ' - p  I_n = [\|\mu\|^2 \ell \ell' +  \ell \mu' Z' + Z \mu \ell'] + H_0. 
\end{equation} 
Since $\xi$ is a left singular vector of $X$, 
$\lambda \xi =   [ \|\mu\|^2 (\xi, \ell) +   (\xi, Z \mu) ] \ell  + (\xi, \ell) Z \mu + H_0 \xi$.   Rearranging it, we have   
\begin{equation} \label{rootnxi} 
\sqrt{n} \xi =   (I_n - (1/\lambda) H_0)^{-1} [b_1 \ell + b_2 Z (\mu / \|\mu\|)], 
\end{equation} 
where $b_1= b_1(\ell, Z, \mu) =  (1/\lambda) \cdot [ \sqrt{n} \|\mu\|^2 (\xi, \ell) + \sqrt{n} (\xi, Z \mu)]$ and $b_2 = b_2(\ell, Z, \mu) = (1/\lambda) \sqrt{n} \|\mu\|  (\xi, \ell)$.  
Therefore,  $\min\{\|\sqrt{n} \xi - \ell \|_{\infty}, \|\sqrt{n} \xi + \ell \|_{\infty} \}$ is no greater than 
\begin{equation} \label{toshow0}
\min\{| b_1  - 1|, |b_1 + 1|\} + |b_1| \|\ell -  (I_n - (1/\lambda) H_0)^{-1} \ell \|_{\infty} +  |b_2| \| (I_n - (1/\lambda) H_0)^{-1}  Z (\mu /\|\mu\|) \|_{\infty}. 
\end{equation} 
To show the claim, it is sufficient to show that  with probability at least $1 - o(p^{-3})$, 
\begin{equation} \label{toshow1} 
\min\{|b_1 - 1|, |b_1 + 1|\}   \leq p^{-\delta}, \qquad |b_2| \leq    p^{-\delta}, 
\end{equation} 
and 
\begin{equation} \label{toshow2} 
\| \ell - (I_n - (1/\lambda) H_0)^{-1} \ell \|_{\infty} \leq p^{-\delta},   \qquad  \| (I_n - \frac{1}{\lambda} H_0)^{-1} Z (\mu / \|\mu\|)   \|_{\infty} 
\leq  C \sqrt{\log(p)}.    
\end{equation}  

We now show (\ref{toshow1}).  
Consider the first item.  Since $Z$ and $\mu$ are independent, we have that with probability at least $1 - o(p^{-3})$, $|(\xi, Z \mu)| \leq  \|\mu\| \cdot  \|Z (\mu / \| \mu\|) \| \leq 2 \|\mu\|  \sqrt{n}$. Combining this with the triangle inequality,    
\begin{align} 
& \min\{b_1- 1,   b_1 + 1\}   \nonumber \\
\leq &  (n \|\mu\|^2 /\lambda) |\ang(\ell, \xi) - 1|  + |1 - (n \|\mu\|^2 / \lambda)| + (\sqrt{n}/\lambda)  |(\xi, Z \mu)|  \nonumber  \\\leq &  (n \|\mu\|^2 /\lambda) |\ang(\ell, \xi) - 1|  + |1 - (n \|\mu\|^2 / \lambda)|  +  2 n\|\mu\|  /\lambda. \label{PCApf1}   
\end{align} 
At the same time, we rewrite (\ref{XX'}) as   
\begin{equation} \label{XX'2} 
XX' - p I_n = A  + H_0, \;  \mbox{where $A =  \|\mu\|^2 \ell \ell' +  \ell \mu' Z' + Z \mu \ell'$ for short}.
\end{equation}  
Note that $A$ is a symmetric matrix of rank $2$. For short, write 
$\nu = \|\mu\|^{-2} \mu$ and $a = a(\ell, \mu, Z) = (1 + 4 n^{-1} [ \ell' Z \nu   + \| Z \nu \|^2])^{1/2}$. 
Let $\lambda_{\pm}$ be the two nonzero eigenvalues of $A$, and let $\eta_{\pm}$ be the corresponding eigenvectors. By elementary algebra,   
\begin{equation} \label{lambdas} 
\lambda_\pm(A) = n \|\mu\|^2[(1/2) (1 \pm a)   +  n^{-1} \ell' Z \nu], \qquad  \eta_{\pm}  \propto   (1/2) (1 \pm a)  \ell +     Z \nu. 
\end{equation} 
By elementary statistics, it is seen that with probability at least $1 - o(p^{-3})$ that $n^{-1} [ |\ell' Z \nu|  + \| Z \nu\|^2]$ does not exceed  
\begin{equation}\label{lznu}
C \sqrt{\log(p)}  n^{-1}  [(\sqrt{n} \|\mu\|^{-1})    +  n  \|\mu\|^{-2}]  =  C \sqrt{\log(p)} [(\sqrt{n} \|\mu\|)^{-1} + \|\mu\|^{-2}]. 
\end{equation}
Note that for $(\alpha, \beta, \theta)$ in our range of interest, $p \eps_p  \tau_p^2 \geq \sqrt{p/n} = p^{(1 - \theta)/2}$. 
By the way $\mu$ is generated, $\|\mu\|^2 \sim p \eps_p \tau_p^2$. Therefore, with probability at least $1 - o(p^{-3})$, 
\begin{equation} \label{PCApf2} 
\|\mu\|^2 \sim  p \eps_p \tau_p^2 \geq p^{\delta} \sqrt{p/n}, \qquad n \|\mu\|^2 \sim n p \eps_p \tau_p^2 \geq p^{\delta} \sqrt{pn}. 
\end{equation} 
Inserting (\ref{PCApf2}) into (\ref{lznu}) gives that  with  probability at least  $1 - o(p^{-3})$,  
$|a-1|  \leq  C p^{-\delta}$.  %, and $|(\ell / \|\ell\|)' Z \nu| \leq   p^{-\delta}$.  
Combining this with (\ref{lambdas}),  
\begin{equation} \label{lambdas2} 
|(n \|\mu\|^2 / \lambda_+) - 1| \leq p^{-\delta},  \qquad  (\lambda_{-}/\lambda_+)  \leq p^{-\delta},   \qquad |\ang(\ell, \eta_+) - 1| \leq p^{-\delta}. 
\end{equation} 
At the same time,  by a direct use of the elementary Random Matrix Theory \cite{Vershynin}, $\|H_0 \|  = \| ZZ' - p I_p\| \leq C \sqrt{pn}$. Combining these with (\ref{PCApf2})-(\ref{lambdas2}) gives  
\begin{equation} \label{PCApf4}
\|(1/\lambda_+) H_0\| \leq C \sqrt{pn} / (n \|\mu\|^2)  \leq C p^{-\delta}. 
\end{equation} 
This says that in (\ref{XX'2}), the leading eigenvalue of $A$ is larger than that of $H_0$ by $p^{\delta}$ times. 
By matrix perturbation theory, we have that  with probability at least $1 - o(p^{-3})$, 
\begin{equation} \label{PCApf3} 
|\lambda_+ / \lambda -1|  \leq p^{-\delta}, \qquad  |\ang(\eta_+, \xi) - 1| \leq p^{-\delta}. 
\end{equation} 
Combining (\ref{lambdas2}) and (\ref{PCApf3})  gives 
\begin{equation} \label{PCApf5}
|(n \|\mu\|^2 / \lambda)  -1|  \leq p^{-\delta}, \qquad  |\ang(\ell, \xi) - 1| \leq p^{-\delta}. 
\end{equation} 
In particular, combining  (\ref{PCApf2}), (\ref{PCApf4}),  and (\ref{PCApf3}) gives that with probability at least $1 - o(p^{-3})$,   
\begin{equation} \label{PCApf6}
\|(1/\lambda) H_0\|  \leq   C p^{-\delta}, \qquad \sqrt{pn} / \lambda \leq p^{-\delta}. 
\end{equation} 
Inserting (\ref{PCApf5}) into (\ref{PCApf1})  gives the first item of (\ref{toshow1}).  

Consider the second item of (\ref{toshow1}).  Note that $|b_2| \leq (n\|\mu\|/\lambda)$,  where by (\ref{PCApf5}), the right hand side $\leq \|\mu\|^{-1}$. The claim follows directly from (\ref{PCApf2}). 

We now show (\ref{toshow2}).  Since the proofs are similar, we only show the first item.   Let $e_1$ be the first base vector of $R^n$. 
Note that by symmetry and by using the union bound, it is sufficient to show that  with probability at least $1 - o(p^{-4})$,  
\begin{equation} \label{PCApf10} 
\bigl| e_1' (I_n - \frac{1}{\lambda} H_0)^{-1} e_1 -1 \bigr| \leq  p^{-\delta}, \qquad  |e_1'  (I_n - \frac{1}{\lambda} H_0)^{-1}  (\ell - \ell_1 e_1) |      \leq C \sqrt{\log(p)}. 
\end{equation} 
The first claim follows easily by (\ref{PCApf6}) and basic algebra. For the second claim, write $\ell = (\ell_1, \tilde{\ell})'$, and 
let $\tilde{Z}$ be the $(n-1) \times p$ matrix consisting all but the first row of $Z$, and let $\tilde{H}_0 = \tilde{Z} \tilde{Z}' - p I_{n-1}$. It follows that 
\[
I_n -  (1/\lambda) H_0 = 
\left( 
\begin{array}{ll} 
1 - (1/\lambda) [\|Z_1\|^2 - p],  & - (1/\lambda) Z_1' \tilde{Z} \\
-(1/\lambda) \tilde{Z} Z_1,  &I_{n-1} - (1/\lambda) \tilde{H}_0 
\end{array}
\right), 
\]
and 
\begin{equation} \label{PCApf11}  
e_1' (I_n - \frac{1}{\lambda} H_0)^{-1}  (\ell - \ell_1 e_1)  =  (e_1' [I_n - (1/\lambda) H_0]^{-1} e_1) \cdot (1/\lambda)  Z_1'  \tilde{Z}'  [I_{n-1} - (1/\lambda) \tilde{H}_0]^{-1} \tilde{\ell}. 
\end{equation} 
Now,  since rows of $Z$ are independent, $Z_1$ and $\tilde{Z}[I_{n-1} - (1/\lambda) \tilde{H}_0]^{-1} \tilde{\ell}$ are two vectors that  almost independent of each other; the only issue is that $Z_1$ is correlated with $\lambda$. To overcome the difficulty, we write 
\begin{equation} \label{PCApf7} 
\frac{ Z_1'  \tilde{Z}'  [I_{n-1} - (1/\lambda) \tilde{H}_0]^{-1} \tilde{\ell} }{\lambda}=  \sum_{k = 0}^{\infty}  \frac{ Z_1' \tilde{Z} \tilde{H}_0^k \tilde{\ell} }{\lambda^{k+1}}= \sum_{k = 0}^{\infty} \frac{ \| \tilde{Z} \tilde{H}_0^{k} \tilde{\ell}\|}{\lambda^{k+1}} \cdot   \frac{Z_1' \tilde{Z} \tilde{H}_0^k \tilde{\ell}}{\| \tilde{Z} \tilde{H}_0^{k} \tilde{\ell}\|}.     
\end{equation} 
Now, for each $k$, $Z_1$ and $\tilde{Z} \tilde{H}_0^k \tilde{\ell}$ are independent, and so 
\[
Z_1'  (\tilde{Z}' \tilde{H}_0^k \tilde{\ell}  / \| \tilde{Z}' \tilde{H}_0^k \tilde{\ell} \|)   \sim N(0,1). 
\]
For $k$-th term, with probability $1 - o(p^{-4(k+1)})$, there is $|Z_1'  (\tilde{Z}' \tilde{H}_0^k \tilde{\ell}  / \| \tilde{Z}' \tilde{H}_0^k \tilde{\ell} \|)| \leq \sqrt{8(k+1) \log(p)}$. 
Additionally, by basics in RMT \cite{Vershynin}, with probability at least  $1 - o(p^{-4})$,  $\|\tilde{Z} \tilde{H}_0^k\| \leq \sqrt{p}   (C\sqrt{np})^k$ for all $k$. 
\[
(1/\lambda)^{k+1} \| \tilde{Z} \tilde{H}_0^k \tilde{\ell}\| \leq \sqrt{(n-1)} (1/\lambda)^k  \|\tilde{Z} \tilde{H}_0^k\| \leq  (C \sqrt{n p} / \lambda)^{k+1}.    
\]
Combining these with (\ref{PCApf7}) and the second term of (\ref{PCApf6}),  
it is seen  that with probability at least $1 - o(p^{-4})$, 
\[
|(1/\lambda)  Z_1'  \tilde{Z}'  [I_{n-1} - (1/\lambda) \tilde{H}_0]^{-1} \tilde{\ell}| \leq p^{-\delta}. 
\] 
Inserting this into (\ref{PCApf11}) and using the first item of (\ref{PCApf10}), the second item of (\ref{PCApf10}) follows. 

For a general $B$, the proof is similar by noting that $\|ZB\|\leq L_p \|Z\|$ and the following lemma, which is proved below. 
%%%%%%%%
\begin{lemma}  \label{lem:rmtcolor}
	As $n,p\goto\infty$ and $p/n\goto\infty$, for an $n \times p$ random matrix $Z$ where $Z(i,j) \overset{iid}{\sim} N(0, 1)$ and any non-random matrix $B\in R^{p,p}$ such that $\max\{ \|B\|, \|B^{-1}\|\} \leq L_p$, with probability $1 - O(p^3)$, $\|ZBB'Z' - \mathrm{tr}(BB') I_n\| \leq C\sqrt{n p}$. 
\end{lemma}
%%%%%%%%

%%%%%%%%%%%%
%%%%%%%%%%%%
%%%%%%%%%%%%
%%%%%%%%%%%%
\subsection{Proof of Lemma \ref{lemma:FN}}  
Letting $\Phi$ be the CDF of $N(0,1)$, denote the mean and variance of $|z_i + h|$ by $u(h)$ and $\sigma^2(h)$, respectively. It is seen that  
\begin{equation} \label{musigma} 
u(h)  =   \sqrt{2/\pi} e^{-h^2/2} + h[1 - 2 \Phi(-h)], \qquad   \sigma^2(h) = 1 + h^2 - \mu^2(h).  
\end{equation} 
By Jensen's inequality,   $E|z_i + h| \geq | E (z_i  + h)| = h$. It follows that 
\begin{equation} \label{elementary1}
u(h) \geq  h, \qquad  \sigma^2(h) \leq 1, 
\end{equation} 
At the same time, 
we claim that as $n \goto \infty$, for any $0  \leq x \leq \sqrt{n} / \log(n)$,   
\begin{equation} \label{FNadd}  
P \bigl(  \bigl| \sum_{i = 1}^n  \bigl( |z_i + h| - u(h) \bigr)  \bigr|  \geq \sqrt{n} x \bigr)  \leq 
2 \mathrm{exp}\bigl(-  (1 + o(1))  \frac{x^2}{2 \sigma^2(h)}  \bigr),   
\end{equation} 
where $o(1) \goto 0$ as $n \goto \infty$, uniformly for all $h > 0$ and $0 <  x \leq \sqrt{n} / \log(n)$. 
%\end{lemma} 
Combining (\ref{elementary1}) and (\ref{FNadd})  gives Lemma \ref{lemma:FN}. 
%\begin{equation} \label{FN} 
%P \bigl(  \bigl| \sum_{i = 1}^n  \bigl( |X_i + h| - u(h) \bigr)  \bigr|  \geq \sqrt{n} x \bigr)  \leq 
%2 \mathrm{exp}\bigl(- (1 + o(1)) \frac{x^2}{2} \bigr). 
%\end{equation} 
%At the same time, 
%\begin{equation} \label{mudelta} 
%u(h)  \sim \left\{ 
%\begin{array}{ll}
%\sqrt{2 / \pi},  &\qquad h \goto 0,  \\
%h,  &\qquad h \goto \infty,    
%\end{array} 
%\right.  
%\qquad 
%\sigma^2(h)  \sim 
%\left\{ 
%\begin{array}{ll}
%(\pi - 2)/ \pi,  &\qquad h \goto 0,  \\
%1,  &\qquad   h \goto \infty.    
%\end{array} 
%\right.
%\end{equation} 
%In  particular, there is a universal constant $c_ 0> 0$ so that  
%\begin{equation} \label{elementary2} 
%u(h) \geq   c_0
%\end{equation} 

We now show (\ref{FNadd}). Write for short  $Y_i = |z_i + h|$.  
It is sufficient to show that 
\begin{equation} \label{FN1a} 
P \bigl(   \sum_{i = 1}^n  Y_i  \geq  n u(h)    +  \sqrt{n} x \bigr)  \leq 
\mathrm{exp}\bigl(-  (1 + o(1))  \frac{x^2}{2 \sigma^2(h)}  \bigr),   
\end{equation} 
and 
\begin{equation} \label{FN1b} 
P \bigl(   \sum_{i = 1}^n  Y_i     \leq   n u(h)  -  \sqrt{n} x \bigr)  \leq 
\mathrm{exp}\bigl(-  (1 + o(1))  \frac{x^2}{2 \sigma^2(h)}  \bigr).    
\end{equation} 
Since the proofs are similar, we only show (\ref{FN1a}). 
By elementary calculations,  the moment generating function of $Y_i$  is 
\begin{equation}  \label{mgf} 
M_Y(s) = E[e^{sY}] = e^{s^2/2}[e^{h s} \Phi(s + h) + e^{- h s} \Phi(s - h)],
\end{equation} 
By Cramer-Chernoff Theorem (\cite{Lugosi}), for any $s > 0$ and any $y$,  
%%%%%%%%%%%
%%%%%%%%%%%
%%%%%%%%%%%
\begin{equation}\label{chernoff}
P(\sum_{i = 1}^n Y_i \geq ny) \leq e^{-n (ys - \log M_Y(s))}. 
\end{equation}
We now show this (\ref{FN1a}) for the cases of $h < 2\log(\sqrt{n}/x)$ and $h \geq  2\log(\sqrt{n}/x)$ 
separately. 

Consider the case where $h < 2\log(\sqrt{n}/x)$.  We wish to use (\ref{chernoff})  with 
\[
s = \frac{1}{\sigma^2(h)} \frac{x}{\sqrt{n}}, \qquad   y = u(h) + x/\sqrt{n}. 
\]  
By our assumptions of $h < 2 \log(\sqrt{n} / x)$ and $0 < x \leq \sqrt{n} /\log(n)$,  
\[
s  = O(x / \sqrt{n}) = o(1), \qquad  h s  \leq   2\log(\sqrt{n} / x)  (x / \sqrt{n})  = o(1). 
\]
Now, on one hand, since  $y \log^3(1/y) \goto 0$ as $y \goto 0+$,  $h^3 s  = o(1)$ and $h^3 s^3 = o(s^2)$.  Combining this with elementary Taylor expansion, 
\begin{equation} \label{Taylor1} 
e^{\pm h s} = 1 \pm h s  + \frac{(h s)^2}{2} + o(s^2). 
\end{equation} 
On the other hand, applying Taylor expansion to $\Phi(s \pm h)$ and noting that $\phi$ is a symmetric function,  
\begin{equation} \label{Taylor2} 
\Phi(s \pm h) = \Phi(\pm h) +  \phi(h) s - h \phi(h) s^2 + o(s^2).  
\end{equation} 
where we have used that the third derivative of $\Phi$ is a  bounded function.  Combining (\ref{Taylor1})-(\ref{Taylor2}) and re-arranging, 
\begin{align} 
& e^{h s} \Phi(s + h) + e^{- h s} \Phi(s - h)   \\ 
= & 1 + 2 s \phi(h) + h s[\Phi(h) - \Phi(-h)] + h^2 s^2/2 + o( s^2) \nonumber \\
= &  1 +  u(h) s + h^2 s^2/2  + o(s^2),  \label{Taylor3} 
\end{align} 
where in the first step,  we have used $\Phi(h) + \Phi(-h) = 1$, and in the second step, we have used  the expression of $u(h)$ given in (\ref{musigma}). 

We now analyze $\log[e^{h s} \Phi(s + h) + e^{- h s} \Phi(s - h)]$. Write for short $w = e^{h s} \Phi(s + h) + e^{- h s} \Phi(s - h) -1$. By (\ref{Taylor3}) and $|u(h)| \leq h + 1$ from (\ref{musigma}), $|w| \leq C\max\{(h+1) s, h^2 s^2\}$, and so  
\[
|\log(1 + w) - w  + w^2/2| \leq C |w|^3  \leq C  \max\{(h+1)^3 s^3, h^6 s^6\},   
\]
where by similar argument as above, $\max\{(h+1)^3 s^3, h^6 s^6 \}   = o(s^2)$.  Combining this with (\ref{Taylor3}), 
\begin{align*} 
& \log[e^{h s} \Phi(s + h) + e^{- h s} \Phi(s - h)]   \\
= &   \log(1 + w) \\
= & w  - w^2/2 + o(s^2) \\
= &  u(h) s + h^2 s^2/2 - [u(h)  s + h^2 s^2/2]^2/2 + o(s^2) \\
= &   u(h) s +  (h^2  - u(h)^2) s^2/2  -    [u(h) h^2 s^3 + h^4 s^4/4]/2  + o(s^2),  
\end{align*} 
where we note $|u(h) h^2 s^3 + h^4 s^4/4| \leq C (h+1) h^2 s^3 + h^4 s^4/4 = o(s^2)$. As a result, 
\[
\log[e^{h s} \Phi(s + h) + e^{- h s} \Phi(s - h)]   = u(h) s + (h^2 - u(h)^2) s^2/2 + o(s^2). 
\]
Combining this with (\ref{mgf}) and the expression of $\sigma(h)$ given in (\ref{musigma}) and rearranging it,  
\[ 
ys - \log[M_Y(s)]     =  (y - u(h)) s - (1 + h^2 - u(h)^2)\frac{s^2}{2} + o(s^2) = (y - u(h)) s - \sigma^2(h)\frac{s^2}{2} + o(s^2).
\] 
Now, invoking $s =  \frac{1}{\sigma^2(h)} x/\sqrt{n}$ and $y = u(h) + x /\sqrt{n}$ gives 
\[
ys - \log[M_Y(s)]  = \frac{1}{2\sigma^2(h)}(x/\sqrt{n})^2(1 + o(1)). 
\]
Combining this with (\ref{chernoff})  gives the claim.

We now consider the case of $h \geq 2\log(\sqrt{n}/x)$.   
We wish to use (\ref{chernoff}) again, with the same $y$ but a different $s$: $s = x/ \sqrt{n}$.  In the current case,  
since $x \leq \sqrt{n}/\log(n)$,   
\[ 
h  \goto \infty, \qquad s \goto 0. 
\] 
By the assumptions of $h \geq 2 \log(\sqrt{n} / x)$ and $s = x/\sqrt{n}$,  and 
\[ 
\phi(h/2) \leq C \mathrm{exp}(-(\log(\sqrt{n}/x))^2/2) = o(s^2), 
\] 
it follows that  $\max\{\Phi(- s - h), \Phi(s - h)\} \leq 
\Phi(-h/2)  = o(1) \phi(h/2)$, where the right hand side is $o(s^2)$. 
As a result, 
\begin{equation} \label{Taylor7} 
e^{h s} \Phi(s + h) + e^{- h s} \Phi(s - h)   = e^{h s} [1 - \Phi(-s-h)  + e^{-2h s} \Phi(s - h)]   
=  e^{h s}[1  + o(s^2)], 
\end{equation} 
and so 
\[ 
\log[e^{h s} \Phi(s + h) + e^{- h s} \Phi(s - h)]      =   h s +  o(s^2). 
\] 
Combining this with (\ref{mgf}) and (\ref{mudelta}) and invoking $s = x/\sqrt{n}$ and $y = u(h)  + x/\sqrt{n}$,  
%%%%%%%
%%%%%%%
%%%%%%%
\begin{align}
ys - \log M_Y(s) &= (u(h) +  x/ \sqrt{n}  - h) s - s^2/2 +  o(s^2)  \nonumber  \\ 
& =   s^2/2 +  o(s^2)  \nonumber   \\
& =  s^2/(2 \sigma^2(h)) + o(s^2),    \label{Taylor6} 
\end{align}
where in the last two steps, we have used  
\begin{equation}\label{mudelta}
h - u(h)  =  2 h \Phi(-h) - 2\phi(h)  = o(s),   \qquad \sigma^2(h) = 1 + h^2 - u(h)^2 = 1 + o(s). 
\end{equation} 
Inserting (\ref{Taylor6})  into (\ref{chernoff}) gives the claim.  \qed

%%%%%%%%%%%
\subsection{Proof of Lemma~\ref{lem:u-prop}}
Denote by  $\Phi$ the CDF of $N(0,1)$. By direct calculations,
\[
u(h)  =  \sqrt{2/\pi} e^{-h^2/2} + h[1 - 2 \Phi(-h)]. 
\] 
This implies $u(h)\goto \sqrt{2/\pi}$ when $h\goto 0$ and $u(h)/h\goto 1$ when $h\goto\infty$. Furthermore, 
\begin{align*}
& u'(h)  = -2 h\phi(h) + [1 - 2 \Phi(-h)] +  2h \phi(- h)  = 1 - 2 \Phi(-h),\cr
& u''(h) = 2 \phi(-h)  > 0.
\end{align*}
So $u(h)$ is strictly convex and monotony increasing for $h \in (0, \infty)$. 

Let $h_0$ be the unique solution of $u'(h)=0.9$. Fix $(h_1, h_2)$ such that $h_2>h_1>0$. If $h_1>h_0$, by convexity, 
\[
u(h_2)-u(h_1)\geq u'(h_1)(h_2-h_1)\geq 0.9(h_2-h_1).
\]
If $h_2<h_0$, using the Taylor expansion, for some $\tilde{h}\in[h_1, h_2]$, 
\[
u(h_2)-u(h_1)=u'(h_1)(h_2-h_1) + \frac{1}{2}u''(\tilde{h})(h_2-h_1)^2\geq \frac{1}{2}u''(h_0)(h_2-h_1)^2. 
\] 
If $h_1 < h_0 < h_2$, then we decompose the difference into $u(h_2) - u(h_0) + u(h_0) - u(h_1)$ and combine with the two cases we just dicussed, then we have that 
\[
u(h_2) - u(h_1) \geq 0.9(h_2 - h_0) + C_1 (h_0 - h_1)^2.
\]
When $h_2 - h_0 \geq h_0 - h_1$, then we have $u(h_2) - u(h_1) \geq 0.45(h_2 - h_0) + 0.45 (h_0 - h_1) = 0.45 (h_2 - h_1)$; otherwise, there is $u(h_2) - u(h_1) \geq \frac{C_1}{2}[(h_2 - h_0)^2 + (h_0 - h_1)^2] \geq C (h_2 - h_1)^2$. 
Combining the three cases gives the claim. 
\qed

%%%%%%%%%%%%%%%
%%%%%%%%%%%%%%%
\section{Proof of Secondary Lemmas}\label{app:C}
In this section, we show the proof of Lemmas~\ref{lem:chi-square}, \ref{lem:Frobenius} and \ref{lem:rmtcolor}.

\subsection{Proof of Lemma~\ref{lem:chi-square}} \label{proof:chi-square}
The following lemma is useful, which is proved below. 
%%%%%%%%%%%%
%%%%%%%%%%%%
\begin{lemma}  \label{lem:chi-square-0}
	For any fixed $q>0$, 
	\begin{align*} 
	& \pi_0^{(q)} = \bar{\Phi}\big( \sqrt{2q\log(p)} \big) \big(1+L_pn^{-1/2}\big),\cr
	& \pi_1^{(q)} =\left\{
	\begin{array}{ll}
	1 - L_p p^{-(\sqrt{r}-\sqrt{q})^2}, & r>q,\\
	\bar{\Phi}\big((\sqrt{q}-\sqrt{r})\sqrt{2\log(p)}\big) \big(1+L_pn^{-1/4}\big), & r\leq q. 
	\end{array}\right.
	\end{align*}
\end{lemma}
%%%%%%%%%%%%%
%\noindent
%We only show (a) and (b), and (c) follows from (a)-(b) and the definitions.  

First, we prove (a). Write for short $z_j=z$ and $z^{(q)}=z_j^{(q)}$. Since the distribution of $z^{(q)}$ is spherically symmetric, $v'z^{(q)}$ has the same distribution as $e_1'z^{(q)}$, for any $v\in \mathcal{S}^{n-1}$. It follows that $E[(v'z^{(q)})^2]=E[(z^{(q)}(1))^2]=a_p$. Furthermore, $E(\|z^{(q)}\|^2)=nE[(z^{(q)}(1))^2]=na_p$.  

Consider $E(|v'z^{(q)}|^m)$. Again, by spherical symmetry,
\[
E(|v'z^{(q)}|^m) = E\big(|z^{(q)}(1)|^m\big)= E\big(|z(1)|^m1\{z^2(1) + \|\tilde{z}\|^2>n+2\sqrt{qn\log(p)}\}\big),
\]
where $\tilde{z}=(z(2), \cdots, z(n))'$. Note that $\tilde{z}$ is independent of $z(1)$ and $\|\tilde{z}\|^2\sim \chi^2_{n-1}$. Let $B_1$ be the event that $|z(1)|\leq \sqrt{2\delta_1\log(p)}$, for some $\delta_1$ to determine. From basic properties of the $N(0,1)$ distribution,  $P(B_1^c)=L_pp^{-\delta_1}$ and $E(|z(1)|^mI_{B_1^c})=L_pp^{-\delta_1}$. It follows that
\begin{align*}
E(|v'z^{(q)}|^m) &\leq E\big(|z(1)|^m1\{z^2(1) + \|\tilde{z}\|^2>n+2\sqrt{qn\log(p)}, B_1\}\big) + L_pp^{-\delta_1}\cr
&\leq E\big(|z(1)|^m1\{\|\tilde{z}\|^2>n+2\sqrt{qn\log(p)}-2\sqrt{\delta_1\log(p)}\}\big) + L_pp^{-\delta_1}\cr
& = E(|z(1)|^m)\cdot P(\|\tilde{z}\|^2>n+2\sqrt{qn\log(p)}(1+o(1))) +L_pp^{-\delta_1}\cr
& = E(|z(1)|^m) \pi_0(1+o(1)) + L_pp^{-\delta_1}. 
\end{align*}
By choosing $\delta_1$ appropriately large, we find that the first term dominates. 

Consider $E(\|z^{(q)}\|^{2m})$. %The claim for $a_p$ is then a special case $m=1$. 
Denote by $f_{n}$ the density of $\chi^2_{n}$, where $f_n(y) = \frac{y^{n/2-1}e^{-y}}{2^{n/2}\Gamma(n/2)}$. Note that $y^mf_n(y)=\frac{2^m\Gamma(m+n/2)}{\Gamma(n/2)}f_{n+2m}(y)$. It follows that 
\[
E(\|z^{(q)}\|^{2m}) 
%& = \frac{2^m\Gamma(m+n/2)}{\Gamma(n/2)}  \int_{n+2\sqrt{qn\log(p)}}^{\infty}f_{n+2m}(y)dy\cr
= \frac{2^m\Gamma(m+n/2)}{\Gamma(n/2)} P(\chi^2_{n+2m}>n+2\sqrt{qn\log(p)}).
%= \kappa_m(n)\pi_0(1+o(1)), 
\]
First, by letting $q=0$ on both hand sides, we have $\kappa_{2m}(n)=E(\|z\|^{2m})=\frac{2^m\Gamma(m+n/2)}{\Gamma(n/2)}$. Second, 
since $n+2\sqrt{qn\log(p)}=n_*+2\sqrt{qn_*\log(p)}(1+L_pn^{-1/2})$ for $n_*=n+2m$, Lemma~\ref{lem:chi-square-0} implies that $P(\chi^2_{n+2m}>n+2\sqrt{qn\log(p)})=\pi_0(1+o(1))$.  Together, the above right hand side is 
$\kappa_{2m}(n)\pi_0(1+o(1))$. 

Consider $a_p$. Similarly to the above, for $n_*=n+2$,
\begin{align*}
a_p=n^{-1}E(\|z^{(q)}\|^{2}) &= \frac{2\Gamma(1+n/2)}{n\Gamma(n/2)}P(\chi^2_{n+2}>n_*+2\sqrt{qn_*\log(p)}(1+L_pn^{-1/2}))\cr
& = P(\chi^2_{n+2}>n_*+2\sqrt{qn_*\log(p)}(1+L_pn^{-1/2}))\cr
& = \pi_0(1+L_pn^{-1/2}). 
\end{align*}

Second, we prove (b). We first state an approximation of $\pi_1$. From basic properties of chi-square distributions, for all $q,r\geq 0$, 
\begin{align*}
&P(\chi^2_n(0)>n+2\sqrt{qn\log(p)})=\bar{\Phi}(\sqrt{2q\log(p)})(1+L_pn^{-1/2}),\cr
&P(\chi^2_n(2r\log(p))>n + 2\sqrt{qn\log(p)})=\bar{\Phi}\big((\sqrt{q}-\sqrt{r})\sqrt{2\log(p)}\big)(1+L_pn^{-1/4}). 
\end{align*}
Therefore, we find that
\begin{align}  \label{lem-chisquare-1}
\pi_1& = P(\chi^2_n(2r\log(p))>n+2\sqrt{qn\log(p)})\cr
& = P(\chi^2_n(0)>2(\sqrt{q}-\sqrt{r})\sqrt{n\log(p)})\cdot(1+o(1)). 
\end{align}

Consider $E[(v'z_j^{(q)})^2]$. Fix $v$ and introduce
\[
w_1 = \ell/\|\ell\|, \qquad  w_2 = (1-(v'\ell)^2/\|\ell\|^2)^{-1/2}[ v - (v'\ell)\ell/\|\ell\|^2 ].
\]
Both $w_1$ and $w_2$ are unit vectors and $w_1'w_2=0$. Let $Q$ be any orthogonal matrix whose first two columns are $w_1$ and $w_2$. By direct calculations, $Q'v = (x_0, \sqrt{1-x_0^2},0,\cdots,0)'$ and $Q'\ell = (\sqrt{n}, 0, \cdots, 0)$, where $x_0=(v'\ell)/\|\ell\|$. Since $Q'z$ and $z$ have the same distribution, 
\begin{align}  \label{lem-chisquare-2}
v'z_j^{(q)} & = v'QQ'z\cdot 1\{\|Q'z + \mu(j) Q'\ell\|^2 > n+2\sqrt{qn\log(p)} \}\big]\cr
& \overset{(d)}{=} v'Qz\cdot 1\{\|z+\mu(j) Q'\ell \|^2 > n+2\sqrt{qn\log(p)}\}\big]\cr
& = \bigl[ x_0z(1)+(1-x_0^2)^{1/2}z(2)\bigr] \cdot1\{ \|z+\sqrt{n}\tau_p^* e_1\|> n \}. 
\end{align}
It follows that 
\begin{align*}
E[ (v'z_j^{(q)})^2 ]
&= E\big[ (x_0 z(1) + (1-x_0^2)^{1/2} z(2))^2 1\{\|z+\sqrt{n}\tau_p^*e_1\|^2>n+2\sqrt{qn\log(p)} \}    \big]\cr
&=  (1-x_0^2) E\big[  (z(2))^2 1\{\|z+\sqrt{n}\tau_p^*e_1\|^2>n+2\sqrt{qn\log(p)} \}  \big]\cr
& + x_0^2\ E\big[ (z(1))^2 1\{\|z+\sqrt{n}\tau_p^*e_1\|^2> n+2\sqrt{qn\log(p)}\}\big] \cr
&= b_p +  (c_p - b_p) (v'\ell)^2/\|\ell\|^2,
\end{align*}
where the second equality comes from the symmetry on $z(2)$ (so the cross term disappears). 
% and the last equality is because $y_0^2=(\mu(j)\|\ell\|)^2=\sqrt{2n}\tau_p(r)$. 

Consider $b_p$ and $c_p$. Let $\tilde{z}=(z(2),\cdots,z(n))'$, where $\|\tilde{z}\|^2\sim\chi^2_{n-1}$ and it is independent of $z(1)$. We write
\[
c_p= E\big[(z(1))^2 1\{\|\tilde{z}\|^2 > n+2\sqrt{qn\log(p)}-g(z(1)) \}\big], \quad g(x)\equiv (x+ \sqrt{n}\tau_p^*)^2. 
\]
For a constant $\delta_2>0$ to be determined, let $B_2$ be the event that $|z(1)|\leq \sqrt{2\delta_2\log(p)}$. From basic properties of normal distributions, $P(B_2^c)=L_p p^{-\delta_2}$ and $E[z^2(1) I_{B_2^c}]=L_pp^{-\delta_2}$. Over the event $B_2$, we have $g(z(1))= [z(1)- (2\sqrt{nr\log(p)})^{1/2}]^2=2\sqrt{rn\log(p)}(1+L_pn^{-1/4})$. It follows that  
\begin{align*}
c_p &\leq E\big[(z(1))^2 \cdot P(B_2\cap \{\|\tilde{z}\|^2> n + 2\sqrt{qn\log(p)} -g(z(1))\} | z(1)) \big] + L_pp^{-\delta_2}\cr
&\leq  E[(z(1))^2]\cdot P\big( \chi^2_{n-1}> n + 2(\sqrt{q}-\sqrt{r})\sqrt{n\log(p)}(1+L_pn^{-1/4}) \big) + L_pp^{-\delta_2} \cr
&= \pi_1(1+L_p n^{-1/4}), 
\end{align*}
where the last inequality comes from \eqref{lem-chisquare-1} and that $\delta_2$ is chosen appropriately large.  To compute $b_p$, we write 
\begin{align*}
b_p & = E\big[ (z(2))^2 1\{ \|\tilde{z}\|^2 > n+2\sqrt{qn\log(p)} -g(z(1))\} \big]\cr
& =  (n-1)^{-1} E\big[ \|\tilde{z}\|^2 1\{ \|\tilde{z}\|^2 > n+2\sqrt{qn\log(p)} -g(z(1))\}  \big]. 
\end{align*}
Let $B_2$ be the same event. Let $q_*=[(\sqrt{q}-\sqrt{r})_+]^2$. We have
\begin{align*}
b_p 
&=  (n-1)^{-1} E\big[ \|\tilde{z}\|^2 1\{ \|\tilde{z}\|^2> n + 2\sqrt{q_*n\log(p)}(1+L_pn^{-1/4}) \}\big]  + L_pp^{-\delta_2} \cr
&= (n-1)^{-1}E(\|\tilde{z}^{(q_*)}\|^2)(1+L_pn^{-1/4}) = \pi_1(1+L_pn^{-1/4}), 
\end{align*}
where in the last equality, we have applied the result in (a) with $q=q_*$. 

Consider $E(|v'z^{(q)}|^m)$. Let $\tilde{w}=(z(3),\cdots,z(n))'$. Then $\|\tilde{w}\|^2\sim\chi^2_{n-2}$ and it is independent of $(z(1),z(2))$. By \eqref{lem-chisquare-2}, 
\begin{align*}
E(|v'z^{(q)}|^m) & = E\bigl(|x_0z(1)+(1-x_0^2)^{1/2} z(2)|^m\cr
&  \cdot1\{\|\tilde{w}\|^2>n+2\sqrt{qn\log(p)}-g(z(1))-(z(2))^2\}\bigr). 
\end{align*}
Let $B_3$ be the event that $\max\{|z(1)|, |z(2)|\}\leq \sqrt{2\delta_3\log(p)}$. Then $P(B_3^c)=L_pp^{-\delta_3}$ and over $B_3$, $g(z(1))+(z(2))^2=2\sqrt{rn\log(p)}(1+o(1))$. Applying similar arguments as above, we find that
\[
E(|v'z^{(q)}|^m) \leq E(|x_0z(1)+ (1-x_0^2)^{1/2} z(2)|^m)\cdot \pi_1(1+o(1)) =  \kappa_m\pi_1(1+o(1)). 
\]
Here the last inequality is because $x_0z(1)+ (1-x_0)^{1/2} z(2)\sim N(0, 1)$. The claim then follows. 

Consider $E(\|z^{(q)}\|^{2})$ and $E(\|z^{(q)}\|^{2m})$. Using $Q$ defined above (for an arbitrary $v$) 
\begin{align*}
E(\|z^{(q)}\|^{2}) &= E\big( \|Q'z\|^2 1\{\|Q'z+\tau_p^*Q'\ell\|^2>n+2\sqrt{qn\log(p)} \} \big)\cr
&= E\big( \|z\|^2 1\{\|z+\sqrt{n}\tau_p^*e_1 \|^2>n+2\sqrt{qn\log(p)} \} \big)\cr
&= E\big( (z(1))^2 1\{\|z+\sqrt{n}\tau_p^*e_1 \|^2>n+2\sqrt{qn\log(p)} \} \big) \cr
& + (n-1) E\big( (z(2))^2 1\{\|z+\sqrt{n}\tau_p^*e_1 \|^2>n+2\sqrt{qn\log(p)} \} \big)\cr
& = c_p + (n-1)b_p. 
\end{align*}
Recall that $\tilde{z}=(z(2),\cdots, z(n))'$, $q_*=[(\sqrt{q}-\sqrt{r})_+]^2$ and $g(x)=(x+\sqrt{n}\tau_p^*)^2$ for any $x\in R$. Note that $(x+y)^m\leq 2^m(|x|^m+|y|^m)$ for any $x,y\in R$. We have 
\begin{eqnarray} \label{z2m}
E(\|z^{(q)}\|^{2m}) 
&=& E\big( \|z\|^{2m} 1\{\|z+\sqrt{n}\tau_p^*e_1 \|^2>n+2\sqrt{qn\log(p)} \} \big)\nonumber\\
&\leq & 2^m E\big( (z(1))^{2m} 1\{\|\tilde{z}\|^2 >n+2\sqrt{qn\log(p)} -g(z(1))\} \big)\nonumber\\
&& + 2^m E\big( \|\tilde{z}\|^{2m} 1\{\|\tilde{z}\|^2 >n+2\sqrt{qn\log(p)} -g(z(1))\} \big)\nonumber\\
&=& 2^m \kappa_{2m} \pi_1(1+o(1)) + 2^m E\big( \|\tilde{z}^{(q_*)}\|^{2m} \big)(1+o(1))\nonumber\\
&=& 2^m (\kappa_{2m} + \kappa_{2m}(n-1))\cdot \pi_1(1+o(1)).
\end{eqnarray}
Here, we have applied the result in (a) for $E(\|z^{(q)}\|^{2m})$ with $q=q_*$. 

%Last, we prove (c). Note that $z_j^{(q)}$ and $z_k^{(q)}$ are independent. Moreover, $z^{(q)}_j/\|z_j^{(q)}\|$ and $\|z_j^{(q)}\|$ are independent. As a result,
%\begin{align*}
%E\big\{ [(z_j^{(q)})'z_k^{(q)}]^{2m}\big\} & = E\big( \|z_j^{(q)}\|^{2m}  E[ (v'z_k^{(q)})^2m]\big|_{v=z_j^{(q)}/\|z_j^{(q)}\|}. 
%\end{align*}
%The claims then follow immediately from (a) and (b). \qed

Last, we prove (c). Using the spherical symmetry of $z^{(q)}$ and the $Q$ defined above,  we have already seen that $\|z+\tau_p^*\ell\|\overset{(d)}{=} \|z+\sqrt{n}\tau_p^*e_1\|$ and 
\[
z\cdot 1\{\|z+\tau_p^*\ell\|>n+2\sqrt{qn\log(p)} \}\overset{(d)}{=}z\cdot 1\{\|z+\sqrt{n}\tau_p^*e_1\|>n+2\sqrt{qn\log(p)} \}. 
\]
Then the claims follow from the definitions and (a)-(b). \qed

%%%%%%%%%%%
%%%%%%%%%%%
\subsection{Proof of Lemma~\ref{lem:Frobenius}}
%Note that the assumption $n^{-7/4}\|\mu\|_0\pi_p^{(1)}/[N_p^*(\mu)]^{1/2}= L_pn^{-\delta}$ translates to
%\beq  \label{thm-postZ-10}
%L_pn^{-5/4}\|\mu\|_0\pi_p^{(1)} = o(1)\cdot \sqrt{nN_p^*(\mu)}, \qquad \delta_p= \sqrt{nN_p^*(\mu)}(1+o(1)). 
%\eeq
Let $\pi(j)=\pi_0$ for $j\notin S(\mu)$ and $\pi(j)=\pi_1$ for $j\in S(\mu)$, where $\pi_0,\pi_1$ are defined in Section~\ref{subsec:prelim}. Then $\sum_{j=1}^p\pi(j)=m^{(q)}$ by (c) of Lemma~\ref{lem:chi-square}. 

First, consider $\mathrm{tr}(H_0)$. Write $M_j=n^{-1}[\|z_j^{(q)}\|^2 - E(\|z_j^{(q)}\|^2)]$. By definition, 
\beq \label{lem-Frobenius-1}
n^{-1}\mathrm{tr}(H_0)-m_*^{(q)} = \sum_{j=1}^p M_j. 
\eeq
By Lemma~\ref{lem:chi-square}, $E(\|z_j^{(q)}\|^2)\lesssim n\pi(j)$ and $E(\|z_j^{(q)}\|^{2m})\leq 2^{m}\kappa_{2m}(n)\pi(j)\leq C4^m\pi(j)n^m$, where $\kappa_{2m}(n)$ is the $m$-th moment of the $\chi^2_n$ distribution and we have used $\kappa_{2m}(n)\leq C2^mn^m$. Noting that $(a+b)^m\leq 2^m(a^m+b^m)$ for any real values $a$ and $b$, by direct calculations, 
\[
E(M_j)=0, \qquad  \mathrm{var}(M_j) \leq C\pi(j), \qquad E(|M_j|^m)\leq C8^m. 
\]
By Lemma~\ref{lem:Bernstein} (Bernstein inequality), with probability at least $1-O(p^{-3})$,
\beq \label{lem-Frobenius-2}
|\sum_{j=1}^p M_j|\leq C\sqrt{m^{(q)}\log(p)}. 
\eeq
Combining \eqref{lem-Frobenius-1}-\eqref{lem-Frobenius-2} gives the first claim. 

Second, consider $\|H_0\|_F^2$. 
By direct calculations, 
\begin{eqnarray}  \label{thm-postZ-7}
&& \|H_0\|_F^2 - n^{-1}[\mathrm{tr}(H_0)]^2= \sum_{1\leq j,k\leq p} [(z_j^{(q)})'z_k^{(q)}]^2 - n^{-1}(\sum_{j=1}^p \|z_j^{(q)}\|^2)^2\\
&=& \frac{n-1}{n}\sum_{j=1}^p \|z_j^{(q)}\|^4 + 2\sum_{1\leq j<k\leq p}\left( [(z_j^{(q)})'z_k^{(q)}]^2 -\frac{1}{n} \|z_j^{(q)}\|^2 \|z_k^{(q)}\|^2\right) \nonumber\\
& \equiv &  (I) + (II). \nonumber
\end{eqnarray}
We now study $(I)$. Write $U_j=n^{-1}\|z_j^{(q)}\|^4$ for short. By Lemma~\ref{lem:chi-square}, $E(U_j)=n^{-1}\kappa_4(n)\pi_0(1+o(1))$ and $\mathrm{var}(U_{j})\leq Cn^2\pi_0$ for $j\notin S(\mu)$; moreover, $\mathrm{var}(U_{j})\leq Cn^2 \pi_1$ for $j\in S(\mu)$. We also claim that $E(U_j)\geq n^{-1}\kappa_4(n-1)\pi_1(1+o(1))$ for $j\in S(\mu)$. The proof is similar to that for \eqref{z2m}, but in the second line of \eqref{z2m}, we instead use the inequality $\|z\|^{2m}\geq \|\tilde{z}\|^{2m}$. Note that $\kappa_4(n)$ is the second moment of $\chi^2_n$ and so $\kappa_4(n)=n^2+2n$. 
It follows that
\[
\sum_{j=1}^p E(U_j)\gtrsim nm^{(q)}, \qquad \sum_{j=1}^p \mathrm{var}(U_j) \leq Cn^2m^{(q)}. 
\]
Using Lemma~\ref{lem:Bernstein}, with probability at least $1-O(p^{-3})$, 
$\sum_{j=1}^p U_j\gtrsim nm^{(q)} - Cn\sqrt{m^{(q)}\log(p)}\gtrsim Cnm^{(q)}$. Since $(I)= (n-1) \sum_{j=1}^p U_j$, 
\beq  \label{thm-postZ-8}
(I) \geq  C_1n^2 m^{(q)}, \qquad \mbox{for some constant }C_1>0. 
\eeq
We then study $(II)$. Let $V_{jk}=[(z_j^{(q)})'z_k^{(q)}]^2 -\frac{1}{n} \|z_j^{(q)}\|^2 \|z_k^{(q)}\|^2$. Introduce 
\[
W_j(v) = [v'z_j^{(q)}]^2 - n^{-1}\|z_j^{(q)}\|^2, \qquad \mbox{for any }v\in \mathcal{S}^{n-1}. 
\]
Let $v_j=z_j/\|z_j\|$. Then $v_j$ is independent of $\|z_j^{(q)}\|$ and $V_{jk}=\|z_j^{(q)}\|^2 W_k(v_j)$.  
By Lemma~\ref{lem:chi-square}, for any fixed $v\in \mathcal{S}^{n-1}$, 
\[
\begin{array}{lr}
E[W_j(v)]=0, \;\; E[(W_j(v))^2] \leq C\pi_0, & j\notin S(\mu),\\
E[W_j(v)]=(c_p-b_p)[(v'\tilde{\ell})^2 - n^{-1}], \;\; E[(W_j(v))^2]\leq C\pi_1, & j\in S(\mu),  
\end{array}
\]
where $\tilde{\ell}=\ell/\|\ell\|$. 
%For $1\leq j\leq p$, let $\pi(j)=\pi_1$ for $j\in S(\mu)$ and $\pi(j)=\pi_0$ for $j\notin S(\mu)$. 
As a result, if either $j\notin S(\mu)$ or $k\notin S(\mu)$, then $E(V_{jk})=0$; if both $j,k\in S(\mu)$, then $E(V_{jk})=(c_p-b_p)E\{\|z_j^{(q)}\|^2[(v_j'\tilde{\ell})^2 - n^{-1}]\}=(c_p-b_p)E[W_j(\tilde{\ell})]=(1-n^{-1})(c_p-b_p)^2\geq 0$. It follows that
\beq \label{thm-postZ-9(1)}
E[(II)]\geq 0. 
\eeq
To compute $\mathrm{var}((II))$, we calculate $E(V_{jk}V_{j'k'})$ for all $(j,k,j',k')$ such that $j\neq k$ and $j'\neq k'$. Since $V_{jk}=V_{kj}$, we assume $j\neq k'$ and $j'\neq k$ without loss of generality.  
%For notation convenience, write $\pi(j)=\pi_0$ for $j\notin S(\mu)$ and $\pi(j)=\pi_1$ for $j\in S(\mu)$. 
We have the following observations: (1) $E(V_{jk})\leq (c_p-b_p)^2$ if both $j,k\in S(\mu)$ and $E(V_{jk})=0$ otherwise. (2) $E(V_{jk}^2)=E(\|z_j^{(q)}\|^2)E[W^2_k(v_j)]\leq Cn^2\pi(j)\pi(k)$ for any $j\neq k$.  
(3) When $j\neq j'$ and $k\neq k'$, $V_{jk}$ is independent of $V_{j'k'}$, so $E(V_{jk}V_{j'k'})=E(V_{jk})E(V_{j'k'})$. 
(4) When $j=j'$ and $k\neq k'$, $E(V_{jk}V_{jk'})=E[\|z_j^{(q)}\|^4 W_k(v_j)W_{k'}(v_j)]$; as a result, $E(V_{jk}V_{jk'})=0$ when either $k\notin S(\mu)$ or $k'\notin S(\mu)$; if $k,k'\in S(\mu)$, $E(V_{jk}V_{jk'})=(c_p-b_p)^2 E[(W_j(\tilde{\ell}))^2]\leq C(c_p-b_p)^2\pi(j)$. Therefore, 
\begin{align*}
\sum_{\substack{(j,j',k,k'):\\ j\neq k,j'\neq k'}}& E(V_{jk}V_{j'k'})  
\leq \sum_{(j,k):j\neq k} Cn^2 \pi(j)\pi(k) + \sum_{\substack{(j,k,k'):j \notin \{k,k'\} \\ \{k,k'\}\subset S(\mu), k\neq k'}} C(c_p-b_p)^2 \pi(j)\cr
& + \sum_{\substack{(j,j',k,k'): \{j,j',k,k'\}\subset S(\mu)\\j,j',k,k'\text{ are different}}} (c_p-b_p)^4\cr
& \leq Cn^2(m^{(q)})^2 + C(c_p-b_p)^2 m^{(q)}|S(\mu)|^2 + C(c_p-b_p)^4|S(\mu)|^4\cr
&\leq Cn^2(m^{(q)})^2 +  m^{(q)} (|S(\mu)|\pi_1)^2 \cdot o(1)+ C(|S(\mu)|\pi_1)^4\cdot o(1), 
\end{align*}
where the last inequality is due to that $c_p-b_p=o(\pi_1)$. Using \eqref{rate}, when $r<\rho^*_\theta(\beta)$ (``impossibility") and $q<\tilde{q}(\beta,\theta,r)$ (``fat" case), $(|S(\mu)|\pi_1)^2=o(m^{(q)})$ and so the first term in the above dominates the other two. It implies  
\beq \label{thm-postZ-9(2)}
\mathrm{var}((II)) \leq C_2 n^2(m^{(q)})^2, \qquad \mbox{for some constant }C_2>0.
\eeq
We combine \eqref{thm-postZ-9(1)}-\eqref{thm-postZ-9(2)} and apply the Markov inequality. It follows that with probability at least $1-4n^{-2}C_2/C_1^2$, 
\beq \label{thm-postZ-9}
(II)\geq -C_1n^2m^{(q)}/2. 
\eeq
The second claim follows by plugging \eqref{thm-postZ-8} and \eqref{thm-postZ-9} into \eqref{thm-postZ-7}. \qed

%%%%%%%%%%%
\subsection{Proof of Lemma~\ref{lem:rmtcolor}}
% Write
%$\delta_p = L_p \sqrt{np}$. The claim translates to
%\beq  \label{thm-postZ-3-col}
%\|ZBB'Z' - \mathrm{trace}(BB')\|\leq \delta_p. 
%\eeq
The proof is similar to that of Lemma~\ref{lem:postZ-1}. 
By Lemma~\ref{lem:alpha-net}, there exists an $(1/4)$-net of $\mathcal{S}^{n-1}$, 
denoted as $\mathcal{M}^*_{1/4}$, such that $|\mathcal{M}^*_{1/4}|\leq 9^n$ and $\sup_{v\in \mathcal{M}^*_{1/4}} v'Av\geq 2\|A\|$ for any $n\times n$ matrix $A$. Therefore, 
to show the claim, it suffices to show that for each fixed $v\in \mathcal{M}^*_{1/4}$, with probability $\geq 1-O(9^{-n}p^{-2})$, 
\beq  \label{thm-postZ-4-col}
|v'(ZBB'Z' - \mathrm{tr}(BB') I_n)v| \leq  C\sqrt{np}.  
\eeq 

Denote the eigenvalue decomposition of $BB'$ by $V' \Lambda V$, where $\Lambda$ is diagonal matrix with diagonals $\lambda_1 \geq \lambda_2 \geq \cdots \geq \lambda_p$.  
Fix $v$, we can write
\[
v'ZBB'Z'v = v'Z V' \Lambda V Z' v = \sum_{i = 1}^{p} \lambda_i \eta_i^2, \qquad \eta_i \stackrel{iid}{\sim} N(0,1).
\]
The last equation comes from $V Z' v \sim N(0, I_p)$. 
So we have $E[v' ZBB'Z' v] = \mathrm{tr}(BB')$ for any fixed $v$ with $\|v\| = 1$. 
Let $W_j = \lambda_i \eta_i^2/\lambda_1 - \lambda_i/\lambda_1$, then $W_j$'s are independent of each other, $E(W_j) = 0$, $\mathrm{var}(W_j) \leq 2 $ and $E(|W_j|^m)\leq \kappa_{2m}$.
We apply Lemma~\ref{lem:Bernstein} with $\lambda=2\sqrt{n\log(9)+2\log(p)}$. To check the moment conditions, we note that $\kappa_{2m}=E_{z\sim N(0,1)}(|z|^{2m})\leq 2^m m!$ for all $m\geq 1$. 
It follows that with probability $\geq 1-O(9^{-n}p^{-2})$,
\[
\lambda_1 |\sum_{j=1}^p W_j|\leq 2\sqrt{np\log(9)+2p\log(p)} \lambda_1 \leq C\sqrt{np}. 
\]
The last inequality is because $\lambda_1 = \|B\|^2 \leq L_p$. 
This proves \eqref{thm-postZ-4-col}.\qed
% then \eqref{thm-postZ-3-col} follows immediately. 

%%%%%%%%%%%%%%%
%%%%%%%%%%%%%%%
\subsection{Proof of Lemma~\ref{lem:chi-square-0}}

We start from computing $\pi_0$.
Using the density of the $\chi^2_n$ distribution, 
\[
\pi_0 = \int_{n+2\sqrt{qn\log(p)}}^{\infty} \frac{ x^{n/2-1}e^{-x/2}dx}{2^{n/2}\Gamma(n/2)} \equiv \frac{1}{2^{n/2}\Gamma(n/2)} \cdot (I). 
\]
Now, we calculate the integral $(I)$.  Write for short
\[
t = \sqrt{2q\log(p)} \qquad\mbox{and}\qquad  x_0=n+\sqrt{2n}t.
\]
With a variable change $x=n+\sqrt{2n}y$, we have
\begin{align}  \label{moment-equ5}
(I) &= \sqrt{2n}  x_0^{n/2-1} e^{-x_0/2}     \int_0^{\infty} (1 + \sqrt{2n}  y /x_0)^{n/2-1} e^{- \sqrt{2n}  y/2} dy\cr
& = \sqrt{2n}  x_0^{n/2-1} e^{-x_0/2}     \int_0^{\infty} \exp\Big\{ (n/2-1)\log(1 + \sqrt{2n}y /x_0)- \sqrt{2n}  y/2\Big\} dy\cr
& \equiv \sqrt{2n}  x_0^{n/2-1} e^{-x_0/2}  \big[ (I_1) + (I_2) + (I_3)\big],
\end{align}
where $(I_1)$ contains the integral from $0$ to $ct$, $(I_2)$ contains that from $ct$ to $x_0/\sqrt{2n}$ and $(I_3)$ contains that from $x_0/\sqrt{2n}$ to infinity. We will determine the constant $c>0$ later.

Consider $(I_1)$. From the Taylor expansion, $\log(1+a)=a-a^2/2+O(a^3)$ for small $a$. Moreover, $n/x_0= 1-\sqrt{2n}t/x_0+O(t^2/n)$, $x_0=O(n)$ and $y=O(t)$. As a result, for $0<y<ct$, by simple calculations,  
\[
(n/2-1) \log(1 + \sqrt{2n} y /x_0)  - \sqrt{2 n}y / 2= -t y -y^2/2 + O(t^3/\sqrt{n}). 
\]
Noting that $e^{a}=1+O(a)$ for small $a$, so $(I_1)$ is equal to $\int_{0}^{ct}e^{-  t y   - y^2/2}dy\cdot [ 1+ O(t^3/\sqrt{n})]$. By direct calculation, $\int_{0}^{ct}e^{-  t y   - y^2/2}dy=e^{t^2/2}\int_{t}^{(1+c)t}e^{-y^2/2}dy = \sqrt{2\pi} e^{t^2/2}\big[\bar{\Phi}(t)-\bar{\Phi}((1+c)t)\big]$. 
By Mills' ratio, $\bar{\Phi}(t)=L_pp^{-q}$ and $\bar{\Phi}((1+c)t)=L_pp^{-(1+c)^2q}$. Therefore, when $c$ is chosen large enough, 
$\bar{\Phi}((1+c)t)=o(1)\cdot \bar{\Phi}(t) L_pn^{-1/2}$. It follows that
\beq  \label{moment-equ2}
(I_1) = \sqrt{2\pi} e^{t^2/2}\bar{\Phi}(t)\big( 1+ L_p/\sqrt{n}\big).
\eeq
Consider $(I_2)$. Since $\log(1+a)-a\leq a^2/4$ for $a\in [0,1]$, when $ct<y<x_0/\sqrt{2n}$,
\[
(n/2-1) \log(1 + \sqrt{2n} y /x_0)  - \sqrt{2 n}y / 2 
= -t y -y^2/4 + O(tx_0^2/(\sqrt{n})^3)
\]
As a result,
$(I_2) \leq (1 + L_p n^{-1/2}) \int_{ct}^{\infty}e^{-yt-y^2/4}dy$, where 
$\int_{ct}^{\infty}e^{-yt-y^2/4}dy= 2\sqrt{\pi}e^{t^2}\bar{\Phi}((c+2)t/\sqrt{2})= e^{t^2/2}L_p \sqrt{n}p^{-[(c+2)^2/2-1]q}$. 
By choosing $c$ appropriately large, we have
\beq \label{moment-equ3}
(I_2) = o(1)\cdot e^{t^2/2}L_pn^{-1/2}.
\eeq
Consider $(I_3)$. Since $\log(1+a)\leq t$ for all $a\geq 0$, when $t>x_0/\sqrt{2n}$,
\[
(n/2-1) \log(1 + \sqrt{2n} y /x_0)  - \sqrt{2 n}y / 2 \leq  - (n / x_0)  t y  \leq -t y/2 .
\]
It follows that $(I_3) \leq \int_{x_0/\sqrt{2n}}^{\infty} e^{-ty/2}dy = (2/t) e^{-x_0^2/(2n)} = o(1)\cdot e^{t^2/2}L_pn^{-1/2}$. 

Combining the above results for $(I_1)$-$(I_3)$, we obtain that 
\[
\pi_0 =R_n(t) \cdot \bar{\Phi}(t)\big(1+L_pn^{-1/2}\big), \qquad \mbox{where } R_n(t) \equiv \frac{2\sqrt{\pi n}  x_0^{n/2-1} e^{-x_0/2+t^2/2}}{2^{n/2}\Gamma(n/2)}.
\]
We plug in $x_0=n+\sqrt{2n}t$ and rewrite $R_n(t) = \frac{n}{x_0}\frac{\sqrt{\pi} (n/e)^{n/2}}{\Gamma(n/2)}\big( 1+ t\sqrt{2/n}\big)^{n/2} e^{-t\sqrt{n/2} +t^2/2}$. Note that by Taylor expansion, $\log(1+a)=a-a^2/2+O(a^3)$ for $a=t\sqrt{2/n}$. Therefore, we have  $R_n(t) = \frac{n}{x_0} \frac{\sqrt{\pi} (n/e)^{n/2}}{\Gamma(n/2)} \exp\big\{ O(t^3/\sqrt{n}) \big\}= 1+L_pn^{-1/2}$. This gives 
\[
\pi_0 = \bar{\Phi}(t)\big(1+L_pn^{-1/2}\big). 
\]

Next, we compute $\pi_1$. Define 
\[
\tilde{r} = \frac{(z(1)-\sqrt{n}\tau_p)^2}{2\sqrt{n\log(p)}}, \qquad W = \sum_{i=2}^n z^2(i). 
\]
Then $\tilde{r}$ and $W$ are independent; furthermore, $W$ has a $\chi^2_{n-1}$ distribution. We rewrite
\[
\pi_1 = E\left[ P\left( \frac{ W- n}{\sqrt{2n}} > (\sqrt{q} - \sqrt{\tilde{r}}) \sqrt{2\log(p)} \middle| \tilde{r}\right) \right]. 
\]
For a constant $c>0$ to be determined, let $B_1$ be the event that $|z(1)|\leq \sqrt{2c\log(p)}$. Then $P(B_1^c)=L_pp^{-c}$. Over the event $B_1$, 
$\tilde{r} =  r + L_p n^{-1/4}$. 
When $r>q$, utilizing the results for $\pi_0$, we get 
\begin{align*}
\pi_1 & = \Phi\big( (\sqrt{r}-\sqrt{q}+o(1))\sqrt{2\log(p)} \big)(1+L_pn^{-1/2}) + L_p p^{-c} \cr
&= 1 - L_p p^{-(\sqrt{r}-\sqrt{q})^2} + L_p p^{-c}. 
\end{align*}
When $r<q$. 
\begin{align*}
\pi_1 & = \bar{\Phi}\big( (\sqrt{q}-\sqrt{r}+L_pn^{-1/4})\sqrt{2\log(p)} \big)(1+L_pn^{-1/2}) + L_p p^{-c} \cr
&=  \bar{\Phi}\big( (\sqrt{q}-\sqrt{r})\sqrt{2\log(p)} \big)(1+L_pn^{-1/4}) + L_p p^{-c}. 
\end{align*}
We choose $c$ large enough so that $L_p p^{-c}$ is always dominated by any other term. 
This gives the claim for $\pi_1$. \qed

\bibliographystyle{imsart-number}
\bibliography{ifpca}
%%%%%%%%%%%%
%%%%%%%%%%%%

\end{document}